\documentclass[graybox]{svmult}

\usepackage{type1cm}        % activate if the above 3 fonts are
\usepackage{makeidx}         % allows index generation
\usepackage{graphicx}        % standard LaTeX graphics tool
\usepackage{multicol}        % used for the two-column index
\usepackage[bottom]{footmisc}% places footnotes at page bottom

\usepackage{newtxtext}       % 
\usepackage[varvw]{newtxmath}       % selects Times Roman as basic font

\spnewtheorem{thm}{Theorem}[section]{\bf}{\it}
\spnewtheorem{cor}[thm]{Corollary}{\bf}{\it}
\spnewtheorem{prop}[thm]{Proposition}{\bf}{\it}
\spnewtheorem{lem}[thm]{Lemma}{\bf}{\it}
\spnewtheorem{ex}[thm]{Example}{\bf}{\rm}
\spnewtheorem{rem}[thm]{Remark}{\bf}{\rm}
\spnewtheorem{OP}[thm]{Problem}{\bf}{\rm}

\newcommand{\C}{\ensuremath{\mathbb{C}}}
\renewcommand{\D}{\ensuremath{\mathbb{D}}}
\newcommand{\G}{\ensuremath{\Gamma}}
\newcommand{\mc}{\ensuremath{\mathcal}}
\newcommand{\mf}{\ensuremath{\mathfrak}}
\newcommand{\N}{\ensuremath{\mathbb{N}}}
\newcommand{\ov}{\ensuremath{\overline}}
\newcommand{\oline}{\ensuremath{\overline}}
\newcommand{\R}{\ensuremath{\mathbb{R}}}
\renewcommand{\Re}{\ensuremath{\mathrm{Re}}}
\renewcommand{\Im}{\ensuremath{\mathrm{Im}}}

\DeclareMathOperator{\Ad}{Ad}
\DeclareMathOperator{\ad}{ad}

\DeclareMathOperator{\Kern}{ker}
\DeclareMathOperator{\Span}{span}
\DeclareMathOperator{\grad}{grad}
\DeclareMathOperator{\Id}{Id}
\DeclareMathOperator{\RT}{RT}
\DeclareMathOperator{\im}{im}

\usepackage{enumerate}

\renewcommand{\a}{\mathfrak a}
\newcommand{\Abb}[4]{\begin{cases} #1 & \rightarrow  #2 \\ #3 &\mapsto  #4\end{cases}}
\newcommand{\bbar}{\overline}
\newcommand{\g}{\mathfrak g}

\newcommand{\M}{\mathcal{M}}
\newcommand{\n}{\mathfrak n}
\newcommand{\tu}[1]{\textup{#1}}

\newcommand{\X}{\mathcal{X}}

\newcommand{\demi}{\tfrac{1}{2}}

\DeclareMathOperator{\Diffeo}{Diffeo}
\DeclareMathOperator{\ran}{ran}
\DeclareMathOperator{\supp}{supp}

\newcommand{\cfcn}{\mathbf{c}} % HC c function

\newcommand{\la}{\langle}
\newcommand{\ra}{\rangle}
\newcommand{\1}{\mathbf{1}}

\newcommand{\psdiff}{pseudo-differential}
\DeclareMathOperator{\Op}{Op}
\newcommand{\Dprime}{\ensuremath{\mathcal{D}^\prime}}
\newcommand{\intd}{\mspace{0.5mu}\operatorname{d}\mspace{-2.5mu}}

\newcommand{\Ccinfty}{\ensuremath{C_c^{\infty}}}
\newcommand{\bigoh}{\mathcal{O}}

\newcommand{\bq}{\begin{equation}}
\newcommand{\eq}{\end{equation}}
\newcommand{\bqn}{\begin{equation*}}
\newcommand{\eqn}{\end{equation*}}
\newcommand{\Cinft}{{\rm C^{\infty}}}

\newcommand{\on}[1]{\ensuremath{\operatorname{#1}}}
\newcommand{\rhoa}{\ensuremath{\rho}}

\newcommand{\pst}[1]{\ensuremath{H_{#1}}}

\newcommand{\Ga}{\Gamma}
\newcommand{\bs}{\backslash}

\makeindex             % used for the subject index
\begin{document}

\title*{Quantum-Classical Correspondences for Locally Symmetric Spaces}
\author{Joachim Hilgert}
\institute{Joachim Hilgert \at Institut für Mathematik, Universität Paderborn, D-33095 Paderborn, Germany,\\
\email{hilgert@upb.de}}
%\and Name of Second Author \at Name, Address of Institute %\email{name@email.address}}

%\footnote{Funded by the Deutsche Forschungsgemeinschaft (DFG, German Research Foundation) -- SFB-TRR 358/1 2023 -- 491392403}
\maketitle

\abstract{We review the status of a program, outlined and motivated in the introduction, for the study of correspondences between spectral invariants of partially hyperbolic flows on locally symmetric spaces and their quantizations. Further, we briefly sketch an extension of the program to graphs and quotients of higher rank affine buildings. Finally, we formulate a number of concrete problems which may be viewed as possible further steps to be taken in order to complete the program.}

\begin{acknowledgement} 
While this article does not explicitly refer to the work of Toshiyuki Kobayashi, it still owes a lot to him. His broad way of looking at representation theory with an open mind enabling him to merge algebraic, analytic and geometric methods such as derived functor modules, microlocal analysis and complex geometry, has always been a shining example for me.

As for specific mathematical insights pertaining to the subject matter of this article I thank my  coauthors for patiently sharing their knowledge with me and enabling me to attack problems with tools I could not have handled by myself. 

This work was partially funded by the Deutsche Forschungsgemeinschaft (DFG, German Research Foundation) -- SFB-TRR 358/1 \emph{Integral Structures in Geometry and Representation Theory} 2023 -- 491392403
\end{acknowledgement}

%\tableofcontents

\section{Introduction}
\label{sec:Introduction}
There are very different mathematical models for the dynamics of physical systems depending on whether one is looking at classical or quantum mechanics. Nevertheless, relations between mathematical objects on both sides can be found. In particular, if the quantum model contains the Planck constant $\hbar$ explicitly, one studies the behavior of the model for $\hbar$ tending to zero in order to understand to which extent the classical model is a limit of the quantum model. There is no such canonical procedure to describe a quantum version of a classical physical system -- if there were such a canonical quantization, then quantum mechanics would be unnecessary in the first place. 

A paradigmatic class of examples comes from Riemannian manifolds $\X$. Their cotangent bundles $T^*\X$ carry canonical symplectic forms so that smooth functions $f$ on $T^*\X$ define Hamiltonian vector fields $X_f$ which can be integrated to give flows. Interpreting $f$ as the energy observable on the state space $T^*\X$, classical Hamiltonian mechanics says that the flow of $X_f$ is the time evolution of the system. On the other hand, one has linear differential operators $D$ with smooth coefficients on $\X$ which define unbounded operators on the space $L^2(\X,\mathrm{vol})$ of functions square integrable with respect to the volume form coming with the Riemannian metric. If such a $D$ is (essentially) selfadjoint, it can be interpreted as the energy operator of a quantum system whose time evolution is given by the unitary one-parameter group $e^{itD}$ with infinitesimal generator $iD$. Associated with $D$ is its symbol $\sigma_D$ which is a function on $T^*\X$. 

It is a legitimate question to ask whether one can read off properties of a quantum system given by $D$ from properties of the classical system given by $\sigma_D$. This question received particular attention in the context of classically chaotic dynamical systems when it was observed in computer experiments that quantizations of chaotic systems have different spectral properties than quantizations of completely integrable classical systems. Mathematically this is difficult to pin down as quantum systems are linear by definition and linear systems do not feature chaotic behavior in the established sense. Trying to give rigorous explanations for the observed phenomena is the subject of the field of \emph{quantum chaos}.

A class of examples studied very closely in quantum chaos are hyperbolic surfaces together with their Laplace-Beltrami operators and the geodesic flows generated by their symbols. These surfaces have negative curvature which tends to produce chaotic behavior. Moreover, they occur as moduli spaces of complex tori which are of great relevance in number theory, whence they have been studied very intensely. So there is a lot of information available. In particular, it is a well-established fact that hyperbolic surfaces are a meeting point of many mathematical fields such as  dynamical systems, complex and harmonic analysis, number theory, algebraic and differential geometry. The wealth of documented knowledge and available mathematical tools made hyperbolic surfaces an excellent toy model to study  quantum chaos.  

Studying examples like hyperbolic surfaces, one realizes that the presence of symmetry restricts the possibilities to quantize a classical system. Actually, one can take the success of the orbit method in representation theory as an argument supporting this view. We take it as a guideline to look for correspondences between dynamical invariants of quantum systems with many symmetries and their classical counterparts. It turns out that this leads to very interesting results even far beyond hyperbolic surfaces.  

The major part of this article will deal with such quantum-classical correspondences for systems built from locally symmetric spaces. The symmetry groups in question are (mostly semisimple) real Lie groups. It turns out that one can find non-Archimedean counterparts given in terms of affine buildings which may serve as toy models because of their combinatorial nature. First steps in establishing these toy models will also be described.

\subsection{Some motivating examples}

The first examples of quantum-classical correspondences were typically not designed as such, but need to be reinterpreted to actually show this nature. 

\begin{ex}[Poisson summation formula] The classical Poisson summation formula 
\[\sum_{n\in \mathbb Z^\ell}f(x+n)=\sum_{m\in \mathbb Z^\ell} \hat f(m) e^{2\pi imx},\] 
say for Schwartz functions $f$ on $\mathbb{R}^\ell$ and their Fourier transforms $\hat f$, can be interpreted as follows: The $n\in \mathbb Z^\ell$ parameterize the closed geodesics in the $\ell$-torus $\mathbb{T}^\ell:=\mathbb{R}^\ell/\mathbb{Z}^\ell$ passing through the origin via $n\mapsto \mathbb{R}n+\mathbb{Z}^\ell$. On the other hand, the $m\in \mathbb Z^\ell$ parameterize the eigenvectors $x\mapsto e^{2\pi i m x}$ of the Laplace operator $\Delta=\sum_{j=1}^\ell \frac{\partial^2}{\partial x_j^2}$ in $L^2(\mathbb{T}^\ell)$. Thus the summation formula relates the spectrum of the quantum object $\Delta$ with the geometric ``spectrum'' of oriented closed geodesics (up to translation). The flexibility of this correspondence lies in the fact that one can choose the test functions $f$ freely.  
\end{ex}

Also the next example has its origin in harmonic analysis and is of importance in particular in number theory.

\begin{ex}[Selberg trace formula] In his paper \cite{Se56} Selberg established a formula for compact locally symmetric spaces (with noncompact semisimple symmetry group) which, similar to the torus case, has a geometric side that can be interpreted as a sum over classes of closed geodesics and a spectral side that involves the spectrum of invariant operators. This paper started a spectacular development which has been described in a number of excellent survey articles and books. As trace formulae are not in the focus of this article we refrain from giving any details but refer to \cite{Wa76,Gu77} as good starting points to read up and \cite{DKV79} for the interpretation in terms of geodesics. 
\end{ex}

We come to an example, which so far has not been interpreted as a quantum-classical correspondence, but does motivate such a correspondence (see Example~\ref{ex:Lewis-Zagier}).

%Eichler-Shimura isomorphism

\begin{ex}[Modular cusp forms and period polynomials] We refer to  {\cite[Chap.~2]{MR20}} for details on the material in this example. In particular we follow the notation and conventions laid out in the ``Introductory Roadmap'' of  \cite{MR20} here as well as in the following Example~\ref{ex:Lewis-Zagier}.
A \emph{modular cusp form of weight} $k\in 2\mathbb N$ is a holomorphic function $F$  on the upper half-plane $\mathbb H:=\{z\in \mathbb C\mid \Im z>0 \}$ such that 
\[\forall a,b,c,d\in \mathbb Z, ad-bc=1, z\in \mathbb H:\quad F\Big(\frac{az+b}{cz+d} \Big)=(cz+d)^k F(z),\]
$F(z)\longrightarrow 0$ for $z$ tending to $i\infty$ in a strip of bounded width, and $F$ admits a Fourier expansion of the form
\[F(z)=\sum_{n=1}^\infty a(n) e^{2\pi inz}.\]
We denote the space of modular cusp forms of weight $k$ by $S_k$.

The \emph{period polynomial} $P$ of a modular cusp form $F$ of weight $k$ is given by the integral
\[P(X):=\int_0^{i\infty} (z-X)^{k-2} F(z)\,\mathrm{d}z.
\]
with $X\in \mathbb C$. It satisfies the identities
\begin{eqnarray}\label{eq:period polynomials}
\label{eq:period polynomials1}
P+ (P\vert_{2-k} T)&=&0,\\
\label{eq:period polynomials2}
P+ (P\vert_{2-k} TS)+(P\vert_{2-k} (TS)^2) &=&0,
\end{eqnarray}
where 
\[T:=\begin{pmatrix}
0&-1\\ 1&0
\end{pmatrix}
\quad\text{and}\quad
S:=\begin{pmatrix}
1&1\\ 0&1
\end{pmatrix}
\]
and 
\[\Big(f\Big\vert_m \begin{pmatrix}
a&b\\ c&d
\end{pmatrix}\Big)(z):= (cz+d)^{m} f\Big(\frac{az+b}{cz+d}\Big)
\]
for  $\begin{pmatrix}
a&b\\ c&d
\end{pmatrix}\in \mathrm{SL}_2(\mathbb Z)$, $m\in \mathbb Z$ and $f$ a function on $\mathbb H$. The set of polynomials occuring as period polynomials of modular cusp forms of weight $k$ is the space of polynomials of degree at most $k-2$ satisfying the equations \eqref{eq:period polynomials1} and \eqref{eq:period polynomials2}. We denote it by $\mathbb P_{k-2}$. 

As a corollary to the Eichler Cohomology Theorem one finds that 
\[ S_k\oplus S_k \cong \mathbb P_{k-2}/\mathbb C[X^{k-2}-1].\] 
More precisely, \cite[Cor.~2.53]{MR20} gives an isomorphism between $S_k\oplus S_k$ and the \emph{parabolic cohomology} $H^1_\mathrm{par}(\mathrm{SL}_2(\mathbb Z),\mathbb C[X]_{k-2})$, where $\mathbb C[X]_{k-2}$ denotes the space of polynomials of degree less or equal to $k-2$. The cocycles defining \emph{Eichler cohomology classes} in $H^1(\mathrm{SL}_2(\mathbb Z),\mathbb C[X]_{k-2})$ are $\mathbb C[X]_{k-2}$-valued functions $M\mapsto P_M$ on $\mathrm{SL}_2(\mathbb Z)$ satisfying the \emph{cocycle identity}
\[\forall M,N\in \mathrm{SL}_2(\mathbb Z):\quad  P_{MN} = (P_M\big\vert_{2-k} N )+ P_N.\]
The \emph{coboundaries} are the functions $P_M$ satisfying
\[\exists Q\in \mathbb C[X]_{k-2}\  \forall M\in \mathrm{SL}_2(\mathbb Z):\quad P_M=  (Q\big\vert_{2-k}  M)-Q.\]
A cocycle is called \emph{parabolic} if
\[\exists Q\in \mathbb C[X]_{k-2}:\quad P_S=  (Q\big\vert_{2-k}  S)-Q.\]
The parabolic Eichler cohomology $H^1_\mathrm{par}(\mathrm{SL}_2(\mathbb Z),\mathbb C[X]_{k-2})$ consists of the cohomology classes represented by parabolic cocycles. Then one has a short exact sequence 
\[0
\longrightarrow \mathbb C[X^{k-2}-1] 
\longrightarrow \mathbb P_{k-2}
\underset{\tiny P\mapsto(P_M)}{\longrightarrow} H^1_\mathrm{par}(\mathrm{SL}_2(\mathbb Z),\mathbb C[X]_{k-2})
\longrightarrow 0,
\]
where the second arrow is determined by $P_T=P$ and $P_S=0$.
\label{ex:period-polynomials}
\end{ex}

\begin{ex}[Lewis-Zagier correspondence] 
We keep the notation from Example~\ref{ex:period-polynomials} and  refer to  \cite[Chap.~3-5]{MR20} for details on the material in this example. A \emph{Maass cusp form} is an $\mathrm{SL}_2(\mathbb Z)$-invariant function $F$  on the upper half-plane $\mathbb H$ such that there exists a $\lambda\in \mathbb C$ with $\Delta F=\lambda F$, where $\Delta:=-y^2\Big(\frac{\partial^2}{\partial x^2}+\frac{\partial^2}{\partial y^2}\Big)$ is the \emph{Laplace-Beltrami operator} for the upper half-plane, and $F(z)\longrightarrow 0$ for $z$ tending to $i\infty$ in a strip of bounded width. Each Maass cusp form is real analytic and square integrable as a function on the \emph{modular surface} $\mathrm{SL}_2(\mathbb Z)\backslash \mathbb H$. The $\Delta$-eigenvalues are of the form $\lambda=s(1-s)$ with $\Re s=\frac{1}{2}$. The number $s\in \mathbb C$ is sometimes called the \emph{spectral parameter} of corresponding eigenfunction. Lewis and Zagier proved in \cite{LZ01}  that there is a linear bijection between the space of Maass cusp forms with spectral parameter $s$ and the space of real analytic solutions of the \emph{three-term functional equation}
\[\psi(x)=\psi(x+1)+(x+1)^{-2s}\psi\Big(\frac{x}{x+1}\Big)\]
on $\mathbb R_{> 0}$ which satisfy the growth conditions
\[\psi(x)=\begin{cases}
o(x^{-1})&x\longrightarrow 0, \\
o(1)&x\longrightarrow \infty.
\end{cases}
\]
The functions $\psi$ were called \emph{period functions} by Lewis and Zagier as they have various properties analogous to the properties of period polynomials (see \cite[\S\,2]{LZ01}). For example, the functional equation can be rewritten as
\[\psi=\psi\big\vert_{-2s}S+\psi\big\vert_{-2s}STS\]
and for an even Maass cusp form $F$ according to \cite[(2.2)]{LZ01} the associated period function has the integral representation
\[\psi(x)=\int_0^{\infty}\frac{xy^s}{(x^2+y^2)^{s+1}} F(iy)\, \mathrm{d}t.\]
In fact, in the follow-up paper \cite{BLZ15} to \cite{LZ01} Bruggeman, Lewis and Zagier also established an intepretation of period functions in terms of parabolic cohomology.

While we do not try at this point to relate parabolic cohomology to the geodesic flow on the modular surface (that such a relation exists is suggested for instance by the work of Juhl \cite{Ju01} and Bunke-Olbrich \cite{BO99}), we note  that \cite[\S\,3]{LZ01} contains an interpretation of period functions as eigenfunctions of transfer operators which are built from a symbolic dynamics associated with the geodesic flow on the modular surface. This gives the Lewis-Zagier correspondence an interpretation as a quantum-classical correspondence, see also the survey article \cite{PZ20}. 
\label{ex:Lewis-Zagier} 

In view of this interpretation of the Lewis-Zagier correpondence a natural question to ask is, whether one can give a reasonable quantum system having modular cusp forms in the sense of Example~\ref{ex:period-polynomials} as states and such that it quantizes a classical system to which one can naturally associate period polynomials. 
\end{ex}

The Lewis-Zagier correspondence from Example~\ref{ex:Lewis-Zagier} has been set up specifically for the modular surface, but meanwhile period functions have been generalized in various ways (see \cite{MR20} for a fairly up-to-date compilation). Also the transfer operator approach has been extended substantially, mostly through the work of Pohl and her collaborators (see e.g. \cite{Po12}).

\begin{ex}[Distributions invariant under the horocycle flow]
Let $\X$ be a hyperbolic surface of finite volume and $S\X$ its tangent unit sphere bundle. In \cite{FF03} Flaminio and Forni studied the space $\mathcal J(S\X)$ of compactly supported distributions on $S\X$ which are invariant under the \emph{horocycle flow} generated by the matrix 
\[U:=\begin{pmatrix}
0&1\\ 0&0
\end{pmatrix}\in \mathfrak{sl}_2(\mathbb R),
\]
where $S\X$ is realized as a quotient $\Gamma\backslash S\mathbb H\cong \Gamma\backslash \mathrm{PSL}_2(\mathbb R)$ for a discrete subgroup $\Gamma$ of $\mathrm{PSL}_2(\mathbb R)$ acting without fixed points on $\mathbb H$. They showed that $\mathcal J(S\X)$ decomposes as
\begin{equation}\label{eq:U-invariant distributions}
\mathcal J(S\X)=\bigoplus_{\mu\in \mathrm{spec}_\mathrm{pp}(\Delta_\X)}\mathcal J_\mu \oplus \bigoplus_{n\in \mathbb N} \mathcal J_n^d \oplus \bigoplus_{c\in \mathcal C} \mathcal J_c,
\end{equation}
where $\mathrm{spec}_\mathrm{pp}(\Delta_\X)$ is the pure point spectrum of the Laplace-Beltrami operator $\Delta_\X$ of $\X$ and $\mathcal C$ is the (finite) set of cusps of $\X$. The spaces $\mathcal J_n^d$ are closely related to the spaces of holomorphic sections of powers of the canonical line bundle over $\X$. This decomposition is derived from the Plancherel decomposition of $L^2(S\X)\cong L^2(\Gamma\backslash \mathrm{PSL}_2(\mathbb R))$. 

For $0<\mu<\frac{1}{4}$ the  spaces $\mathcal J_\mu$ decompose further as $\mathcal J_\mu=\mathcal J_\mu^+\oplus \mathcal J_\mu^-$  with $\dim(\mathcal J_\mu^\pm)$ equal to the multiplicity of the eigenvalue $\mu$. For $\mu\ge\frac{1}{4}$ the multiplicity of the eigenvalue $\mu$ is equal to  $\dim(\mathcal J_\mu)$. Finally, the space $\mathcal J_0$ is spanned by the  $\mathrm{PSL}_2(\mathbb R)$-invariant volume.  

In the case of compact surfaces, $\mathrm{spec}_\mathrm{pp}(\Delta_\X)$ is the entire spectrum of $\Delta_\X$ and $\mathcal C$ is empty. In the noncompact case, $\mathrm{spec}(\Delta_\X)$ has an absolutely continuous part consisting of $[\frac{1}{4},\infty[$ with multiplicity equal to the number of cusps.
Thus in both cases one can extract quantum information, namely the spectral data of the Laplace-Beltrami operator, from  $\mathcal J(S\X)$. 

%\footnote{Discuss the quantum interpretation of the exceptional points belong to the Casimir spectrum $\mu=-n^2-n$ (discrete series)}   

Flaminio and Forni also observe that the splitting  \eqref{eq:U-invariant distributions} is invariant under the geodesic flow. More precisely, they show that it is essentially the spectral decomposition of the geodesic flow (see \cite[Thm.~1.4]{FF03} for details). Thus, $\mathcal J(S\X)$ yields a spectral quantum-classical correspondence for finite area hyperbolic surfaces. 
 \label{ex:Flaminio-Forni}
\end{ex}

\begin{ex}[Dyatlov-Faure-Guillarmou correspondence]
Let\label{ex:DFG} ~$\X$ be a compact hyperbolic surface and let $X$ be 
the vector field generating the geodesic flow $\varphi_t$ on the unit tangent 
bundle $S\X$ of $\X$. The linear operator
\[
 \mathcal L_t: \left\{ \begin{array}{ccc}
                       C_c^\infty(S\X) &\to&C_c^\infty(S\X)\\
                       f&\mapsto &f \circ \varphi_{-t}
                      \end{array}
\right.
\]
is called the \emph{transfer operator} of the geodesic flow and the vector field
$-X$ is its generator.

For $f_1,f_2\in C_c^\infty(S\X)$ one defines the \emph{correlation functions}
\[ 
C_X(t;f_1,f_2):=\int_{S\X}(\mathcal L_t f_1)\cdot  f_2\intd\mu_L
\]
where $\mu_L$ is the Liouville measure (invariant by $\varphi_t$). By 
\cite{BuLi,FaSj,DyZw} the Laplace transform 
\[
R_X(\lambda;f_1,f_2):=-\int_0^\infty e^{-\lambda t}C_X(t; f_1,f_2)dt
\]
extends meromorphically from ${\Re}(\lambda)>0$ to $\C$. 
For 
${\Re}(\lambda)>0$ we have $R_X(\lambda;f_1,f_2)=\langle(-X-\lambda)^{-1}f_1,f_2\rangle$ and 
$R_X(\lambda)$ gives a meromorphic extension of the Schwartz kernel of the resolvent of $-X$. 
The poles are called \emph{Ruelle resonances} and 
the \emph{residue operator} $\Pi^X_{\lambda_0}:C_c^\infty(S\X)\to \mc{D}'(S\X)$
defined by 
\[ 
\forall f_1,f_2\in C^\infty_c(S\X):\quad \langle\Pi^X_{\lambda_0}f_1,f_2\rangle:=\mathrm{Res}_{\lambda_0}R_X(\lambda;f_1,f_2)
\] 
has finite rank,  commutes with $X$, and $(-X-\lambda)$ is nilpotent on its range. 
The elements in the range of $\Pi_{\lambda_0}^X$ are called \emph{generalized Ruelle 
resonant states}. By results in \cite{BuLi, FaSj, DyZw} the poles
can be identified with the discrete spectrum of $-X$ in certain
Hilbert spaces and the generalized resonant states with generalized eigenfunctions. 

In \cite{DFG15} Dyatlov, Faure and Guillarmou show that if $\mathrm{spec}(\Delta_\X)=\{s_j(1-s_j)\}$ for a sequence  $\{s_j\in[0,1]\cup(\frac{1}{2}+ i\R)\mid j\in \N_0\}$, then the Ruelle resonances of $S\X$ in $\C\setminus (-1-\frac{1}{2}\N_0)$ are given by
\[\{\lambda_{j,m}:= -m-1+s_j\mid j,m\in \N_0\}.\]
The sets $\{\lambda_{j,m}\mid j\in \N_0\}$ are called the \emph{bands} of the Ruelle resonance spectrum. One refers to $\{\lambda_{j,0}\mid j\in \N_0\}$ as the \emph{first band}. It is characterized by the fact that the corresponding resonant states are annihilated by the horocycle vector field generated by  $U\in \mathfrak{sl}_2(\R)$, see Example~\ref{ex:Flaminio-Forni}.
Moreover, the points $\lambda\in -1-\frac{1}{2}\N_0$ excluded from this correspondence between the spectrum of the Laplacian and the Ruelle resonance spectrum of the geodesic flow are called \emph{exceptional points}. 

In fact, \cite{DFG15} contains generalizations of the above results to compact real hyperbolic spaces of arbitrary dimension.
\end{ex}

\subsection{Outline of the program and its purpose}\label{subsubsec: program}

Let $\X=\Gamma\backslash G/K$ be a locally symmetric space and $\mathcal{A}$ the algebra of $G$-invariant smooth functions on the cotangent bundle $T^*(G/K)$. The elements of $\mathcal{A}$ may be considered as Poisson commuting symbols of $G$-invariant pseudo-differential operators on $G/K$ but, when $\Gamma$-invariant, also as functions on $T^*(\X)$. The goal is now to study commutative algebras $\Psi$ of pseudo-differential operators on $\X$ quantizing $\mathcal{A}$ and to relate dynamical invariants of $\mathcal A$ (classical side) and $\Psi$ (quantum side). 

On the classical side the primary dynamical systems to be considered are the Hamiltonian flows associated with the elements of $\mathcal A$. Derived from those one has restrictions of flow-invariant submanifolds such as $\Gamma$-quotients of $G$-orbits or subbundles given e.g. by hyperbolic properties of the Hamiltonian flows. This may lead to restrictions in the flow variables as well (e.g. ``positive time'' only). For these derived dynamical systems, alternative descriptions may be available, for instance via symbolic dynamics. Given any of these dynamical systems, one may consider  dynamical invariants such as compact orbits, dynamical zeta functions and their divisors, transfer (or Koopman) operators and their spectral theory or classical scattering. As in the higher rank situation one typically has multiparameter dynamical actions, one wants to establish multiparameter versions of these invariants rather than simply using the invariants coming with one-parameter restrictions of the dynamical systems. 

On the quantum side one expects the pseudo-differential operators to generate multiparameter flows of linear operators, whose spectral theory needs to be studied. $G$-invariance allows to decompose the dynamical system and leads to questions of representation theory. Here spectral invariants can involve parameters of representations such as infinitesimal characters, dimensions of associated spaces (e.g. consisting of $\Gamma$-invariant distributions) or spectral invariants of the flow operators themselves. The latter could be ordinary spectra, resonances or derived objects such as divisors of meromorphically continued scattering operators. 

The general strategy to relate the quantum and the classical side is to use the $G$-invariance properties to connect the spectral invariants with $G$-representations and then use representation theoretic knowledge about intertwining of representations to relate two such invariants. A particularly successful example in this direction has been to relate the classical dynamics with principal series representations and then use Poisson transform to connect them to eigenspace representations which in turn are intimately related to the quantum dynamics.

Apart from supporting the philosophy that the presence of symmetry narrows down the possibilities to quantize given classical systems, the main motivation for establishing quantum-classical correspondences is to create bridges which allow to solve problems on one side using the tool boxes from the other side. An example of this kind of application is the proof in \cite{HWW21} of uniform spectral gaps of Weyl chamber flows in higher rank. 

Finally, one should note that given specific results for the rigid symmetric situations, one can by deformation sometimes guess or even prove results for situations where the symmetry is broken.   This is the case for instance in \cite{KW20}, where Betti numbers are given as the dimensions of Ruelle resonant states. Another example are the cases of quantum-classical correspondences for finite graphs in \cite{BHW23,AFH23a} which take a similar form as the ones for compact rank one locally symmetric spaces but need no homogeneity.

\section{Locally symmetric spaces}\label{sec:loc sym spaces}

Locally symmetric spaces are quotients of symmetric spaces by discrete subgroups of the group of isometries. In general they have an orbifold structure but we will always make assumptions on the discrete subgroups which guarantee that the quotient is actually a manifold. A natural way to study locally symmetric spaces is to lift objects to the symmetric space, where one can use group invariance as a tool and check results for invariance under the discrete subgroup allowing to interprete them as results for the locally symmetric space. Thus a detailed description of the symmetric spaces and (geometric) objects derived from them will be useful.

\subsection{Symmetric spaces}\label{sec: Classical Dynamics}

Consider a real semi-simple Lie group $G$, connected and of noncompact type, and let $G= K A N$ be an Iwasawa decomposition with $A$ abelian, $K$ a compact maximal subgroup and $N$ nilpotent. Then $ A\cong \R^\ell$ and $\ell$ is called the \emph{real rank} of $G$. Let $\a$ be the Lie algebra of $A$ and consider the adjoint action of $\a$ on $\mathfrak g$ which leads to the definition of a finite set $\Sigma\subset \a^*$ of \emph{restricted roots}. More precisely,
for $0\not=\alpha\in \mf a^*$ set
$$\g_\alpha:=\{X\in \g\mid \forall D\in \mf a:\ [D,X]=\alpha(D)X\}.$$
If $\g_\alpha\not=\{0\}$, then  $\alpha$ is a \emph{restricted root} of $(\mf a,\g)$.  An element $X$ of $\mf a$ is \emph{regular} if $\alpha(X)\not=0$ for all $\alpha\in \Sigma$. For each regular element $X_o$ one can define a set of {\it positive roots} via
$$\Sigma^+:=\{\alpha\in \Sigma\mid \alpha(X_o)>0\}.$$
Then $\Sigma$ is the disjoint union of $\Sigma^+$ and $\Sigma^-:=-\Sigma^+$, the set of \emph{negative roots}. We fix some choice of positive roots and define the corresponding  \emph{positive Weyl chamber} by
\[
\a_+:= \{H\in\a\mid \forall  \alpha\in\Sigma_+: \alpha(H)>0 \}.
\]
Note that  $\mf a$ is a {\it Cartan subspace}, i.e. a  maximal abelian subspace of $\mf p:=\mf k^\perp\subseteq \mf g$, where the orthogonal complement is taken with respect to the Killing form. Further let $M':=N_K(\mf a)$ be the normalizer
of $\mf a$ in $K$ and $M:=Z_K(\mf a)$ the centralizer
of $\mf a$ in $K$.
Then  $Z_K(\mf a)$ is normal in $N_K(\mf a)$. The
quotient group  $W:= W(\mf a,K):=N_K(\mf a)/Z_K(\mf a)$
is the  {\it Weyl group}.

For later use we set
$$\mf k_\alpha:=\mf k\cap(\g_\alpha+\g_{-\alpha}),\quad \mf p_\alpha:=\mf p\cap(\g_\alpha+\g_{-\alpha})$$
for $\alpha\in \mf a^*$.

\begin{rem}[Quotients by discrete torsion free subgroups]~If one now considers  a torsion free discrete subgroup $\Gamma < {G}$ one can define the locally symmetric space  $\X:=\Gamma\backslash G/K$ and the manifold  $\M := \Gamma \backslash  G /  M$. If $G$ is of real rank one, then $\M$ can be identified with the sphere bundle $S\X$. In higher rank the situation is more complicated as we will explain below.
\label{rem:Gamma-quotients}
\end{rem}

\subsubsection{Equivariant differential geometry on $G/K$}

We review some preparatory material on differential analysis of the natural $G$-action on $G/K$ which can found in \cite{Hi05}. The action of $G$ on $G/K$ induces natural actions of $G$ on the
tensor bundles of $G/K$. Viewing $g\in G$ as a diffeomorphism on
$G/K$, the action on the tangent bundle $T(G/K)$ is given by the
derivative $g'\colon T(G/K)\to T(G/K)$.
Then the action on the cotangent bundle $T^*(G/K)$ is defined via the canonical 
pairing
\[\forall
   \xi\in T^*_x(G/K), v\in T_x(G/K):\quad \langle g\cdot \xi, v\rangle_{G/K}=\langle \xi, g^{-1}\cdot v\rangle_{G/K}=
   \langle \xi, (g^{-1})'(x)(v)\rangle_{G/K}.\]
%The action on the higher tensor bundles $T^r_s\mathcal{N}$ is then defined via tensor products.
Note that in this way all the canonical projections are $G$-equivariant
and $G$ acts on all bundles by bundle maps covering the original
action, so we obtain induced actions on the spaces of sections. For example
$G$ acts on the space ${\cal V}(G/K)$ of smooth vector fields via
\begin{equation}\label{vfactiondef}
(g\cdot X)(x)=g\cdot (X(g^{-1}\cdot x)).
\end{equation}
It is now clear what is meant by a $G$-invariant vector field, metric,
symplectic form, and so on.

\begin{rem}[Homogeneous vector bundles and invariant vector fields] 
Let $H$ be a closed subgroup of $G$. Then $T(G/H)\cong G\times_H(\g/\mf h)$, where $G\times_H(\g/\mf h)$ is the homogeneous vector bundle over $G/H$ with respect to the $H$-action on $\mf g/\mf h$. More precisely, $G\times_H(\g/\mf h)$ is the set of $H$-orbits on $G\times(\g/\mf h)$ under the right  $H$-action 
\[(g,X+\mf h)\cdot h=(gh,h^{-1}\cdot X+\mf h),\]
where $h\cdot X$ denotes the adjoint action. We denote the equivalence class of $(g,X+\mf h)$ by ${[g,X+\mf h]}$. Then $g_1\cdot[g_2,X+\mf h]:=[g_1g_2,X+\mf h]$ defines a left action of $G$ on $G\times_H (\g/\mf h)$.  The derivative of the left translation
\[\lambda_g\colon G/H\to G/H,\quad    g_1H\mapsto gg_1H\]
is an isomorphism  denoted by  $\lambda_g'(x)\colon T_x(G/H)\to T_{g\cdot x}(G/H)$ and we identify $G\times_H (\g/\mf h)$ with $T(G/H)$ via the equivariant correspondence
\begin{equation}\label{invhom1}
         [g,X+\mf h]\longleftrightarrow \lambda'_g(o)(X+\mf h),
\end{equation}
where $o=H\in G/H$ is the canonical base point of $G/H$.

The dimensions in the stratification of $T(G/H)$ by $G$-orbits are
\begin{eqnarray*}
\dim(G\cdot[g,X+\mf h])
%&=&\dim(G/K)+\dim(gHg^{-1}\cdot [g,X+\mf k])\\
%&=&\dim(G/K)+\dim(K\cdot [\mathbf{1},X+\mf k])\\
&=&\dim(G/H)+\dim(K\cdot (X+\mf h))
\end{eqnarray*}
Moreover, $[g,X+\mf h]$ and $[\tilde g,\tilde X+\mf h]$ belong to the same $G$-orbit
in $G\times_H (\g/\mf h)$ if and only if $X+\mf h$ and $\tilde X+\mf h$ belong to the same
 $H$-orbit in $\g/\mf h$.

Let $\mathrm{pr}_H: G\to G/H$ be the canonical projection. Its  derivative $\mathrm{pr}_H': TG\to T(G/H)$ is surjective and given by $\mathrm{pr}_H'(g,X)=[g,X+\mf h]$ in the notation of homogeneous vector bundles. Suppose that we have an $H$-invariant decomposition $\g=\mf h+\mf q$. Then we can identify  $T_o(G/H)\cong \g/\mf h$ with $\mf q$. This gives an identification $T(G/H)\cong G\times_H\mf q$ via
\[[g,X]\longleftrightarrow \lambda'_g(o)X.\]

The vector fields on $G/H$ correspond to smooth functions $\oline{X}\colon G\to \g/\mf h$ satisfying
$$\forall g\in G, h\in H: \quad \oline{X}(gh)=h^{-1}\cdot \oline{X}(g).$$
The vector field is then given by $gH\mapsto [g,\oline{X}(g)]$.
%\item[(v)] 
A vector field $gH\mapsto [g,\oline{X}(g)]$ is $G$-invariant
if and only if $\oline X$ is a constant map:
\begin{eqnarray*}
[x,\oline X(x)]
&=&\left(g\cdot[\bullet, \oline X(\bullet)]\right)(x)
\ \overset{(\ref{vfactiondef})}{=}\  g\cdot[g^{-1}x, \oline X(g^{-1}x)]
\ =\ [x,\oline X(g^{-1}x)].
\end{eqnarray*}
In that case the value of
$\oline X$ has to be a $H$-fixed point in $\g/\mf h$.
Therefore the map
$${\cal V}(G/H)^G\to (\g/\mf h)^H,\quad \oline{X}\mapsto \oline{X}(o),$$
where ${\cal V}(G/H)^G$ denotes the $G$-invariant vector fields and
$(\g/\mf h)^H$ the $H$-fixed points, is a linear isomorphism.
If $\oline{X}(o)=X+\mf h$ with $X\in \g$, then the flow of
$\oline{X}$ is given by
\begin{equation}\label{invhom2}
\varphi^X_t(g\cdot o)=g(\exp tX) \cdot o.
\end{equation}
In particular, the invariant vector fields are complete.
\label{rem:homogeneous vector bundles}
\end{rem}

Using the description of $T(G/K)$ given in Remark~\ref{rem:homogeneous vector bundles} we can show the following proposition.

\begin{prop}[\small{cf. \cite[Lemmas~8.4 \& 8.6]{Hi05}}]
~Let $G$ act on the left of $G/M\times \mf a$ via left translation on the first factor. Then the following holds.
\begin{enumerate}
\item[{\rm (i)}] The map
$$
\Phi\colon G/M\times \mf a\to T(G/K)\cong G\times_K \mf p,\quad
(gM,X)\mapsto\lambda_g'(o)X=[g,X]
$$
is $G$-equivariant, surjective, and the following two statements
are equivalent:
\begin{enumerate}
\item[{\rm(1)}] ~$\Phi(gM,X)=\Phi(\tilde gM,\tilde X)$.
\item[{\rm(2)}] ~There exists $k\in K$ with $gk=\tilde g$ and
                $\Ad(k)X=\tilde X$.
\end{enumerate}
In particular, $\tilde X$ is contained in the Weyl group orbit
of $X$ and we have a bijection of orbit spaces
$$G\backslash T(G/K) \longleftrightarrow  W(\mf a,K)\backslash \mf a.$$

\item[{\rm (ii)}] $\Phi$ is smooth and surjective. Nevertheless, $\Phi$ is not a covering.
            The derivative $\Phi'(gM,X)$ is bijective if and only if
            $X$ is regular. More precisely, 
            $$\Kern \Phi'(e,X)=\sum_{X\in \Kern\alpha} \mf k_\alpha\quad\mbox{ and }\quad
              \im\Phi'(e,X)=\mf p\times (\mf a+\sum_{X\not\in\Kern\alpha}\mf p_\alpha).$$

\end{enumerate}
\label{orbitstructure}\label{a-action}
\end{prop}

Set $P:=\exp {\mf p}$ and note that the Cartan decomposition $G=PK$
 yields a retraction
$$j: G/K\to P,\quad pK\mapsto p$$
for the canonical quotient map  $\mathrm{pr}_K$. We identify  $G/K$ with $P\subseteq G$, so that $T(G/K)$ can be viewed as a subset of $TG\cong G\times \g$. Then the derivative $\exp':T{\mf p}\cong {\mf p}\times{\mf p}\to TP\subseteq TG\cong G\times \g$ is given by
\begin{eqnarray*}
\exp'(X,Y)
&\overset{(\ref{invhom1})}{\cong}&\Big(\exp X, \sum_{k=1}^\infty \textstyle\frac{(-\ad X)^{k-1}}{k!} Y\Big)=: (\exp X,\Psi(X)Y).
\end{eqnarray*}
and yields a coordinate system for $T(G/K)$.  

We recall descriptions the $G$-actions on $G/K$ and $T(G/K)$ in our various identifications. We consider the projection
$$
\mathrm{pr}_P :  G=PK\to P,\quad pk\mapsto p
$$
and note that $p^2=(pk)\theta(pk)^{-1}$, where $\theta\colon G\to G$ is the Cartan involution corresponding to the Cartan decomposition $G=PK$. The $G$-action on $G/K$ by left translation becomes $g\cdot p=\mathrm{pr}_P(gp)$ under the above identification of $G/K$ with $P=\exp\mf p\cong \mf p$. Thus
$$(g\exp X)\theta(g\exp X)^{-1}=g(\exp X)\theta(\exp X)^{-1}\theta(g)^{-1}= g(\exp 2X)\theta(g)^{-1}$$
for $g\in G$ and $X\in {\mf p}$ shows that the induced action on ${\mf p}$ is given by
\begin{equation}
\label{Gonp}
g\cdot X=\textstyle{\frac{1}{2}}\log(g(\exp 2X)\theta(g)^{-1}).
\end{equation}

\begin{lem}[\small{\cite[Lemma~8.2]{Hi05}}]
~The action of $G$ on $T(G/K)\cong {\mf p}\times {\mf p}$ is given by
\begin{eqnarray*}
k\cdot(X,Y)&=&(\Ad(k)X,\Ad(k)Y)\\
p\cdot(X,Y)&=&\left(\textstyle{\frac{1}{2}}\log(p(\exp 2X)p),
               \big(\Ad(p)\circ\Psi(\log(p(\exp 2X)p))\big)^{-1}\Psi(2X)Y\right)
\end{eqnarray*}
for $k\in K$,  $p\in P$ and $X,Y\in {\mf p}$.
\label{G-action}
\end{lem}

The Riemannian metric on $G/K$  allows us to identify $T(G/K)$ and the
cotangent bundle $T^*(G/K)$. The cotangent bundle carries a natural symplectic structure which we will have to study in some detail. In this context the iterated tangent bundle $T(T(G/K))$ will be of importance. \cite[\S~1.3]{Pa99} contains a detailed geometric description. Here we give a description in group theoretic terms which allows explicit calculations for Weyl chamber flows. To this end  
we consider the derivative
$$f'\colon T(T{\mf p})=({\mf p}\times {\mf p})\times ({\mf p}\times {\mf p})\to
T(TG)\cong(G\times \g)\times (\g\times \g)$$
of the map $f:=\exp'$. It is given by 
$$f'(X,Y)(A,B)\cong(\exp X,\Psi(X)Y,\Psi(X)A,(\Psi'(X)A)Y+\Psi(X)B).$$
For $X=0$ the formula for $f'$ simplifies to
$$f'(0,Y)(A,B)=(\1,Y,A,B-\textstyle{\frac{1}{2}}[A,Y]).$$
Under the identification of the tangent bundle $T(TG)$ with $(G\times \g)\times(\g\times \g)$ and the corresponding identification of $T(T(G/K))$ as a subset of $(G\times \g)\times(\g\times \g)$ for $Y=(0,Y)\in T_o (G/K)\cong {\mf p}$, we have
\begin{equation}
\label{tangembed}
T_Y(T(G/K))=\{(\1,Y,A,-\textstyle{\frac{1}{2}}[A,Y]+B)\mid A,B\in {\mf p}\}.
\end{equation}
The other tangent spaces  are then determined via $G$-invariance.

We consider the canonical projection 
$\mathrm{pr}_{T^*}: T^*(G/K)\to G/K$ and its derivative
$\mathrm{pr}_{T^*}': T(T^*(G/K))\to T(G/K)$.
Then the set
\[\{\xi\in T(T^*(G/K))\mid \mathrm{pr}_{T^*}'(\xi)=0\}\]
is called the \emph{vertical subbundle} of $T(T^*(G/K))$. Analogously one has a vertical subbundle of $T(T(G/K))$ using the canonical projection $\mathrm{pr}_{T}: T(G/K)\to G/K$. It is clear that the identification $T(G/K)\cong T^*(G/K)$ preserves the vertical bundles. Since $T_v(T(G/K))$ is canonically identified
with $T_{\mathrm{pr}_{T}(v)}(G/K)\times T_{\mathrm{pr}_{T}(v)}(G/K)$,
the Riemannian metric on $G/K$ induces a Riemannian metric on $T(T(G/K))$. Therefore it
makes sense to speak about the \emph{horizontal subbundle}, i.e. the bundle of
elements orthogonal to the elements of the vertical subbundle.

Using the above identifications the vertical and horizontal bundles can be described explicitly in group theoretic terms. 

\begin{lem}[\small{\cite[Lemma~8.1]{Hi05}}]
~Viewing $T_Y(T(G/K))$ as a subset of $T(TG)=(G\times \g)\times (\g\times \g)$ we have:
\begin{enumerate}
\item[{\rm (i)}] The space of vertical vectors is
$\{(\1,Y,0,B)\mid B\in {\mf p}\}.$
\item[{\rm (ii)}] The space of horizontal vectors is
$\{(\1,Y,A,-\textstyle{\frac{1}{2}}[A,Y])\mid A\in {\mf p}\}$.
\end{enumerate}
\end{lem}

\begin{lem}[\small{\cite[Lemma~8.3]{Hi05}}]
Under the identification 
\begin{eqnarray*}
T_{(0,Y)}({\mf p}\times {\mf p})={\mf p}\times{\mf p}&\to& \g\times \g\cong T_{(\1,Y)}(G\times \g),\ %\\
(A,B)\mapsto(A,B-\textstyle{\frac{1}{2}}[A,Y])
\end{eqnarray*}
we have
\begin{enumerate}
\item[{\rm(i)}]
The vector fields $\tilde Z\in{\cal V}(T(G/K))$
induced by the action of $\exp \R Z$ are given by
\begin{eqnarray*}
\forall Z\in \mf k,Y\in {\mf p}:&& \tilde Z(\1,Y)=(\1,Y,0,[Z,Y])\quad \\
\forall Z\in {\mf p}, Y\in {\mf p}:&&\tilde Z(\1,Y)=(\1,Y,Z,-\textstyle{\frac{1}{2}}[Z,Y]).
\end{eqnarray*}
In particular, $\tilde Z$ is a section of the vertical bundle for $Z\in \mf k$
and a section of the horizontal bundle for $Z\in{\mf p}$.

\item[{\rm(ii)}] \label{G-orbitkor}
The tangent space $T_{(\1,Y)}(G\cdot(\1,Y))$ is given
by
$$\{(\1,Y,A,[B,Y]-\textstyle{\frac{1}{2}}[A,Y])\mid A\in {\mf p},B\in \mf k\}
\cong {\mf p}\times \ad(Y)\mf k.$$
It contains all the horizontal vectors in $T_{(\1,Y)}(T(G/K))$.
\end{enumerate}
\label{G-orbit}
\end{lem}

\begin{rem}[The right $A$-action on $\M$]
As $ A$ commutes with $M$, the space $\mathcal M$ from Remark~\ref{rem:Gamma-quotients} carries a right $A$-action. The tangent bundle $T(G/M)$ can be canonically identified with the homogeneous bundle $G\times_M (\mf g/\mf m)$. With respect to this identification the derived $A$-action on $T(G/M)$ is given by $[g,X]\cdot a=[ga,\Ad(a^{-1})X]$. To describe the properties of this action in more detail, we consider the Lie algebra $\n=\sum_{\alpha\in\Sigma+}\g_\alpha$ of $N$ and its ``opposite'' $\bbar\n:=\sum_{-\alpha\in\Sigma^+}\g_\alpha$.  The decomposition $\mf g=\bbar{\mf n}+\mf m+\mf a+\mf n$ allows us to identify $\mf g/\mf m$ with  
$\bbar{\mf n}+\mf a+\mf n$. Note that each of these spaces is $M$-invariant, so $T(G/M)$  is decomposed as the Whitney sum
\[T(G/M)= (G\times_M \bbar{\mf n}) \oplus  (G\times_M \mf a) \oplus (G\times_M \mf n).\]
As the left $\Gamma$-actions do not interfere with the right $A$-action this sum descends to $\M$:
\begin{equation}
T\M=(\Gamma\backslash G\times_M \bbar{\mf n})\oplus (\Gamma\backslash  G\times_M \mf a) \oplus (\Gamma\backslash G\times_M \mf n).
\end{equation}
This decomposition will show that the $A$-action  is an Anosov flow in the sense of Subsection~\ref{subsec: Anosov flows}. Note that $T_e(G/M)$, where $e:= M\in G/M$ is the canonical base point, can be identified with $\g/\mf m$. This in turn can also be identified with $\sum_{\alpha\in \Sigma^+} \mf k_\alpha$.
\label{rem:A-action on G/M}
\end{rem}

\subsubsection{Invariant differential operators}\label{sec:harishchandra}
Let $\mathbb D(G/K)$ be the algebra of \emph{$G$-invariant differential operators} on $G/K$, i.e. differential operators commuting with the natural $G$-action on $C^\infty(G/K)$. Then we have an algebra isomorphism $\mathrm{HC}\colon \mathbb D(G/K)\to I(\mf a_\C^*)$ from $\mathbb D(G/K)$ to the $W$-invariant complex polynomials on $\mf a^\ast_\C$ which is called \emph{Harish-Chandra homomorphism} (see \cite[Ch.~II, Theorem 5.18]{gaga}). For $\lambda\in \mf a^\ast_\C$ let $\chi_\lambda$ be the character of $\mathbb D(G/K)$ defined by $\chi_\lambda(D)\coloneqq \mathrm{HC}(D)(\lambda)$. Obviously, $\chi_\lambda= \chi_{w\lambda}$ for $w\in W$. Furthermore, the $\chi_\lambda$ exhaust all characters of $\mathbb D(G/K)$ (see \cite[Ch.~III, Lemma~3.11]{gaga}). We define the space of joint eigenfunctions 
\begin{equation}\label{eq:Elambda}
E_\lambda \coloneqq\{f\in C^\infty(G/K)\mid \forall D\in  \mathbb D(G/K):\ Df = \chi_\lambda (D) f \}.
\end{equation}
Note that $E_\lambda$ is $G$-invariant.

\begin{rem}[Harish-Chandra isomorphism]
Combining the Harish-Chandra isomorphism with a result of Chevalley we find isomorphisms
$$\D(G/K)\cong{\cal E}_K'(G/K)\cong U(\mf a)^W\cong I(\mf a_\C^*)
\cong \C[X_1,\ldots,X_\ell],$$
where $\ell=\dim_\R\mf a$ is the rank of $G/K$, ${\cal E}_K'(G/K)$ denotes the space of distributions supported in the base point, and $U(\mf a)$ is the universal enveloping algebra of $\mf a_\mathbb C$. In fact, the generators of $I(\mf a_\C^*)$ can be chosen real, i.e.
as real polynomials on $\mf a^*$, which then come from differential operators
with real valued coefficients. % (see \cite{Hu90}, \S 3.5).
\label{Chevalley}
\end{rem}

\subsection{Representation theory and harmonic analysis} 

In the context of quantum-classical correspondences for locally symmetric spaces three aspects of representation theory  and harmonic analysis play key roles. According to the Plancherel formula for symmetric spaces the \emph{spherical principal series representations} are the building blocks of $L^2(G/K)$, so they may be considered as quantum objects. At the same time they are realized on the boundary of the symmetric space. In rank one the set of geodesics can be parameterized by two boundary points, namely their limits for time $t$ tending to $\pm\infty$. This ties the boundary to the underlying classical geometry of the sphere bundle. Moreover, relating geodesics to pairs of boundary points allows to interprete \emph{standard intertwining operators} (often called \emph{Knapp-Stein operators}) as classical scattering operators. It is possible to move functions on the boundary to functions on the symmetric space preserving the symmetries and hence relating associated spectral invariants. The objects doing this are the \emph{Poisson tranforms}. They are instrumental in the  quantum-classical correspondences we will describe. 

The intuitive picture just drawn for rank one symmetric spaces does have higher rank generalizations which, however, are a little less intuitive. Moreover, $G$-equivariance properties of the maps involved will allow to move from symmetric spaces to locally symmetric spaces. 

\subsubsection{Principal series representations for minimal parabolics}\label{sec:principalreps}

The concept of a principal series representation is an important tool in representation theory of semisimple Lie groups. It can be described using different  pictures. We start with the \emph{induced picture}: Pick $\lambda\in \mf a_\C^\ast$ and  an irreducible unitary representation $(\sigma,V_\sigma)$ of $M$. We define 
\[
V_{\sigma,\lambda} \coloneqq \Big\{f\colon G\overset{\text{cont.}}{\to} V_\sigma \ \Big| \  \begin{matrix}\forall g\in G,m\in M, a\in A,n\in N:\\
f(gman)=e^{-(\lambda+\rho)\log a}\sigma(m)^{-1} f(g)
\end{matrix}\ \Big\}
\]
endowed with the norm $\|f\|^2 := \int_K \|f(k)\|_\sigma^2\mathrm{d}k$, where $\mathrm{d}k$ is the normalized Haar measure on $K$, $2\rho$ is the sum of all positive restricted roots weighted by multiplicity, and $\|\cdot\|_\sigma$ is the norm on $V_\sigma$. The group $G$ acts on $V_{\sigma,\lambda}$ by the left regular representation. The completion $H_{\sigma,\lambda}$ of $V_{\sigma,\lambda}$ with respect to the norm is called \emph{induced picture of the (not necessesarily unitary) principal series representation} with respect to $(\sigma, \lambda)$. We also write $\pi_{\sigma,\lambda}$ for this representation. If $\sigma$ is the trivial representation, then we write $H_\lambda$ and $\pi_\lambda$ and call it the \emph{spherical principal series} with respect to $\lambda$. Note that for equivalent irreducible unitary representations $\sigma_1,\sigma_2$ of $M$ the corresponding principal series representations are equivalent as representations as well. In particular, the Weyl group $W$ acts on the unitary dual of $M$ by $m_w\sigma (m) = \sigma (m_w^{-1} m m_w)$, where $w\in W$ is given by a representative $m_w\in M'$ and therefore $H_{\lambda,w\sigma}$ is well-defined up to equivalence.

The \emph{compact picture} is given by restricting the function $f\colon G\to V_\sigma$ to $K$, i.e. a dense subspace is given by 
$$ \{f\colon K\to V_\sigma \text{ cont.}\mid \forall k\in K,m\in M:\ f(km)=\sigma(m)^{-1} f(k)\}$$ 
with the same norm as above. In this picture the $G$-action is given by $$\pi_{\sigma,\lambda}(g)f(k)= e^{-(\lambda+\rho)H(g^{-1} k)}f (k_{KAN} (g^{-1} k)),\quad g\in G, k\in K,$$ where $k_{KAN}$ is the $K$-component in the Iwasawa decomposition $G=KAN$. 

Recall the associated homogeneous vector bundle $\mc V_\sigma:= G\times_M V_\sigma$ over $G/M$ built in analogy to $G\times_K\mf p$  considered in Remark~\ref{rem:homogeneous vector bundles}. We can identify the principal series representation $H_{\sigma,\lambda}$ with $L^2$-sections of an associated bundle. If $\mc V_\sigma^{K/M}$ denotes the restriction of  the vector bundle $\mc V_\sigma$ over $G/M$ to $K/M\subseteq G/M$ we obtain the principal series representation with parameters $(\sigma,\lambda)$ as the Hilbert space of  $L^2$-sections of the bundle $\mc V_\sigma^{K/M}$ over the Furstenberg boundary $K/M$ of $G/K$ with the action $$\ov {\pi_{\sigma,\lambda}(g)f}(k)= e^{-(\lambda+\rho)H(g^{-1} k)}\ov f (k_{KAN} (g^{-1} k)).$$

\begin{rem}[Extension of continuous $G$-representations] For a continuous $G$-representation on some locally convex complete Hausdorff topological vector space $V$ one considers the spaces $V^\infty$  of \emph{smooth vectors}  consisisting of the $v\in V$ for which the map $g\mapsto gv$ is smooth. If one equips $V^\infty$ with the topology induced by the embedding $V^\infty\hookrightarrow C^\infty(G,V)$, then $V^\infty$ is a $G$-invariant Fr\'echet space and the $G$-representation on $V^\infty$ is continuous (see e.g. \cite[\S~4.4.1]{Wa72}). Recall the \emph{contragredient representation} $\check V$ of $G$ from \cite[\S~4.1.2]{Wa72}, where
\[\check V:=\{\nu\in V'\mid G\to V_\mathrm{b}',\ g\mapsto \nu\circ g^{-1} \text{ cont.}\}\]
and $V_\mathrm{b}'$ is the topological dual of $V$ equipped with the topology of bounded convergence. Then the strong dual of $\check V^\infty$ is denoted by $V^{-\infty}$ and called the space of \emph{distribution vectors} of $V$. 

Similarly one defines the space $V^\omega$ of \emph{analytic vectors} in $V$ and the space of \emph{hyperfunction vectors} $V^{-\omega}$ of $V$ (see \cite[\S~4.4.5]{Wa72}). Together one finds the following chain of continuous inclusions (cf. \cite[p.~220]{Ju01})
\[V^\omega\hookrightarrow V^\infty \hookrightarrow V \hookrightarrow V^{-\infty} \hookrightarrow V^{-\omega}.\]
\label{rem:extension of G-representations}
\end{rem}

\begin{rem}[Poisson transforms]
The representation of $G$ on the eigenspace $E_\lambda$ from \eqref{eq:Elambda} can be described via the \emph{Poisson transform}:  The Poisson transform $\mc P_\lambda$ is a $G$-equivariant map $H_{-\lambda}^{-\omega}\to E_\lambda$. It is given by 
\[\mc P_\lambda f (xK):= \int _{K} f(k)e^{-(\lambda+\rho)H(x^{-1} k)} \intd k\] 
if $f$ is a sufficiently regular function in the compact picture of the principal series. If $f$ is given in the induced picture, then $\mc P_\lambda f(xK)$ simply is  $\int_K f(xk) dk$.

It is important to know for which values of $\lambda\in \mf a_\C^\ast$ the Poisson transform is a bijection.
By \cite{KKMOOT} we have that $\mc P_\lambda\colon H_{-\lambda}^{-\omega}\to E_\lambda$ is a bijection if \begin{align}\label{eq:assumptionA}
-\frac{2\langle \lambda,\alpha\rangle}{\langle \alpha,\alpha\rangle}\not \in \N_{>0}  \quad \text{for all} \quad \alpha \in \Sigma^+.
\end{align}
Here $\langle\cdot,\cdot\rangle$ is the inner product on $\mf a^*$ induced from the Killing form. So, in particular, $\mc P_\lambda$ is a bijection if $\Re \lambda$ is in the closure $\ov{\mf  a_+^\ast}$ of ${\mf  a_+^\ast}=\{\lambda\in\mf a^*\mid \forall \alpha\in \Sigma^+: \langle \lambda,\alpha\rangle>0\}$.
\label{sec:Poissontransform}
\end{rem}

\subsubsection{Intertwining operators}\label{subsubsec: intertwiner}

For $\lambda \in i\a^{\ast}$ and $w\in W$, we know that $\pi_{\lambda}$ and  $\pi_{w\lambda}$ are equivalent. An \textit{intertwining operator} $\mathsf A(\lambda,w)$ can be defined as follows.
For $w\in W$, we set $\bar{N}_w:=\bar{N}\cap m_w^{-1}Nm_w$, where as before $m_w \in M'$ is a representative of $w$, and $\bar{N}=\theta(N)$. Recall here that $\theta$ is the Cartan involution on $G$. We normalize the Haar measure on $\bar{N}_w$ by $\int_{\bar N_w}e^{-2\rho(H(\bar n_w))}\intd\bar n_w=1$. For  $\lambda \in \a_{\C}^{\ast}$ and
\bq
\label{eqn: D-lambda}
f \in\mathcal{D}_{\lambda}(G):=\Big\{f \in \Cinft(G)\ \Big|\ 
\begin{matrix}
\forall g \in G, m\in M, a\in A, n\in N:\\
f(gman)=e^{(\lambda-\rho)\log a}f(g)
\end{matrix}\ \Big\},
\eq
we define the intertwining integrals
\bqn
\label{eqn: A-lambda-w}
\mathsf{A}(\lambda, w)f(g):=\int_{\bar N_w}f(gm_w\bar n_w) \intd\bar n_w.
\eqn
This integral is absolutely convergent in the region 
$$C(w):=\{\lambda \in \a_{\C}^{\ast}\; | \; \forall \; \alpha \in \Sigma(w):\ 
\Re\langle \lambda, \alpha\rangle > 0 \;  \},$$ 
where $\Sigma(w):=\{\alpha \in \Sigma^+\; | \; w\cdot \alpha<0\}$.

 With the topology induced from that of $\Cinft(G)$, $\mathcal{D}_{\lambda}(G)$ is Fr\'echet, and is isomorphic to $\Cinft(G/MAN)$.
Using the decomposition $G=KAN$ we see that this linear isomorphism is just the restriction map $\phi\mapsto \phi|_{K}$ with the inverse mapping  given by $\psi \longmapsto \phi(kan):=e^{(\lambda-\rho)\log a}\psi(k)$, $(k,a, n)\in K \times A\times N$. 
To $\phi \in C^\infty_c(G)$, we can associate $\phi_{\lambda} \in \Cinft(G/MAN)$ which is defined by
\bqn
\phi_{\lambda}(kM):=\int_{MAN}f(kman)e^{(\lambda+\rho)\log a}\intd m \intd a \intd n.
\eqn
Observe that the mapping $\phi\mapsto \phi_{\lambda}$ is surjective (see \cite[Chapter 7, \S 2, Proposition 2]{Bourbaki04int2}).
This then leads to the corresponding identification at the level of distributions as follows.

We denote by $\mathcal{D}'(G/MAN)$ the space of distributions on $G/MAN\cong K/M$. Given $T \in \mathcal{D}'(G/MAN)$, we define $T_{\lambda} \in \mathcal{D}'(G)$ by  $T_{\lambda}(\phi)=T(\phi_{\lambda}),$ where $\phi_{\lambda} \in \Cinft(G/MAN)$ is as above. Then 
$$T_{\lambda} \in \mathcal{D}'_{\lambda}(G):=\{T \in \mathcal{D}'(G)\;|\; \forall m \in M, a \in A, n\in N:\  R_{man}T=e^{(\rho-\lambda)\log a}T \},$$ 
the space of distributions on $G$ that are equivariant under the right regular action $R_{man}$ of $MAN$ in the specified way, and $T \longmapsto  T_{\lambda}$ is an isomorphism between $\mathcal{D}'(G/MAN)$ and $\mathcal{D}_{\lambda}'(G).$ Thus, we can and we will move freely between the two spaces.

It can be shown (see e.g. \cite[\S\S VII.3--7]{knapp} and \cite[Cor. to Prop.~5]{STS76})
that we have the following properties.

\begin{prop}[Standard intertwining operators]
For $\lambda \in C(w)$,
\begin{enumerate}
\item[{\rm(i)}]~ $\mathsf{A}(\lambda, w): \mathcal{D}'_{\lambda}(G) \longrightarrow \mathcal{D}'_{w \cdot \lambda}(G)$ is a continuous $G$-homomorphism.
\item[{\rm(ii)}]~ $\mathsf{A}(\lambda, w_1w_2)=\mathsf{A}(w_2\cdot\lambda, w_1)\mathsf{A}(\lambda, w_2)$ for a minimal decomposition $w=w_1w_2$.
\item[{\rm(iii)}]~ $\mathsf{A}(\lambda, w)\circ \pi_{\lambda}=\pi_{w\cdot \lambda}\circ \mathsf{A}(\lambda, w)$.
\item[{\rm(iv)}]~ By considering the principal series representations in the compact realization on $\mathcal{D}'(G/P)$, so that the space is independent of the parameter $\lambda$, the operators $\mathsf{A}(\lambda, w)$ extend meromorphically in $\lambda$ to all of $\a_{\C}^{\ast}$, and the properties (i), (ii), and (iii) still hold.
\end{enumerate}
\end{prop}

For each $w\in W$ we can define an intertwining operator $\mathsf A(w)$ on $L^2(i\a^\ast;L^2(K/M))$ by
$$(\mathsf A(w)f)(\lambda):=\mathsf A(\lambda,w)f(\lambda),$$
where we interpret $\mathsf A(\lambda,w)$ as an operator on $L^2(K/M)$.

\subsubsection{$\Gamma$-cohomology}\label{subsubsec: Gamma cohomology}
It is obvious that analysis on a locally symmetric space $\mc X=\Gamma\backslash G/K$ leads to the question how to find $\Gamma$-invariant objects in spaces defined from $G/K$ in a $G$-equivariant way. If the constructions are functorial but not exact, cohomological obstructions to the existence of $\Gamma$-invariant objects may occur. Thus it is not surprising that $\Gamma$-cohomology plays a role in the context of locally symmetric spaces.

Recall the definition of \emph{$\Gamma$-cohomology} $H^\bullet(\Gamma,V)$ with coefficients in a $\Gamma$-module $V$, e.g. from \cite[p.~218]{Ju01}. For $p\ge 0$ the space $C^p(\Gamma,V)$ of \emph{homogeneous $p$-cochains} consists of all functions $c:\Gamma^{p+1}\to V$. The \emph{coboundary operator} $\partial: C^p(\Gamma,V)\to C^{p+1}(\Gamma,V)$ is given by
\[(\partial c)(\gamma_0,\ldots,\gamma_{p+1}):= \sum_{j=0}^{p+1}(-1)^j c(\gamma_0,\ldots,\widehat{\gamma_j},\ldots,\gamma_{p+1}),\]
where $\widehat{\gamma_j}$ means that this entry is ommitted. Then $\partial\circ \partial=0$ and 
$\Gamma$ acts on $C^p(\Gamma,V)$ via
\[(\gamma\cdot c)(\gamma_0,\ldots,\gamma_{p}):= \gamma(c(\gamma^{-1}\gamma_0,\ldots,\gamma^{-1}\gamma_{p})).\]
The spaces of $\Gamma$-invariants form a complex
\[0 \longrightarrow  C^0(\Gamma,V)^\Gamma\overset{\partial}{\longrightarrow} 
 C^1(\Gamma,V)^\Gamma\overset{\partial}{\longrightarrow}
  C^2(\Gamma,V)^\Gamma\overset{\partial}{\longrightarrow}\ldots\]
whose cohomology is $H^\bullet(\Gamma,V)$. In particular, $H^0(\Gamma,V)$ coincides with the space $V^\Gamma$ of $\Gamma$-invariants in $V$.

\begin{rem}[Patterson conjecture]  In \cite{Ju01} $\Gamma$-cohomologies are mostly taken with coefficients related to principal series representations with the goal to explain Patterson's conjectures on cohomological descriptions of the divisor of dynamical zeta functions and the proofs of these conjectures by Bunke and Olbrich (see in particular  \cite[\S~3.4 and Chap.~8]{Ju01}). While the conjecture relates two dynamical objects belonging to the classical domain, the proof given by Bunke and Olbrich involves Poisson transforms and thus eigenspaces of Laplacians which belong to the quantum world.
\end{rem}

\begin{rem}[Lewis-Zagier correpondences]
In the paper \cite{BLZ15} Bruggeman, Lewis and Zagier  generalized  results from  \cite{LZ01} to various discrete subgroups $\Gamma$ of $G:=\mathrm{PSL}_2(\R)$. They concentrated on proving linear isomorphisms between spaces of Maass cusp forms and $\Gamma$-cohomology spaces with coefficients in principal series representations. In the case of hyperbolic surfaces with cusps, they showed that as for the modular surface, \emph{parabolic} $\Gamma$-cohomology (see \cite[\S~11.3]{BLZ15}) is a suitable cohomology theory to use for this purpose. 
\end{rem}

Instead of recalling the definition of parabolic $\Gamma$-cohomology used in \cite{BLZ15} we make a little detour which allows us to also introduce \emph{cuspidal $\Gamma$-cohomology} which played a decisive role in \cite{DH05}, a first attempt to extend the Lewis-Zagier correspondence to higher rank situations.  

As we want to consider arithmetic groups we now assume that $G$ is algebraic.  Let $\Ga$ be an arithmetic subgroup of  $G$ and assume that $\Ga$ is torsion-free. Then every $\Ga$-module $V$ induces a local system or locally constant sheaf $\mc C_V$ on $\mc X=\Gamma\bs G/K$. In the \'etale picture the sheaf $\mc C_V$ equals $\mc C_V=\Ga\bs ((G/K)\times V)$, (diagonal action). Let $\overline{\mc X}^\mathrm{BS}$ denote the Borel-Serre compactification \cite{BS73} of $\mc X$, then $\Ga$ also is the fundamental group of $\overline{\mc X}^\mathrm{BS}$ and $V$ induces a sheaf also denoted by $\mc C_V$ on $\overline{\mc  X}^\mathrm{BS}$. This notation is consistent as the sheaf on ${\mc X}$ is indeed the restriction of the one on $\overline{\mc X}^\mathrm{BS}$. Let $\partial^\mathrm{BS}\mc X$ denote the boundary of the Borel-Serre compactification. There are natural identifications
$$
H^j(\Ga,V)\ \cong\ H^j(\mc X,\mc C_V)\ \cong
H^j(\overline{\mc X}^\mathrm{BS},\mc C_V).
$$
Now \emph{parabolic cohomology} of a $\Ga$-module $V$ can be defined as the
kernel of the restriction to the boundary,
$$
H_\mathrm{par}^j(\Ga,V) := \ker\left( H^j(\overline{\mc X}^\mathrm{BS},\mc C_V)\longrightarrow H^j(\partial^\mathrm{BS} \mc X,\mc C_V)\right).
$$
The long exact sequence of the pair $(\overline{\mc X}^\mathrm{BS},\partial^\mathrm{BS}\mc X)$
gives rise to
$$
\dots\longrightarrow H_c^j(\mc X,\mc C_V)\longrightarrow H^j(\mc X,\mc C_V)=
H^j(\overline{\mc X}^\mathrm{BS},\mc C_V)\longrightarrow H^j(\partial^\mathrm{BS}\mc X,\mc C_V)\longrightarrow\dots,
$$
where $H_c^j$ denotes cohomology with compact supports. The image of the cohomology with compact supports under the natural map is called the \emph{interior cohomology} of $\Ga\bs X$ and is denoted by $H_!^j(\mc X,\mc C_V)$. The exactness of the above sequence shows that
$$
H_\mathrm{par}^j(\Ga,V)\ \cong\ H_!^j(\mc X,\mc C_V).
$$
It is shown in \cite[\S~2]{DH05} that  $H^\bullet(\Gamma,V)$ is computed by a standard complex of $(\mf g,K)$-cohomology (see \cite[\S~I.5]{BW00}) if $V$ is a \emph{smooth $G$-representation}, i.e. if $V=V^\infty$. More precisely, one then has
\[H^\bullet(\Gamma,V)\cong H^\bullet(\mf g,K, C^\infty(\Gamma\bs G)\hat \otimes V),\]
where $C^\infty(\Gamma\bs G)\hat \otimes V$ is the completed projective tensor product isomorphic to $C^\infty(\Gamma\bs G,V)$ (see \cite[p.~80/81]{Gr55}). 

Recall the space of automorphc cusp forms $L^2_\mathrm{cusp}(\Gamma\bs G)$ (see e.g. \cite[\S~I.2]{Ha68}), which is a closed $G$-invariant subspace of $L^2(\Gamma\bs G)$. Then, for a smooth $G$-representation $V$ the \emph{cuspidal $\Gamma$-cohomology} is defined by
\[H^\bullet_\mathrm{cusp}(\Gamma,V):= H^\bullet(\mf g,K, L^2_\mathrm{cusp}(\Gamma\bs G)^\infty\hat \otimes V).\] 
Note that there is a natural map $H^\bullet_\mathrm{cusp}(\Gamma,V)\to H^\bullet(\Gamma,V)$ which is not necessarily injective (see \cite[Cor.~5.2]{DH05}). Its image is called the \emph{reduced cuspidal $\Gamma$-cohomology} and denoted by $\tilde H^\bullet_\mathrm{cusp}(\Gamma,V)$. By \cite[Prop.~2.2]{DH05} we have 
\[\tilde H^\bullet_\mathrm{cusp}(\Gamma,V)\subseteq H_\mathrm{par}^\bullet(\Ga,V).
\]

The representation $L^2_\mathrm{cusp}(\Gamma\bs G)$ decomposes discretely into irreducibles with finite multiplicity (see e.g. \cite[Part I]{Ar79}). We denote the multiplicity of an ireducible $G$-representation $\pi$ in $L^2_\mathrm{cusp}(\Gamma\bs G)$ by $N_\Gamma(\pi)$.

\begin{thm}[\small{\cite[Theorems~0.2 \& 4.3]{DH05}}]
~Let $\Ga$ be a torsion-free arithmetic subgroup of a split semisimple Lie group $G$ and let $\pi\in\hat G$ be an irreducible unitary principal series representation. Then
\[
N_{\Ga}(\pi) = \mathrm{dim}~H_\mathrm{cusp}^{d-\ell}(\Ga,\pi^\omega)= \mathrm{dim}~H_\mathrm{cusp}^{\ell}(\Ga,\pi^{-\omega}),
\]
where $d=\dim G/K$ and $\ell$ is the real rank of $G$.
\label{thm:DH05-0.2}  
\end{thm}

For $G$ non-split the assertion of Theorem~\ref{thm:DH05-0.2} remains true for a generic set of representations $\pi$.

\section{Dynamical systems}

In this section we discuss some  mathematical specifics of dynamical systems entering the classical and quantum evolutions of our class of examples.

\subsection{Classical dynamics}\label{sec: ClassDynamics}

The mathematics of modelling  classical mechanics revolves around vector fields and their integration. We start by collecting some background material before we describe our key examples, the geodesic and the Weyl chamber flows.

\subsubsection{Anosov flows} \label{subsec: Anosov flows}

Let $\M$ be a compact manifold, $A\simeq \R^\ell$ an abelian group and let $\tau: A\to \Diffeo(\M)$ be a smooth locally free group action. 
If $\a:= \tu{Lie}(A)\cong \R^\ell$, we can define a \emph{generating map} 
\[
 X:\Abb{\a}{C^\infty(\M;T\M)}{H}{X_H:=\frac{d}{dt}_{|t=0}\tau(\exp(tH)),}
\]
so that for each basis $H_1,\dots,H_\ell$ of $\a$, $[X_{H_j},X_{H_k}]=0$ for all $j,k$. For $H\in \a$ we denote by $\varphi_t^{X_H}$ the flow of the vector field $X_H$. Notice that, as a differential operator, we can view $X$ as a map 
\[X:C^\infty(\M)\to C^\infty(\M;\a^*),\quad (Xu)(H):=X_{H}u.\] 

It is customary to call the action \emph{Anosov} if there is an $H\in \a$ such that there is a continuous $d\varphi_t^{X_H}$-invariant splitting 
\begin{equation}\label{splittingA}
 T\M=E_0\oplus E_\mathrm{u}\oplus E_\mathrm{s},
\end{equation}
where $E_0:={\rm span}(X_{H_1},\dots,X_{H_\ell})$, and there exists a $C>0,\nu>0$ such that for each 
$x\in \M$ 
\begin{align*}
&\forall w\in E_\mathrm{s}(x),\forall t\geq 0:\quad   \|d\varphi_t^{X_H}(x)w\|\leq Ce^{-\nu t}\|w\|, \\
&\forall w\in E_\mathrm{u}(x),\forall t\leq 0:\quad    \|d\varphi_t^{X_H}(x)w\|\leq Ce^{-\nu |t|}\|w\|.
\end{align*}
Here the norm on $T\M$ is fixed by choosing any smooth Riemannian metric on $\M$. We say that  such an $H$ is \emph{transversely hyperbolic}.  The letters s and u in $E_\mathrm{s}$ and $E_\mathrm{u}$ stand for \emph{stable} and \emph{unstable}.
The splitting is invariant by the entire action. This does not mean that all $H\in \a^*$ have this transversely hyperbolic behavior. In fact, there is a maximal open convex cone $\mf a_+^\tau\subset \a$ containing $H$ such that for all $H'\in\mf a_+^\tau$, $X_{H'}$ is also transversely hyperbolic with the same splitting as $H$ (\cite[Lemma~2.2]{GGHW21}); $\mf a_+^\tau$ is called a \emph{positive Weyl chamber}. This name is motivated by the classical examples of such Anosov actions that are the Weyl chamber flows for locally symmetric spaces of rank $\ell$ as described in Section~\ref{subsec:Weyl chamber flow}. For other classes of examples see e.g. \cite{KaSp94, SV19}.

\subsubsection{Symplectic geometry and Hamiltonian actions}\label{subsubsec: symp-geo}

Let $({\cal N},\omega)$ be a symplectic manifold with tangent bundle $T{\cal N}$ and cotangent bundle $T^*{\cal N}$. Then for each smooth function $f\in C^\infty({\cal N})$ one defines its \emph{symplectic gradient}  or \emph{Hamiltonian vector field} $X_f$ via
$$-\omega (X_f,X)=\la df,X\ra_{\cal N},$$
where $X\in {\cal V}({\cal N})$ is a vector field on ${\cal N}$,
and $\la\cdot,\cdot\ra_{\cal N}$ denotes
the canonical pairing between $T^*{\cal N}$ and $T{\cal N}$
(for the basic facts of symplectic geometry see e.g. \cite{AM87}).
The Hamiltonian vector fields satisfy
\begin{equation}
\label{ProdHam}
X_{f_1 f_2}=f_2 X_{f_1}+f_1 X_{f_2}.
\end{equation}
For $f_1,f_2\in C^\infty({\cal N})$ one has the \emph{Poisson bracket} $\{f_1,f_2\}\in C^\infty({\cal N})$ defined by
$$\{f_1,f_2\}:=\omega(X_{f_1},X_{f_2}).$$
It defines a Lie algebra structure on $C^\infty({\cal N})$ and the equality
$\{f_1,f_2\}= d{f_2}(X_{f_1})$ shows that for two Poisson commuting
functions the one function (and hence its level sets) is invariant
under the flow of the other's Hamiltonian vector field.

Suppose that a group $G$ acts on $\cal N$  by diffeomorphisms such that $\omega$ is  $G$-invariant. If  
$f\in C^\infty({\cal N})$ is  $G$-invariant as well, then the Hamiltonian
vector field $X_f\in {\cal V}({\cal N})$
is again $G$-invariant. As a consequence we see that the
flow $\varphi^{X_f}_t(x)$ of $X_f$ is $G$-equivariant, i.e.
$$\varphi^{X_f}_t(g\cdot x)=g\cdot \varphi^{X_f}_t(x)\quad \mbox{for } (x,t)\in {\cal N}\times\R,$$
where $(x,t)$ is restricted to the domain of definition for the flow.

Let ${\cal A}\subseteq C^\infty({\cal N})$ be a finitely generated associative subalgebra consisting of Poisson commuting functions. Let  $f_1,\ldots,f_k$ be a set of generators and define $F\colon {\cal N}\to \R^k$ by $F:=(f_1,\ldots, f_k)$. Then an arbitrary element of ${\cal A}$ can be written as $f=\sum_{\alpha\in \N_0^n} c_\alpha F^\alpha$, where we  use the usual multi-index notation. For any $v\in \R^k$ we consider the closed subset ${\cal N}_{(v)}:=F^{-1}(v)$ of ${\cal N}$. Since the $f_j$ Poisson-commute we see that all the ${\cal N}_{(v)}$ are stable under all the flows of the $X_{f_j}$'s.

If now  ${\cal N}_{(v)}$ is a submanifold of ${\cal N}$, then the
Hamiltonian vector fields $X_{f_j}$ are all tangent to ${\cal N}_{(v)}$, i.e.
restrict to vector fields on ${\cal N}_{(v)}$. Note that all elements of
${\cal A}$ are constant on ${\cal N}_{(v)}$. But then the identity
$(\ref{ProdHam})$ implies that the restriction of any $X_f$ with
$f\in {\cal A}$ is a linear combination of the restrictions of the $X_{f_j}$.

%\item[(iii)]
In general, the ${\cal N}_{(v)}$ will not be manifolds. It is, however,
always possible to find submanifolds of ${\cal N}$ contained in ${\cal N}_{(v)}$
such that all the vector fields $X_{f_k}$ restrict to vector fields of these
submanifolds.
To this end we recall the notion of a ${\cal D}$-orbit for a family
${\cal D}$ of vector fields from \cite{Su73}. We set
\begin{equation}\label{Ddef}
{\cal D}:=\{X_{f}\mid f\in {\cal A}\}
\end{equation}
and note that
according to $(\ref{ProdHam})$ the linear spans of
$\{X_f(x)\mid f\in {\cal A}\}$ and $\{X_{f_1}(x),\ldots,X_{f_k}(x)\}$
in $T_x{\cal N}$ agree for all $x\in {\cal N}$.
We denote this linear span by $\Delta(x)$.
Then $x\mapsto \Delta(x)$ is a smooth distribution in the sense of
\cite{Su73}. It is involutive since ${\cal D}$ is a commutative family of
vector fields. Moreover, this shows that $\Delta$ satisfies the
condition (e) from \cite[Thm. 4.2]{Su73}. Therefore the ${\cal D}$-orbits
${\cal D}\cdot x$ are maximal integral manifolds of $\Delta$. They satisfy
\begin{equation}\label{Dorbitdim1}
\dim({\cal D}\cdot x)=\dim \Delta(x).
\end{equation}

%\item[(iv)]
From the construction of the ${\cal D}$-orbits it is clear that
they are invariant under the flows of the $X_f\in {\cal D}$. Since these
flows preserve the ${\cal N}_{(v)}$ we see also that the ${\cal D}$-orbits are
contained in single ${\cal N}_{(v)}$'s. This, finally, shows that the restriction of
$X_f$ with $f\in {\cal A}$ to any ${\cal D}$-orbit is a linear
combination (with coefficients depending on $v$) of the restrictions
of the $X_{f_1},\ldots,X_{f_k}$.

%\item[(v)] We take up the situation from (i).
Suppose that $y$ is a regular value of $F$, i.e. the derivative
$F'(x)\colon T_x{\cal N}\to \R^k$ is surjective for all $x\in F^{-1}(y)$.
Then ${\cal N}_{(F(x))}$ is a closed
submanifold of ${\cal N}$ with $\dim  {\cal N}_{(F(x))}= \dim {\cal N} - k$.
On the other hand, $k$ is equal to
$$\dim\Span\{f_1'(x),\ldots,f_k'(x)\}
=\dim\Span\{X_{f_1}(x),\ldots,X_{f_k}(x)\}
=\dim ({\cal D}\cdot x)
$$
so that
\begin{equation}\label{Dorbitdim2}
\dim ({\cal D}\cdot x) + \dim  {\cal N}_{(F(x))}= \dim {\cal N}.
\end{equation}

%\item[(vi)]
The commuting vector fields $X_{f_1},\ldots, X_{f_k}$ define a (local) action
of $\R^k$  on ${\cal N}$ which leaves  the ${\cal D}$-orbits invariant.
%\end{enumerate}
%}\hfill\halmos
%\end{bemerkung}

\begin{rem}[Hamiltonian actions on cotangent bundles]
Let $\mathcal{N}$ be any manifold. Recall the \emph{canonical symplectic form} $\omega=d\Theta$ on $T^*\mathcal{N}$, where the $1$-form $\Theta$ is given by the formula
$$\forall v\in T_\xi(T^*\mathcal{N}), \xi\in T^*\mathcal{N}:\quad \la\Theta(\xi),v\ra= \la \xi, \mathrm{pr}_{T^*}'(v)\ra.$$
with the canonical projection $\mathrm{pr}_{T^*}:T^*\mathcal{N}\to \mathcal{N}$. For a diffeomorphism $h\colon \mathcal{N}_1\to \mathcal{N}_2$ the induced map $h^*\colon T^*\mathcal{N}_2\to T^*\mathcal{N}_1$ is a \emph{canonical transformation}, i.e. it preserves the  symplectic structures. Thus the induced map
$$C^\infty(T^*\mathcal{N}_1)\to C^\infty(T^*\mathcal{N}_2),\quad f\mapsto h_*f:=f\circ h^*$$
is a Poisson isomorphism satisfying $h_*X_f=X_{h_*f}$. Therefore the push-forward of vector fields $h_*\colon {\cal V}(T^*\mathcal{N}_1)\to {\cal V}(T^*\mathcal{N}_2)$ maps Hamiltonian vector fields to Hamiltonian vector fields.

If now $G$ acts smoothly on $\mathcal{N}$, these considerations show
that the Hamiltonian vector fields of
$G$-invariant functions on $T^*\mathcal{N}$ are $G$-invariant Hamiltonian
vector fields on $T^*\mathcal{N}$. Moreover,  the action
of $G$ on $T^*\mathcal{N}$ is Hamiltonian with moment map
$J\colon T^*\mathcal{N}\to \g^*$ given by
$$\la J(\eta), X\ra:= \la \eta, \tilde X\circ \mathrm{pr}_{T^*}(\eta)\ra_\mathcal{N},$$
where $X\in \g$ and $\tilde X\in {\cal V}(\mathcal{N})$ is the vector field on
$\mathcal{N}$ induced by the derived action of $G$ on $\mathcal{N}$ via
$$\tilde X(x):=\frac{d}{dt}\vert_{t=0}(\exp tX)\cdot x.$$

If $T^*\mathcal{N}$ and $T\mathcal{N}$ are identified via a Riemannian metric $g$, the formula for $\Theta$ reads
$$\la\Theta(\xi),v\ra=\la \xi, \mathrm{pr}_{T^*}'(v)\ra=
g\big(\xi, \mathrm{pr}_{T}'(v)\big)$$
for $\xi\in T^*\mathcal{N}\cong T\mathcal{N}$,
$v\in T_\xi(T^*\mathcal{N})\cong T_\xi(T\mathcal{N})$ and the canononical projection $\mathrm{pr}_{T}:T\mathcal{N}\to \mathcal{N}$.
\end{rem}

In the case of $\mathcal{N}=G/K$, using $T^*(G/K)\cong G\times_K{\mf p}^*$ and identifying ${\mf p}$ and ${\mf p}^*$ via the Killing form  gives 
$$T^*(G/K)\cong G\times_K{\mf p}\cong T(G/K).$$
Under these identifications the moment map is given by $J([g,X])=\Ad(g)X$.

\subsubsection{Geodesic and Weyl chamber flow}\label{subsec:Weyl chamber flow}

Recall the notation for locally symmetric spaces and the right $A$-action on $T\M$ from Section~\ref{sec:loc sym spaces}. All elements of the positive Weyl chamber $\a_+$ are transversely hyperbolic elements sharing the same stable and unstable distributions given by the associated homogeneous vector bundles
 \[
  E_0 = \Gamma\backslash {G}\times_{{M}}\a,\quad E_\mathrm{s}=\Gamma\backslash {G}\times_{M}\n, \quad E_\mathrm{u} = \Gamma\backslash {G}\times_{M}\bbar{\n} .
 \]
Thus, if $\X=\Gamma\backslash {G}/K$ is compact the $A$-action is Anosov. Below, we will recall the symplectic interpretation of the $A$-action from \cite{Hi05}. We  will not discuss the origin of the name Weyl chamber flow for the $A$-action in detail. Let it suffice to say that the objects involved can be reinterpreted in such a way that $A$ acts on the set of Weyl chambers.

We give a symplectic interpretation of the right $A$-action on $\M=\Gamma\backslash G/M$. The natural symplectic form on $T(T(G/K))$ is given by
\begin{equation}\label{symform}
\omega_{(\1,Y)}\big((A_1,B_1-\textstyle{\frac{1}{2}}[A_1,Y]),
           (A_2,B_2-\textstyle{\frac{1}{2}}[A_2,Y])\big)=
    B(A_1,B_2)-B(A_2,B_1),
\end{equation}
where $B$ denotes the Killing form.

Let $f\in C^\infty(T(G/K))$ be $G$-invariant and $h\in C^\infty({\mf p})$ the restriction of $f$ to $T_o (G/K)= {\mf p}$ (which determines $f$ uniquely). Then the Hamiltonian vector field $X_f$ is determined by
\begin{equation}\label{hamfield}
X_f(X)=\left(\grad h(X),-\textstyle{\frac{1}{2}}[\grad h(X),X]\right),
\end{equation}
where the gradient is taken with respect to the Killing form on ${\mf p}$. In particular, $X_f$ is horizontal, so by Lemma \ref{G-orbitkor}
the flow of $X_f$ preserves
$G$-orbits.  Writing such an orbit as $G/H$, according to (\ref{invhom2}) we obtain a flow $\varphi^{X_f}\colon \R\times G/H\to G/H$
of $X_f$ on  $G/H\subseteq T(G/K)$ of the form
\begin{equation}\label{flowgrad}
\varphi^{X_f}_t([g,X])=[g\exp(t\,\grad h(X)),X]
\end{equation}
and preserves $G$-orbits. Therefore the flow  can be studied on
these orbits separately.
Note further
that (\ref{symform}) and (\ref{hamfield}) show that any two $G$-invariant functions
$f_1,f_2\in C^\infty(T(G/K))$ commute under the Poisson bracket,
i.e.
\begin{equation}\label{Poissoninvcomm}
\forall f_1,f_2\in C^\infty(T(G/K))^G:\quad \{f_1,f_2\}=0.
\end{equation}

\begin{thm}[\small{\cite[Thm.~8.1]{Hi05}}]  
~Consider the  algebra ${\cal A}_\mathrm{pf}$  of $G$-invariant
functions in $C^\infty(T^*(G/K))$ which restrict
to polynomials on the fiber $\mf p$. Then ${\cal A}_\mathrm{pf}$ is finitely
generated, commutative under the Poisson bracket,
and  the joint level sets of these functions are precisely the
$G$-orbits in $T^*(G/K)$.
 \label{levelsets}
\end{thm}

Let $\mathcal A$ be the algebra of $G$-invariant functions in $C^\infty(T^*(G/K))$. In view of Dadok's smooth Chevalley theorem, \cite{Da82}, it is easy to show that restriction yields an isomorphism between $\mathcal A$ and $C^\infty(\mf a^*)^W$. 

By Remark~\ref{Chevalley} the associative subalgebra of $C^\infty(T(G/K))\cong C^\infty(T^*(G/K))$ generated by the principal symbols of elements of $\D(G/K)$ coincides with the complexification of the algebra ${\cal A}_\mathrm{pf}$ introduced in Theorem~\ref{levelsets}.
Note that for a $K$-invariant smooth function $h\colon \mf p\to \R$ we have $\grad h(X)\in \mf a$ for all $X\in \mf a$, so that $\grad h(X)=\grad h\vert_\mf a(X)$ (see e.g. \cite[Lemma~8.5]{Hi05}). Together with (\ref{flowgrad}) this shows that we have a map
$$
{\cal F}\colon{\cal A}\times \mf a\to\mf a,\quad
(f,X)\mapsto\grad f\vert_\mf a(X)
$$
such that the flow $\varphi^{X_f}\colon \R\times T(G/K)\to T(G/K)$
of $X_f$ on  $T(G/K)$ is given by
\begin{equation}\label{Faction}
\varphi^{X_f}_t([g,X])=[g\exp(t\,{\cal F}(f,X)),X].
\end{equation}
In particular,  $\mathcal{A}$ is Poisson-commutative as well.
Extending the right $A$-action on $G/M$ trivially to $G/M\times \mf a$ via
$$(gM,X)\cdot a:=\rho(a)(gM,X):=(gaM,X)$$
we obtain the following equivariance property
of $\Phi\colon G/M\times \mf a\to T(G/K)$ under  $\varphi^{X_f}$:
\begin{equation}
\label{equivA}
\Phi\circ \rho\left(e^{t\,\grad f(X)}\right)
=\varphi^{X_f}_t\circ \Phi.
\end{equation}

%\begin{ex} \label{gedesicflow}
Recall, e.g. from \cite[Ex.~8.1]{Hi05}, that the geodesic flow $\varphi_t$ on the tangent bundle
$T(G/K)\cong G\times_K \mf p$  can be written
$$\varphi_t([g,X])=[g\exp tX,X].$$
This implies in particular that each geodesic, viewed as a curve in $T(G/K)$, i.e. as an orbit of the geodesic flow, is completely contained in a  $G$-orbit.

If $G/K$ has rank $1$, i.e. if $\dim\mf a=1$, the $K$-invariant polynomials on $\mf p$ are generated by $X\mapsto \|X\|^2$, so that the $G$-orbits in $T(G/K)$ are simply the sphere bundles (and the zero section). As this function essentially represents the kinetic energy, the $G$-orbits can be viewed as energy shells being preserved by the flow. The relevant $T^*(G/K)_{(v)}$ for nonzero $v$ is (up to scaling) the sphere bundle $S(G/K)$ of $G/K$ in $T(G/K)$. Left translation gives an identification of $S(G/K)$ with $G/M$. The geodesic flow restricted to the sphere bundle is then given by the right multiplication of $A$. For $v=0$ the bundle is the zero section and the geodesic flow restricted to this section is trivial.

In higher rank we have not only the kinetic energy function which is invariant by the $A$-action. All the elements of $\mathcal A_\mathrm{pf}$ are invariant under the $A$-flow and the joint level sets may rightfully be called \emph{generalized energy shells}. By Theorem~\ref{levelsets} they agree with the $G$-orbits. Studying the multi-parameter flow on the $G$-orbits separately, allows to give more explicit decriptions. To this end  we fix $X\in \mf a$ and recall from \cite{Hi05} the description of the sets $\{\grad h(X)\mid h\in \R[\mf a]^W\}$, where $\R[\mf a]^W$ denotes the $W$-invariant real polynomialson $\mf a$.

Note that the centralizer  $\mf k_X$ of $X$ in $\mf k$ is given by
$$\mf z_\mf k(X):=\{Y\in\mf k\mid [X,Y]=0\}.$$
Analogously, given $X\in\mf a$  consider the centralizer
$\g_X:=\mf z_\g(X):=\{Y\in \g\mid [Y,X]=0\}.$
If $\Sigma_X=\{\alpha\in \Sigma\mid \alpha(X)=0\}$,
then $\g_X=\mf m+\mf a+\sum_{\alpha\in\Sigma_X} \g_\alpha$.
Consider  the center $\mf z(\g_X)$  of the reductive
$\theta$-invariant algebra $\g_X$ and
$\mf a(X):=[\g_X,\g_X]\cap \mf a$. Then $\mf a(X)$ is a Cartan
subspace for $[\g_X,\g_X]$ and $\mf a=\mf a(X)\oplus \mf a_X$, where
\begin{equation}\label{EXdef}
\mf a_X:=\mf p\cap \mf z(\g_X).
\end{equation}

\begin{prop}[\small{\cite[Lemma~8.9]{Hi05}}] 
~Let $f_1,\ldots,f_\ell$ be a set of algebraically independent
generators of ${\cal A}_\mathrm{pf}$ and $X\in \mf a$. Then
$$\mf a_X=\sum_{j=1}^\ell \R\,\grad f_j(X).$$
 \label{singorbit2}
\end{prop}

\begin{thm}[\small{\cite[Thm.~8.2]{Hi05}}]
The smooth distribution $\Delta$ associated with the family of vector fields
${\cal D}=\{X_f\mid f\in {\cal A}_\mathrm{pf}\}$ introduced in (\ref{Ddef}) is given by
$$\Delta([\mathbf{1},X])=\mf a_X\times \{0\}\in \mf p\times\mf p\cong T_{[\mathbf{1},X]}(G\times_K\mf p)\cong T_{[\mathbf{1},X]}(T(G/K))$$
for $X\in\mf a$.
\label{thm:Bialo8.2}  
\end{thm}

\begin{rem}[Singular flows]~ Recall that each $G$-orbit in $T^*(G/K)=G\times_K\mf p^*$ is of the form $G\cdot[\mathbf 1,\xi]$ with $\xi\in \mf a^*$. The stabilizer of  $[\mathbf 1,\xi]$ in $G$ coincides with the stabilizer $K_\xi$ of $\xi$ with respect to the coadjoint action.
Fix some  $\xi\in \mf a^*$ and identify  $\mf g$ with $\mf g^*$. Under this identification $\mf a$ gets identified with $\mf a^*$, so we can set $A_\xi:=\exp \mf a_\xi$. Then $A_\xi$ commutes with $K_\xi$ so that we have a right $A_\xi$-action on $G\cdot[\mathbf 1,\xi]\cong G/K_\xi$ via
\[ (gK_\xi,a)\mapsto gaK_\xi =:(gK_\xi)\cdot a.\]
If $\xi$ is regular, then $\mf a_\xi=\mf a$ and $K_\xi=M$. Thus we recover the Weyl chamber flow. If $\xi=0$, then $K_\xi =K$ and $A_\xi=\{\mathbf 1\}$, so that we find the trivial action on $G/K$.
\label{rem:singular flows}
\end{rem}

\subsection{Quantum dynamics}

While classical dynamics is based on flows of Hamiltonian vector fields, quantum dynamics is based on unitary linear one-parameter groups of operators generated by skew-adjoint densely defined operators via Stone's theorem, see e.g. \cite[\S~35]{La02}.

Multiplication by $i$ allows to shift attention from skew-adjoint to self-adjoint operators. As quantization procedures often start with differential operators an important question is to decide which of the differential operators one considers can actually be completed to self-adjoint operators. Such operators are called \emph{essentially self-adjoint}. There are theorems relating essential self-adjointness of a quantized Hamiltonian to the geodesic completeness of the underlying classical Hamiltonian vector field. This is why some people call essential self-adjointness \emph{quantum completeness}. A very nice discussion of this subject can be found in Tao's blog \cite{Ta11}.

In the situation we are interested in, namely viewing invariant differential operators as generating the quantum evolution, we may refer to \cite{vdB87}, where van den Ban shows that formally self-adjoint invariant differential operators are essentially self-adjoint, i.e. quantum complete. This should be viewed as a counterpart of the fact that the Hamiltonion vector fields associated with functions from the algebra ${\mc A}_\mathrm{pf}$ introduced in Subsection~\ref{subsec:Weyl chamber flow} are complete.

\section{Quantum-classical correspondences}

In this section we will reach the quantum-classical correspondences for locally symmetric spaces which are in the focus of this paper. Before we can describe these correspondences we have to introduce the dynamical invariants on the quantum and the classical side which will be connected.

%\subsection{Dynamical invariants}

Recall that both kinds of dynamical systems we consider are built from linear operators. In the quantum case this is built into the mathematical formulation of quantum mechanics. In the classical case we concentrate on Hamiltonian mechanics which is based on vector fields, i.e. first order differential operators. Thus considering spectral properties of these operators is a reasonable starting point. 

For the Hamiltonian operators generating the quantum dynamics on Hilbert spaces, there is a canonical spectral theory to consider. But already the example of the Dyatlov-Faure-Guillarmou correspondence for compact hyperbolic surfaces (see Example~\ref{ex:DFG}) shows that setting up the right kind of spectral invariants for the first order differential operators is subtle. The situation becomes even more complicated if we consider noncompact surfaces. In \cite{GHW18} it was shown that it is possible to extend the theory of Ruelle resonances to noncompact surfaces without cusps, but in order to obtain a correspondence one has to replace the no longer discrete ordinary $L^2$-spectrum of the hyperbolic Laplacian by the discrete \emph{resonance spectrum} of this operator. 

The dynamical invariants we will discuss here are all of a spectral nature. Resonances will play a distinguished role, but we will also consider decompositions of representations as in Example~\ref{ex:Flaminio-Forni}. Finally, we will also briefly look at zeta functions.

\subsection{Quantum invariants}

The kind of dynamical invariants that comes to mind first in the context of quantum mechanics is the spectrum of the (often unbounded) linear operators involved. Things get more complicated when  one also looks at commutants of Hamiltonians and noncommutative harmonic analysis enters because of their symmetries. 

\subsubsection{Joint spectrum of the algebra of invariant differential operators in the cocompact case}\label{subsubsec: quantum spectrum cocompact}

In the rank one case the quantization of the geodesic flow is given by the Laplacian on $G/K$. In the higher rank case we have to consider the algebra of $G$-invariant differential operators on $G/K$ which we denote by $\mathbb D(G/K)$. As an abstract algebra this is a polynomial algebra with $\ell=\mathrm{rank}(G/K)$ algebraically independent operators, among them the Laplace operator of Remark~\ref{Chevalley}. These operators descend to $\X=\G\backslash G/K$ and in analogy to \eqref{eq:Elambda} we can define the joint eigenspace 
\begin{equation}\label{eq:GammaElambda}
{}^\G E_\lambda:=\{f\in C^\infty(\G\backslash G/K)\mid \forall D\in\mathbb D(G/K):\ Df =\chi_\lambda(D)f\},
\end{equation}
where $\chi_\lambda$ is a character of $\mathbb D(G/K)$ parametrized by $\lambda\in\mf a_\C ^\ast/W$ with the Weyl group $W$. Here $\chi_\rho$ is the trivial character. Let $\sigma_\mathrm{Q}$ denote the corresponding \emph{quantum spectrum} $\{\lambda\in\mf a_\C^\ast\mid {}^\G E_\lambda\neq \{0\}\}$.

\subsubsection{Laplace resonances for noncompact locally symmetric spaces} \label{subsubsec:Resolvent poles}

Resonance spectra for the Laplace-Beltrami operators $\Delta_\X$ play a role as soon as our locally symmetric space $\X$ is no longer compact. This applies in particular to the global symmetric space $G/K$, i.e. the case of trivial $\Gamma$. 

One way to introduce \emph{resonances} of the Laplace-Beltrami operator $\Delta_\X$  for a noncompact Riemannian manifold $\X$ is to consider the Schwartz kernel $R(z)$ of the resolvent $(\Delta_\X-z^2-\|\rho\|^2)^{-1}$. Under good circumstances it is possible to extend $z\mapsto R(z)$ to a meromorphic family of Schwartz kernels whose poles are then called the resonances of  $\Delta_\X$.

In the case of (global) Riemannian symmetric spaces noncommutative harmonic analysis allows for an explicit description of the resolvent as an integral operator.  For real rank $1$ and $2$ the task of doing the meromorphic continuation turns ot to be managable. 

\begin{thm}[\small {\cite[Thm.~3.8]{HP09}}] Let $\X=G/K$ be a Riemannian symmetric space of rank $1$.
\begin{enumerate}
\item[{\rm(i)}] If $\Sigma^+=\{\alpha\}$ with $m_\alpha$ even, then $R(z)$ extends to an entire function (no poles).
\item[{\rm(ii)}] In the other cases let $\alpha$ be the simple root in $\Sigma^+$. Then $R(z)$ extends to a meromorphic function with simple poles at the points $z_k=\lambda_k\|\alpha\|$ for $k\in \N_0$ with $\lambda_k=-i(\rho+j_\X k)$, where $\rho=\frac{1}{2}m_\alpha+m_{2\alpha}$ is the weighted halfsum of positive roots and $j_\X\in \{1,2\}$ depends on $\X$.
\item[{\rm(iii)}] The residue operator $R_k:=\mathrm{res}_{z=z_k}R(z)$ is of finite rank and can be expressed explicitely as a convolution with the spherical function associated with $\lambda_k$.
\end{enumerate}
\end{thm}

The paper \cite{HP09} was written following a suggestion of M.~Zworski, who had done the calculation for the upper halfplane. It later turned out that Miatello and Will had already done (in \cite{MW00}) the case of Damek-Ricci spaces, which include the rank one Riemannian symmetric spaces. The methods of \cite{HP09}, however, are applicable in higher rank as well -- at least in principle. The technical effort to do the meromorphic continuation via contour shifts in higher dimension turns out to be formidable. In fact, it took a series of papers just to deal with rank $2$ (see  \cite{HPP16,HPP17a,HPP17b}).   

Among the noncompact locally symmetric spaces the class of spaces for which one has good control over the resonances obtained as resolvent poles is even more restricted. The following result was established in the context of asymptotically hyperbolic manifolds (see \cite{MM87}). For hyperbolic surfaces this means one needs to exclude cusps, i.e. restrict attention to the so-called convex cocompact surfaces.

\begin{thm}[\small cf. {\cite[Thm.~4.2]{GHW18}}] ~Let  $\X=\Gamma\backslash\mathbb H^{2}$ be a convex co-compact hyperbolic surface.
\begin{enumerate}
\item[{\rm(i)}] The nonnegative Laplacian $\Delta_\X$  on $\X$ has a resolvent $(\Delta_\X-s(1-s))^{-1}$ that admits a meromorphic extension from  $\{s\in \C\mid {\Re}(s)>1/2\}$ to $\C$ as a family of bounded operators $R(z): C_c^\infty(\X)\to C^{\infty}(\X)$.
\item[{\rm(ii)}]~ For a pole $z_0\not=\frac{1}{2}$ the polar part of the resolvent is of the form 
\[
\sum_{j=1}^{J(z_0)} \frac{(\Delta_\X-z_0(1-z_0))^{j-1}(1-2z_0)}{(z(1-z)-z_0(1-z_0))^j}\mathrm{res}_{z=z_0}R(z)
\]
for some $J(z_0)\in \N$.
\item[{\rm(ii)}]~If  $z_0=\frac{1}{2}$ is a  pole, the polar part of the resolvent at $z_0$ is of the form 
\[
 \frac{1}{z(1-z)-z_0(1-z_0)}\mathrm{res}_{z=z_0}R(z).
\]
\end{enumerate} 
\label{thm:GHW18-quantumresonances}
\end{thm}

Note that in this result the poles are not determined explicitly, but one can say more about the location of the resonances and the order of the poles (see \cite[Prop.~4.3]{GHW18}).

Theorem~\ref{thm:GHW18-quantumresonances} is used in the quantum-classical correspondence for convex cocompact hyperbolic surfaces established in \cite{GHW18}. Motivated by that correspondence, Hadfield in the process of proving a quantum-classical correspondence for convex cocompact real hyperbolic spaces obtained  similar but more complicated results for Bochner Laplacians  acting on symmetric tensors   (see \cite{Ha17,Ha20}).

Coming from a different motivation also Bunke and Olbrich have a result on the meromorphic continuation of resolvents. They consider Casimir operators of $G$ acting on homogeneous vector bundles over a rank one locally symmetric space $\X=\Gamma\backslash G/K$ with convex cocompact $\Gamma$, see \cite[\S~6]{BO12}. Their result is that the resolvent of the Casimir $C$ meromorphically extends to a branched cover of $\C$ with branching points given by the values of $C$ on certain principal series representations.

\subsubsection{Laplace resonances vs. scattering poles}

In Subsection~\ref{subsubsec:Resolvent poles} we introduced resonances of Laplacians as poles of meromorphic continuations of the resolvent. Physicists typically define resonances of Laplacians as poles of the scattering matrix of an associated scattering problem. They tend to take these as synonymous concepts which, strictly speaking, is not the case. However, under certain circumstances one can show in a mathematically rigorous way that resolvent poles and scattering poles are closely related. 

The reason such connections are of interest in our context is the desire to generalize quantum resonances to higher rank. The case of compact locally symmetric spaces (see Subsection~\ref{subsubsec:HWW21}) shows that instead of considering just the Laplace-Beltrami operator one should rather consider the commutative algebra of all invariant differential operators. While it is unclear what the resolvent of this commutative (associative) algebra should be, there are candidates for scattering operators derived from a suitable multitemporal wave equation. For the symmetric spaces these are the Knapp-Stein intertwining operators (see Subsection~\ref{subsubsec: intertwiner}), which are well-known to depend meromorphically on as many parameters as the invariant differential operators (via the Harish-Chandra isomorphism, see \cite[Chap.~10]{Wa92}).

Scattering operators for surfaces were studied in the context of automorphic form by Lax and Phillips in \cite{LP76} as an application of their general scattering theory. Guillop\'e and Zworski proved analytic extension of families of scattering operators  for very general hyperbolic surfaces in \cite{GZ95,GZ97}. Patterson and Perry treated the case of higher dimensional real hyperbolic spaces in \cite{PP99}. In \cite{JS00} Joshi and S\'a Barreto provided a generalization to all conformally compact spaces. On the other hand, for classical rank one convex cocompact locally symmetric spaces, Bunke and Olbrich constructed scattering operators which also have meromorphic continuation to all of $\C$, see \cite[Thm.~5.10]{BO00}.

Suppose that $\X$ is a rank one Riemannian symmetric space. Then $\a^*_\C$ can be identified with $\C$ and the spherical principal series representations $\pi_\lambda$ and $\pi_{-\lambda}$ for $\lambda\in \a^*_\C$ are intertwined by the Knapp-Stein intertwiner $\mathsf{A}(\lambda, w)$. Since we are in rank one the only nontrivial Weyl group element  is $w=-1$. Then the \emph{scattering operator} can be defined by 
\begin{equation}
\label{eq-Scatt-Intertwiner}
S(z) := \cfcn(\lambda)^{-1} \mathsf{A}(\lambda, -1), \quad \lambda=i z,
\end{equation}
where $\cfcn$ is the Harish-Chandra c-function. The scattering operator depends meromorphically on $z$.
We have the following result on \emph{scattering poles}, that is the possible
poles of the $S(z) f$, $f\in\C^\infty(K/M)$.

\begin{thm}[\small{\cite[Thm.~7.1]{HHP19}}]
Let $\X$ be a rank one Riemannian symmetric space.
The scattering poles are at most simple and located on the imaginary axis minus the origin.
The scattering poles in $\Im z>0$ are precisely the resonances.
At a resonance $z_0$, $\beta_{\rho+iz_0}$ maps the range of
$\mathrm{res}_{z=z_0}R(z)$ isomorphically onto the range of $\mathrm{res}_{z=z_0}S(z)$;
in particular, the residues of the resolvent
and of the scattering matrix have the same finite rank.
The non-resonance scattering poles are the poles of the standard intertwiner $\mathsf{A}(i z, -1)$;
they are contained in $(2i)^{-1}\N$.
At $\Im z<0$, the residue $\mathrm{res}_{z=z_0}S(z)$ is a nonzero multiple of the residue $\mathrm{res}_{z=z_0}\mathsf{A}(i z, -1)$.
\label{thm-scatt-poles-are-resonances}
\end{thm}

As for resolvent poles one also has results on the relation between resolvent poles and scattering poles for asymptotically hyperbolic manifolds (see \cite{BP02,Gu05}).

\subsubsection{Microlocal lifts, Wigner and Patterson-Sullivan distributions}\label{Lifts and PS-distributions}

For compact  hyperbolic surfaces  Anantharaman and Zelditch in \cite{AZ} considered dynamical zeta functions as meromorphic families of distributions by integrating over the primitive closed geodesics. The resulting residues are what they call Patterson-Sullivan distributions. Their study was motivated by quantum ergodicity as the Patterson-Sullivan distributions asymptotically agree with the Wigner distributions, which can be seen by giving an alternative construction in terms of boundary values and Poisson transforms. This alternative construction was generalized in \cite{Sch10} to compact locally symmetric spaces spaces of rank one and extended to higher rank in  \cite{HHS}.

For $\X=\Gamma\backslash G/K$ we consider sequences, $(\varphi_h)_h\subset L^2(\X)$,
of normalized joint eigenfunctions which belong to the principal spectrum
of the algebra of invariant differential operators.
Using a $h$-\psdiff\ calculus $\mathrm{Op}_h$ on $\X$, in \cite{HHS} we define and study \emph{lifted quantum limits}
or \emph{microlocal lifts} as weak$*$-limit points of \emph{Wigner distributions}
\[ W_h:a\mapsto \big(\Op_{h}(a)\varphi_h\mid\varphi_h\big)_{L^2(\X)}. \]
Here, $h^{-1}$ is the norm of a spectral parameter associated with $\varphi_h$,
and $h\downarrow 0$ through a strictly decreasing null sequence.
%Here, $h^{-2}$ is asymptotic to the eigenvalue sequence of the negative Laplacian.
Lifted quantum limits are positive Radon measures supported in the cosphere bundle.
The problem of quantum ergodicity asks for a description of the lifted quantum limits.
Using the boundary values of the $\varphi_h$ on $G/MAN=K/M$, we construct \emph{Patterson-Sullivan distributions} on $\M=\Gamma\backslash G/M$.
In the context of quantum ergodicity, Patterson-Sullivan distributions are relevant
because they are asymptotically equivalent to lifted quantum limits and satisfy invariance properties. In fact, for compact hyperbolic surfaces $\mc X=\Gamma\backslash\mathbb{H}$,
the asymptotic equivalence of lifted Wigner distributions and Patterson--Sullivan distributions
was observed by Anantharaman and Zelditch \cite{AZ}. While it was known from earlier work that lifted quantum limits on compact hyperbolic surfaces are invariant under geodesic flows it turned out that Patterson--Sullivan distributions are themselves invariant under the geodesic flow.
Moreover, in \cite{AZ} it is shown that they have an interpretation in terms of
dynamical zeta functions which can be defined completely in terms of the geodesic flow.

We assume that  $\Gamma$ is  co-compact and torsion free.
Under the diagonal action, there is a unique open $G$-orbit $(G/MAN)^{(2)}\cong G/MA$ in $G/MAN\times G/MAN$.
For rank one spaces, $(G/MAN)^{(2)}$ is the set of pairs of distinct boundary points.
In this case each geodesic of $\mc X$ has a unique forward limit point and a unique backward limit point in $G/MAN$.
In particular, one can identify $(G/MAN)^{(2)}$ with the space of geodesics.
In higher rank the geometric interpretation is more complicated.
It involves the Weyl chamber flow rather than the geodesic flow.

Joint eigenfunctions come with a spectral parameter $\lambda\in \mf a^*_\C$.
The principal part of the spectrum comes from the purely imaginary spectral parameters.
We assume that the spectral parameter of $\varphi_h$ is $i\nu_h/h \in i\mf a^*$, $|\nu_h|=1$.
The Patterson--Sullivan distribution $PS^\Gamma_h\in\Dprime(\Gamma\backslash G/M)$
associated with $\varphi_h$ is constructed as follows.
The Poisson transform  (see Remark~\ref{sec:Poissontransform}) allows us to write
$\varphi_h(x)=\mc P_{i\frac{\nu_h}{h}}(T_h)$,
where $T_h\in\Dprime(G/MAN)$ is the boundary value of $\varphi_h$.
Consider the \emph{weighted Radon transform} $\mathcal R_h:\Ccinfty(G/M)\to \Ccinfty(G/MA)$ defined by
\[ (\mathcal R_hf)(gMA) =\int_A d_h(gaM,\nu_h)f(gaM)\intd a \]
with a weight function $d_h$ given explicitly in terms of the Iwasawa decomposition, \cite[Def.~4.1]{HHS}.
Denote by  $\mathcal R_h':\Dprime(G/MAN\times G/MAN)\to\Dprime(G/M)$ the dual of $\mathcal R_h$.
The Patterson--Sullivan distribution $PS_h^\Gamma\in \Dprime(\Gamma\backslash G/M)$
is defined as the $\Gamma$-average of ${\mathcal R}_h' (T_h\otimes \tilde T_h)$, where $\tilde T_h$ is given by the complex conjugate $\overline{\varphi_h}=\mc P_{i\frac{w_0\cdot\nu_h}{h}}(\tilde T_h)$ of $\varphi_h$. Here $w_0$ is again the longest element of $W$. 

Let $W_0:=\lim_h W_h\in\Dprime(T^*\X)$ be a lifted quantum limit which, after passing to
a subsequence if necessary, has a regular limit direction $\nu_0=\lim_h \nu_h$.
In addition, assume
\[\nu_h=\nu_0+\bigoh(h)\quad\text{as $h\downarrow 0$.} \]
To link $W_0$ to the sequence $(PS_h^\Gamma)_h$ of Patterson--Sullivan distributions,
we make use of the natural $G$-equivariant map $\Phi\colon G/M\times\mf a^*\to T^*X$ from Proposition~\ref{orbitstructure}.
For regular $\nu_0\in \mf a^*$ this induces a push-forward of distributions,
\[
\Phi(\cdot,\nu_0)_*:\Dprime(\Gamma\backslash G/M)\to \Dprime(T^* \X).
\]
Then a simplified version of the main result on Patterson-Sullivan distributions in \cite{HHS} can be stated as follows:
\begin{equation}
\label{W-PS-diagonal}
W_0 = \kappa(w_0\cdot\nu_0) \lim_{h\downarrow 0} (2\pi h)^{\dim N/2} \Phi(\cdot,\nu_0)_* PS^\Gamma_h
\quad\text{in $\Dprime(T^*\X)$.}
\end{equation}
Here $\kappa$ is a normalizing function defined in terms of structural data of $G/K$.
We point out that \cite[Thm~7.4]{HHS} is more general as it also describes the situation arising from off-diagonal Wigner distributions
$\big(\Op_{\Gamma,h}(a)\varphi_h\mid \varphi'_h\big)_{L^2(\X)}$.

In the course of proving \eqref{W-PS-diagonal} the following theorem, which is of relevance in the context of \emph{quantum unique ergodicity} (see \cite{SV07} and in particular \cite[\S~5]{Si15}) was established.  

\begin{thm}[{\cite[Thm.~6.7]{HHS}, \cite[Theorem 1.6(3)]{SV07}, \cite[Theorem 1.3]{AS}}]
Assume that $(\varphi_h)_h$ has the lifted quantum limit $W_0$. Then $\supp(W_0)\subset S^*\X$, and $W_0$ is invariant under the geodesic flow.
Moreover, $\supp W_0$ is contained in a joint level set of $\mathcal A$, i.e. in a $G$-orbit in $S^*\X$. Moreover, for every $f\in \mathcal A$, $W_0$ is invariant under the Hamiltonion flow generated by $f$. If the direction $\nu_0\in\mf a^*$ of $W_0$ is regular, then $W_0$ is $A$-invariant.
\label{Kap5-Prop-supp-invar-Wigner}
\end{thm}

\subsection{Classical invariants}

In contrast to the situation in quantum mechanics, it is clear from the very beginning that classical dynamics offers a variety of different invariants studied in the theory of dynamical system. Examples are the lengths of closed geodesics, invariant measures, dynamical zeta functions and so on. But also spectral invariants of vector fields, viewed as differential operators, come in. As for the quantum case extra room for investigation is added when one looks at commutants, which here enter as flow invariant functions, and symmetry properties of Hamiltonians.

\subsubsection{Ruelle-Taylor resonances for compact Riemannian manifolds}\label{subsubsec: RT-resonances}

Let $\tau:A\to \mathrm{Diffeo}(\M)$ be a smooth abelian Anosov action on a compact Riemannian  manifold $\M$ with positive Weyl chamber $\mf a_+^\tau$ as in Subsection~\ref{subsec: Anosov flows}. Associated with a basis $H_1,\ldots, H_\ell$ for $\mf a$ we have commuting vector fields $X_{H_1},\dots,X_{H_\ell}$ which we view as first order differential operators. It is natural to consider a joint spectrum for this family. We may choose the $H_j$'s to be transversely hyperbolic with the same splitting.  Guided by the case of a single Anosov flow (done in \cite{BuLi,FaSj,DyZw}), we define $E_\mathrm{u}^*\subset T^*\mc{M}$ to be the subbundle such that\footnote{This may look strange as $E_\mathrm{u}^*$ is not the dual space of $E_\mathrm{u}$, but it should be read rather as $(E^*)_\mathrm{u}$, the unstable part of $E^*$. } $E_\mathrm{u}^*(E_\mathrm{u}\oplus E_0)=0$. We shall say that $\lambda=(\lambda_1\dots,\lambda_\ell)\in \C^\kappa$ is a \emph{joint Ruelle resonance} for the Anosov action if there is a nonzero distribution $u\in \mathcal{D}'(\M)$ with wavefront set ${\rm WF}(u)\subset E_\mathrm{u}^*$  such that
\begin{equation}
 \label{eq:intro_def_joint_res}
\forall j=1,\dots,\ell: \quad  (X_{H_j}+\lambda_j)u=0.
\end{equation}
The distribution $u$ is called a \emph{joint Ruelle resonant state}.  We will denote  the space of distributions $u$ with ${\rm WF}(u)\subset E_\mathrm{u}^*$ by $\mathcal{D}'_{E_\mathrm{u}^*}(\M)$. 

A basis free way to define joint Ruelle resonances is to call an element $\lambda\in \a_\C^*$ joint Ruelle resonance of $\tau$ if there is a nonzero $u\in \mathcal{D}'_{E_\mathrm{u}^*}(\M)$ with 
\[
\forall H\in\a: \quad  (X_{H}+\lambda(H))u=0.
\] 
We denote the Ruelle resonance spectrum of $\tau$ by $\sigma_\mathrm{R}(\tau)$. 
It is a priori not clear that the set of joint Ruelle resonances is discrete -- or nonempty for that matter -- nor that the dimension of joint resonant states is finite. In the case of compact Riemannian manifolds $\M$ this, however, turns out to be the case. In that case one has normalized volume (total volume $1$) which allows to identify $C^{\infty}(\M)$ with the space of smooth densities on $\M$ and hence $\mathcal{D}'(\M)$ with the space $C^{-\infty}(\M)$ of generalized functions. Note that $C^{\infty}(\M)\subseteq C^{-\infty}(\M)$.

\begin{thm}[\small{\cite[Thm.~1]{GGHW21}}]~Let $\tau:A\to \mathrm{Diffeo}(\M)$ be a smooth abelian Anosov action on a compact Riemannian  manifold $\M$ with positive Weyl chamber $\mf a_+^\tau$. Then the set $\sigma_\mathrm{R}(\tau)$ of joint Ruelle resonances $\lambda \in\a^*_\C$ is a discrete set contained in 
\begin{equation}\label{notaylorintro}
 \bigcap_{H\in \mf a_+^\tau} \{\lambda\in\a_\C^*\mid \Re(\lambda(H))\leq 0\}.
 \end{equation}
Moreover, for each joint Ruelle resonance $\lambda\in \a^*_\C$ the space of joint Ruelle resonant states is finite dimensional.
\label{Theo1intro}
\end{thm}

The Ruelle resonance spectrum $\sigma_\mathrm{R}(\tau)$ always contains $\lambda=0$ (with $u=1$ being the joint eigenfunction).  When $\M=\Gamma\backslash G/M$ for some compact locally symmetric space $\Gamma\backslash G/K$ and $\tau$ is the Weyl chamber action, it contains infinitely many joint Ruelle resonances (see \cite[Theorem 1.1]{HWW21}).

This theorem is by no means a straightforward extension of the case of a single Anosov flow. 
It relies on a deeper result based on the theory of joint spectrum and joint functional calculus developed by  Taylor \cite{Tay70, Tay70a}. This theory allows to set up a good Fredholm problem on certain functional spaces by using Koszul complexes. In fact, one defines $X+\lambda$, for $\lambda\in \a_\C^*$,  as an operator 
\[X+\lambda :C^\infty(\M)\to C^\infty(\M;\a_\C^*),\quad ((X+\lambda)u)(A):=(X_{A}+\lambda(A))u.\] 
For each $\lambda\in\a_\C^*$ this yields differential operators 
\[d_{(X+\lambda)}: C^\infty(\M;\Lambda^j\a_\C^*)\to C^\infty(\M;\Lambda^{j+1}\a_\C^*)\] 
via $d_{(X+\lambda)}(u\otimes\omega):= ((X+\lambda)u)\wedge\omega$  for $u\in C^\infty(\M)$ and $\omega\in \Lambda^j\a^*_\C$.  Due to the commutativity of the family of vector fields $X_{H}$ for $H\in \a$, it can be checked that $d_{(X+\lambda)}\circ d_{(X+\lambda)}=0$. Moreover, as a differential operator, it extends to a continuous map
\[
d_{(X+\lambda)}: C^{-\infty}_{E_\mathrm{u}^*}(\M;\Lambda^j\a_\C^*)\to C^{-\infty}_{E_\mathrm{u}^*}(\M;\Lambda^{j+1}\a_\C^*)
\]
and defines an \emph{associated Koszul complex}
\begin{equation}\label{introtaylorcomplex}
  0\longrightarrow C^{-\infty}_{E_\mathrm{u}^*}(\mc{M})\overset{d_{(X+\lambda)}}{\longrightarrow}C^{-\infty}_{E_\mathrm{u}^*}\otimes \Lambda^1\a_\C^* \overset{d_{(X+\lambda)}}{\longrightarrow} \dots \overset{d_{(X+\lambda)}}{\longrightarrow} C^{-\infty}_{E_\mathrm{u}^*}(\mc{M})\otimes \Lambda^\ell \a_\C^*\longrightarrow 0.
\end{equation}
In this setting one has the follwing theorem.

\begin{thm}[\small {\cite[Thm.~2]{GGHW21}}]
~Let $\tau$ be a smooth abelian Anosov action on a closed manifold $\mc M$ with generating map $X$. Then for each $\lambda\in\a_\C^*$ and $j=0,\dots,\ell$, the cohomology 
\[\Big(\ker d_{(X+\lambda)}|_{C^{-\infty}_{E_\mathrm{u}^*}(\M)\otimes \Lambda^j\a^*_\C}\Big)\Big/
\Big(\ran\, d_{(X+\lambda)}|_{C^{-\infty}_{E_\mathrm{u}^*}(\M)\otimes \Lambda^{j-1}\a^*_\C}\Big)\]
is finite dimensional. It is nontrivial only at a discrete subset of $\{\lambda\in \a^*_\C\mid \forall H\in \mf a_+^\tau:\  {\rm Re}(\lambda(H))\leq 0 \}$.
\label{Theo2intro}
\end{thm}

The statement about the cohomologies in Theorem~\ref{Theo2intro} is not only a stronger statement than Theorem~\ref{Theo1intro}, but the cohomological setting is in fact a fundamental ingredient in proving the discreteness of the resonance spectrum and its finite multiplicity. The proof in \cite{GGHW21} relies on the theory of joint \emph{Taylor spectrum}, defined using such Koszul complexes carrying a suitable notion of Fredholmness and furthermore provides a good framework for a parametrix construction via microlocal methods. More precisely, the parametrix construction is not done on the topological vector spaces $C^{-\infty}_{E_\mathrm{u}^*}(\mc M)$ but on a scale of  Hilbert spaces $\mc H_{NG}$, depending on the choice of an escape function $G\in C^\infty(T^*\mc{M})$ and a parameter $N\in \R_{\ge0}$, by which one can in some sense approximate $C^{-\infty}_{E_\mathrm{u}^*}(\mc{M})$.  The spaces $\mc H_{NG}$ are \emph{anisotropic Sobolev spaces} which roughly speaking allow $H^{N}(\mc{M})$ Sobolev regularity 
in all directions except in $E_\mathrm{u}^*$ where one allows for $H^{-N}(\mc{M})$ Sobolev regularity. They can be rigorously defined using microlocal analysis, following the techniques of Faure-Sj\"ostrand \cite{FaSj}. 
By further use of pseudodifferential and Fourier integral operator theory one can then construct a parametrix $Q(\lambda)$, which is a family of bounded operators on $\mathcal H_{NG}\otimes\Lambda \a_\C^*$ depending holomorphically on $\lambda\in\a^*_\C$ and satisfying
\begin{equation}
 \label{eq:parametrix_intro}
d_{(X+\lambda)}Q(\lambda) + Q(\lambda)d_{(X+\lambda)} = \Id + K(\lambda).
\end{equation}
Here $K(\lambda)$ is a holomorphic family of compact operators on $\mc H_{NG}\otimes\Lambda\a_\C^*$ for $\lambda$ in a suitable domain of $\a_\C^*$ that can be made arbitrarily large letting $N\to\infty$. 
Even after having this parametrix construction, the fact that the joint spectrum is discrete and intrinsic (i.e. independent of the precise construction of the Sobolev spaces) is more difficult than for an Anosov flow (the rank $1$ case): 
this is because holomorphic  functions in $\C^\ell$ do not have discrete zeros when $\ell\geq 2$ and we are lacking a good notion of resolvent, while for one operator the resolvent is an important tool. 
Due to the link with the theory of the Taylor spectrum, we call $\lambda\in \a^*_\C$ a \emph{Ruelle-Taylor resonance} for the Anosov action if for some $j=0,\ldots,\kappa$ the $j$-th cohomology is nontrivial
\[
\Big(\ker d_{(X+\lambda)}|_{C^{-\infty}_{E_\mathrm{u}^*}(\M)\otimes \Lambda^j\a^*_\C}\Big)\Big/
\Big(\ran\, d_{(X+\lambda)}|_{C^{-\infty}_{E_\mathrm{u}^*}\otimes \Lambda^{j-1}\a^*_\C}\Big)\not=0,
\]
and we call the nontrivial cohomology classes \emph{Ruelle-Taylor resonant states}. Note that the definition of joint Ruelle resonances precisely means that the $0$-th cohomology is nontrivial. Thus, any joint Ruelle resonance is a Ruelle-Taylor resonance. The converse statement is not obvious but turns out to be true, see \cite[Prop.~4.15]{GGHW21}: if the cohomology of degree $j>0$ is not $0$, then the cohomology of degree $0$ is not trivial.

We conclude this subsection with some applications of Ruelle-Taylor resonances to dynamical systems.  In view of \eqref{notaylorintro}, such a resonance is called a \emph{leading resonance} when its real part vanishes. The leading resonances carry important information about the dynamics as they are related to a special type of invariant measures as well as to mixing properties of these measures.

Let $v_g$ denote the Riemannian measure associated with the given Riemannian metric $g$ on $\mc{M}$. A $\tau$-invariant probability measure $\mu$ on $\M$ is called a \emph{physical measure} if there is $v\in C^\infty(\mc{M})$ nonnegative such that for any continuous function $f$ and any  open cone $\mathcal{C}\subset \mf a_+^\tau $,
\begin{equation}\label{eq:Cone-averaging-physical-measure}
\mu(f) =\lim_{T\to \infty} \frac{1}{{\rm Vol}(\mc{C}_T)} \int_{A\in \mc{C}_T} \int_\M f(\varphi^{-X_A}_1(x)) v(x) \intd v_g(x) \intd A
\end{equation}
where $\mc{C}_T:=\{A\in \mc{C}\, |\, |A|\leq T\}$, and here $|\cdot|$ denotes a fixed Euclidean norm on $\a$.  In other words, $\mu$ is the weak Cesaro limit of a Lebesgue type measure under the dynamics. As an application of the methods developed for the proof of Theorem~\ref{Theo2intro} one can prove the following result.

\begin{thm}[\small {\cite[Thm.~3]{GGHW21}}]~Let $\tau:A\to \mathrm{Diffeo}(\M)$ be a smooth abelian Anosov action on a compact Riemannian  manifold $\M$ with generating map $X$ and positive Weyl chamber $\mf a_+^\tau$.
\begin{enumerate}
\item[{\rm(i)}] \label{it:physM_are_res_states} The linear span over $\C$ of the physical measures is isomorphic (as a $\C$-vector space) to $\ker d_X|_{C^{-\infty}_{E_\mathrm{u}^*}(\M)}$, the space of joint Ruelle resonant states at $\lambda=0\in \a^*_\C$. In particular, it is finite dimensional. The dimension can be expressed in dynamical terms, see \cite[Theorem 3]{BGW21}.
\item[{\rm(ii)}] \label{it:WF_phys_meas} A probability measure $\mu$ is a physical measure if and only if it is $\tau$-invariant and $\mu$ has wavefront set
${\rm WF}(\mu)\subset E_\mathrm{s}^*$, where $E_\mathrm{s}^*\subset T^*\M$ is defined by $E_\mathrm{s}^*(E_\mathrm{s}\oplus E_0)=0$.
\item[{\rm(iii)}] Assume that there is a unique physical measure $\mu$ (or by (i) equivalently that the space of joint resonant states at 0 is one dimensional). Then the following are equivalent:
\begin{itemize} 
 \item[{\rm(1)}] ~The only Ruelle-Taylor resonance on $i\a^*$ is zero.
 \item[{\rm(2)}] ~There exists $H\in \a$ such that $\varphi_t^{X_H}$ is \emph{weakly} mixing with respect to $\mu$.
 \item[{\rm(3)}] ~For any $H\in \mf a_+^\tau$, $\varphi_t^{X_H}$ is \emph{strongly} mixing with respect to $\mu$.
\end{itemize}
\item[{\rm(iv)}] \label{it:complex_measures} $\lambda\in i\a^*$ is a joint Ruelle resonance if and only if there is a complex measure $\mu_\lambda$ with 
${\rm WF}(\mu_\lambda)\subset E_\mathrm{s}^*$ satisfying for all $H\in \mf a_+^\tau, t\in \R$ the following equivariance under push-forwards of the action: $(\varphi^{X_{H}}_t)_*\mu_\lambda=e^{-\lambda(H)t}\mu_\lambda$. Moreover, such measures are absolutely continuous with respect to the physical measure obtained by taking $v=1$ in \eqref{eq:Cone-averaging-physical-measure}.
\item[{\rm(v)}] \label{it:res_at_zero} If $\mc M$ is connected and if there exists a smooth invariant measure $\mu$ with $\supp(\mu)=\mc M$, we have for any $j=0,\ldots,\ell$
\[
\dim \left(\Big(\ker d_{X}|_{C^{-\infty}_{E_\mathrm{u}^*}(\M)\otimes \Lambda^j\a^*_\C}\Big)\Big/
\Big(\ran\, d_{X}|_{C^{-\infty}_{E_\mathrm{u}^*}\otimes \Lambda^{j-1}\a^*_\C}\Big)\right)={\binom{\ell}{j}}.
\] 
\end{enumerate}
 \label{Theo3intro}
\end{thm}

The isomorphism stated in (i) and the existence of the complex measures in (iv) can  be constructed explicitly in terms of spectral projectors built from the parametrix \eqref{eq:parametrix_intro}. 

In the case of a single Anosov flow, physical measures are known to coincide with SRB-measures (Sinai-Ruelle-Bowen, see e.g. \cite{You02} and references therein). The latter are usually defined as invariant measures that can locally be disintegrated along the stable or unstable foliation of the flow with absolutely continuous conditional densities.

In \cite{BGW21} it is proved that the microlocal characterization Theorem~\ref{Theo3intro}(ii) of physical measures via their wavefront sets implies that the physical measures of an Anosov action are exactly those invariant measures that allow an absolutely continuous disintegration along the stable manifolds. Moreover, \cite[Theorem 3]{BGW21} says that for each physical/SRB measure, there is a basin $B\subset \mc{M}$ of positive Lebesgue measure such that for all $f\in C^0(\M)$, all proper open subcones $\mathcal{C}\subset\mf a_+^\tau$ and all $x\in B$, one has the convergence
\begin{equation}\label{eq:Cone-averaging-ergodic-SRB-measure}
\mu(f) =\lim_{T\to \infty} \frac{1}{{\rm Vol}(\mc{C}_T)} \int_{A\in \mc{C}_T}f(\varphi^{-X_A}_1(x)) \intd A.
\end{equation}
Moreover,  \cite[Theorem 2]{BGW21} says that the measure $\mu$ can be written as an infinite weighted sum of the Dirac measures on the periodic tori of the action, showing an equidistribution of periodic tori in the support of $\mu$.  
This measure has full support in $\mc{M}$ if the action is positively transitive in the sense that there is a dense orbit $\cup_{H\in \mf a_+^\tau} \varphi_1^{X_H}(x)$ for some $x\in \mc{M}$. The existence of such a measure is in fact considered as an important step towards the resolution of the Katok-Spatzier rigidity conjecture, \cite{KaSp94}.

\begin{rem}[Cocompact locally symmetric spaces and first band resonant states]
The theory of Ruelle-Taylor resonances sketched in this subsection applies to compact locally symmetric spaces $\mc X=\Gamma\bs G/K$ with the Weyl chamber flow, i.e. $A$-action from Remark~\ref{rem:A-action on G/M} on $\mc M=\Gamma\bs G/M$.

We denote the space $ \{u\in \mathcal D'_{E^*_\mathrm{u}}(\mc M)\mid 
\forall  H\in\mf a:\  (X_H+\lambda(H))u=0\}$ of resonant states by $\textup{Res}_X(\lambda)$. Then the space $\textup{Res}^0_X(\lambda)$ of \emph{first band resonant states} is defined as those resonant states that are in addition horocyclically invariant, 
\[
 \textup{Res}^0_X(\lambda):= \{u\in \mathcal D'_{E^*_\mathrm{u}}(\mc M)\mid 
\forall  H\in\mf a,  U\in C^\infty(\mc M; E_\mathrm{u}):\  (X_H+\lambda(H))u=0, U u= 0 \}
\]
We call a Ruelle-Taylor resonance a \emph{first band resonance} iff $\textup{Res}^0_X(\lambda)\neq 0$.

Working with horocycle operators and vector valued Ruelle-Taylor resonances it is possible to show that all resonances with real part in a certain neighborhood of zero in $\mf a^\ast$ are always first band resonances (see \cite[Prop.~3.7]{HWW21}). In view of the symplectic interpretation of the Weyl chamber flow, we consider the set $\sigma_\mathrm{RT}$ of Ruelle-Taylor resonances for that flow as a \emph{classical spectrum}.
\label{rem:first band resonant states}
\end{rem} 

\subsubsection{Ruelle resonances for noncompact locally hyperbolic spaces}\label{subsubsec:Ruelle Res for convec ccpt surfaces}

We state a special case of a result of Dyatlov and Guillarmou from \cite{DyGu} which yields Ruelle resonances and resonant states for convex cocompact hyperbolic surfaces.

\begin{thm}[\small {cf. \cite[Thm.~4.1]{GHW18}}]~ If $\X=\Gamma\backslash\mathbb{H}^{2}$ is a convex co-compact hyperbolic surface, then 
the generator $X$ of the geodesic flow on $\M=S\X$ has a resolvent 
$R_X(\lambda):=(-X-\lambda)^{-1}$ that admits a meromorphic extension from 
$\{\lambda\in \C\mid  {\rm Re}(\lambda)>0\}$ to 
$\C$ as a family of bounded operators $C_c^\infty(S\X)\to \mc{D}'(S\X)$. The 
resolvent $R_X(\lambda)$ has finite rank polar part at each pole $\lambda_0$ and the 
polar part is of the form 
\[ -\sum_{j=1}^{J(\lambda_0)} \frac{(-X-\lambda_0)^{j-1}\Pi^X_{\lambda_0}}{(\lambda-\lambda_0)^{j}}, \quad J(\lambda_0)\in \N \]
for some finite rank projector $\Pi^X_{\lambda_0}$ commuting with $X$. Moreover, 
$u\in \mc{D}'(S\X)$ is in the range 
of $\Pi^X_{\lambda_0}$ if and only if $(X+\lambda_0)^{J(\lambda_0)}u=0$ with $u$ supported in the outgoing tail
$\Lambda_+$ of the geodesic flow and ${\rm WF}(u)\subset E_\mathrm{u}^*$.  
\label{dygu}
\end{thm}

Now one can define \emph{Ruelle resonance},  \emph{generalized Ruelle resonant state} and 
\emph{Ruelle resonant state} as, respectively, a pole $\lambda_0$ of $R(\lambda)$, an 
element in ${\rm Im}(\Pi^X_{\lambda_0})$, and an element in 
${\rm Im}(\Pi^X_{\lambda_0})\cap \ker(-X-\lambda_0)$.
Define the spaces 
\begin{equation}\label{Resk2}
{\rm Res}^j_X(\lambda):=\{ u\in \mc{D}'(S\X)\mid {\rm supp}(u)\subset \Lambda_+, \, 
{\rm WF}(u)\subset E_\mathrm{u}^*,\,  (-X-\lambda)^ju=0\},
\end{equation}
\begin{equation}\label{Vm2}
V_m^j(\lambda):=\{ u\in {\rm Res}_X^j(\lambda)\mid  U_-^{m+1}u=0\}
\end{equation}
with $U_-:=\begin{pmatrix}
0&0\\ 1&0
\end{pmatrix}$.
As in the compact case, the operator $(-X-\lambda_0)$ is nilpotent on the finite dimensional space 
${\rm Res}_X(\lambda_0):=\cup_{j\geq 1}{\rm Res}^j_X(\lambda_0)$ and 
$\lambda_0$ is a Ruelle resonance if and only if ${\rm Res}_X^1(\lambda_0)\not=0$. The presence of Jordan blocks
for $\lambda_0$ is equivalent to having ${\rm Res}_X^k(\lambda_0)\not={\rm Res}_X^1(\lambda_0) $ for some $k>1$.

\begin{rem}[Asymptotically hyperbolic spaces]~ In \cite{Ha17}
Hadfield  constructed Ruelle resonances for negatively curved manifolds which asymptotically at infinity behave as real hyperbolic spaces. This includes the class of convex cocompact quotients of hyperbolic space.
\label{rem:Ruelle resonance - Hadfield}
\end{rem}

\begin{rem}[Hyperbolic manifolds with cusps]~
Bonthonneau and Weich in \cite{BW22} constructed Ruelle resonances for Riemannian manifolds with finite volume, negative curvature and the unbounded part consisting of finitely many real hyperbolic cusps.
\label{rem:Ruelle resonance - Bonthonneau-Weich}
\end{rem}

\subsubsection{$\Gamma$-invariant distributions for principal series representations}

Consider compact locally symmetric spaces and recall the first band resonant states from Remark~\ref{rem:first band resonant states}. It turns out that these are closely related to $\Gamma$-invariant distribution vectors in principal series representations, Remark~\ref{rem:extension of G-representations}. In rank one
the following proposition was established in \cite{GHW21}.

\begin{prop}[{\small \cite[Prop.~3.8]{GHW21}}] 
 There is an isomorphism of finite dimensional vector spaces
 \[
  \tu{Res}_{X}^0(\lambda) \cong {^\Gamma(H_{-(\lambda+\rho)}^{-\infty})},
\]
where ${^\Gamma(H_{-(\lambda+\rho)}^{-\infty})}$ denotes the spaces of $\Gamma$-invariant distributional vectors in the spherical principal series with spectral parameter $\mu=-(\lambda+\rho)$.
\label{prop:ruelle_gamma_inv_distr}
\end{prop}

In \cite{HWW21} it was shown that this result extends to higher rank. In fact, using  lifts to vector bundles (as they were discussed in \cite{GGHW21}) one even gets a version for nonspherical principal series representations, see \cite[Lemma~3.10]{HWW21} for a precise formulation.

\subsubsection{Dynamical zeta functions}\label{subsec:zeta}

Dynamical zeta functions are generating functions counting dynamical objects such as closed orbits.   The dynamical zeta functions associated with the geodesic flow count closed geodesics using different weight functions.

\begin{enumerate}
\item[(i)]~\emph{Ruelle zeta function}:
$$Z_{\mathrm R}(\lambda)=\prod_{c}(1-e^{-\lambda \ell(c)})^{\pm1}= \exp\Big(\pm \sum_c\sum_{k\ge 1} \frac{e^{-k \lambda \ell(c)}}{k}\Big),$$ 
where the first summation is over the primitive closed geodesics $c$ and $\ell(c)$ is the length of $c$.
\item[(ii)]~\emph{Selberg zeta function}:
\begin{eqnarray*}
Z_{\mathrm S}(\lambda)
&=&\prod_{c}\prod_{k\ge 0}\det\Big(1- S^k(P_c^+)e^{-(\lambda+\rho) \ell(c)})^{\pm1}\Big)\\
&=& \exp\Big(-\sum_c\sum_{j\ge 1} \frac{e^{-j (\lambda+\rho) \ell(c)}}{j\det(1-(P_c^+)^j)}\Big),
\end{eqnarray*}
where $P_c=P_c^+\oplus \mathrm{id}\oplus P_c^-\in \mathrm{End} (T_v(\Gamma\backslash G/M))$ is the monodromy (i.e. Poincar\'e) map and $S^k$ denotes the $k$-th symmetric power.
\item[(iii)]~\emph{Dynamical zeta function}:
$$Z_{\mathrm{dyn}}(\lambda)= \exp\Big(-\sum_c\sum_{j\ge 1} \frac{e^{-j \lambda \ell(c)}}{j\big|\det(1-(P_c^+)^j)\big|}\Big).$$
\end{enumerate}

There are more versions with additional twists and weights (see e.g. \cite[Chap.~3]{BO95} and  \cite{AZ}). For the relation between dynamical and Selberg zeta functions we refer to \cite{FT17}. A key feature of such zeta functions is that they often admit meromorphic continuation (see e.g. \cite{GLP13,DyZw} and the literature cited there). The resulting divisors may serve as a classical spectral invariant of the underlying flows.

\subsubsection{Transfer operators}\label{subsuubsec: transfer operators}

Transfer operators were introduced in statistical mechanics as a means to find equilibrium distributions. In our context they were used primarily to study symbolic dynamical systems derived from the geodesic flow on hyperbolic surfaces. As is observed in \cite{PZ20}, for suitable symbolic systems they do resemble weighted Laplacians on graphs which gives them a double meaning, one classical, one quantum mechanical. For the modular surface this eventually lead to the first dynamical interpretation of Maass cusp forms in \cite{LZ01}. We refer to the introduction in \cite{Po1} for the extended history of this development.

In general, for a discrete dynamical system $(\mc N,f)$ with $f: \mc N\to \mc N$ any map, one associates a \emph{transfer operator} $\mathcal L_\psi$ on suitable spaces function spaces on $\mc N$ to a \emph{weight map} $\psi:\mc N\to \C$  by
\begin{equation}\label{eq:transfer operator}
\mathcal L_\psi\phi(x)= \sum_{y\in f^{-1}(x)} e^{\psi(y)} \phi(y).
\end{equation}
Here ``suitable‘‘ depends on the convergence of the sum in \eqref{eq:transfer operator}.

For hyperbolic surfaces the construction of the transfer operator typically is a two-step process. First one constructs a discrete dynamical system with invertible $f$ whose suspension flow (see \cite[\S~0.3]{KH95}) is basically conjugate to the geodesic flow. Then one applies a reduction leading to a discrete dynamical system for which $f$ is no longer invertible. In \cite[\S~2.2]{CM99} the process for the modular surface goes as follows. They quote the symbolic dynamics for the geodesic flow from earlier work and start with
\begin{eqnarray*}
\tilde f: ]0,1[\times]0,1[\times \{\pm1\} &\to& ]0,1[\times]0,1[\times \{\pm1\}\\
 (a,b,\epsilon)&\mapsto& \big(a^{-1}\,\mathrm{mod}\,1, (b+\lfloor a^{-1}\rfloor)^{-1},-\epsilon\big).
\end{eqnarray*}
and then work with the reduced dynamical system given by 
\[f: ]0,1[\to ]0,1[,\quad x\mapsto x^{-1}\,\mathrm{mod}\,1.\]
As weight functions they use $x\mapsto x^{2\beta}$ with $\Re \beta >\frac{1}{2}$. The constructions in \cite{Po1} are much more general, but also more involved. While it is quite clear that the reduction in \cite{CM99} means that the authors concentrate on the expanding part of the dynamics, it is less transparent how to relate the reduction in \cite[\S~9]{Po1} to the hyperbolic nature of the geodesic flow.

In the transfer operator approach to quantum-classical correspondences it is the spectral theory of the transfer operators which take the role of Ruelle resonances and resonant states. It should be noted that while the scope of locally symmetric spaces that can be dealt with using the transfer operator approach is almost exclusively restricted to surfaces (see, however, \cite{Po20} for some higher rank examples), it works well in many noncompact situations.

The analogs of geodesic and Weyl chamber flows for graphs and higher rank non-Archimedean counterparts of locally symmetric spaces, i.e. quotients of affine buildings, are of a combinatorial nature. So it is not surprising that transfer operators seem to be suitable tools to establish quantum-classical correspondences in those contexts. In the long run, if one wants to deal with adelic situations, microlocal and transfer methods may have to be used simultaneously in a coordinated way.

For yet another approach to quantum-classical correspondences via transfer operators  we refer to the book \cite{FT15} by Faure and Tsujii.

\subsection{Correspondences}\label{subsec:QCC}

In this section we finally come back to the quantum-classical correspondences for locally symmetric spaces that were alluded to in the motivating examples. It turns out that there are partial results of varying generality depending on the specifications of the locally symmetric spaces considered. The most general results are available for compact locally symmetric spaces. We start with that case and then describe what can be done at the moment for noncompact locally symmetric spaces.

\subsubsection{The compact case for generic spectral parameters}\label{subsubsec:HWW21}

Recall the right $A$-action on $T\M$ from Remark~\ref{rem:A-action on G/M} and assume that  $\Gamma$ is cocompact. Then the $A$-action is Anosov as we remarked in Subsection~\ref{subsec: Anosov flows}. For $H\in \mf a $ let $X_H$ again be the vector field on $\mathcal M$ defined by the right $A$-action. As was explained in Subsection~\ref{subsubsec: RT-resonances} the \emph{Ruelle-Taylor resonances} of this Anosov action are then given by
\[
\sigma_{\RT} = \{\lambda \in \mf a^*_\C \mid \exists u\in \mc D'_{E_\mathrm{u}^\ast}(\mathcal M)\setminus \{0\}\,\forall H\in\mf a : (X_H + \lambda(H))u=0\},
\]
where $\mc D'_{E_\mathrm{u}^\ast}(\M)$ is the set of distributions with wavefront set contained in the annihilator $E_\mathrm{u}^\ast\subseteq T^\ast \mc M$ of $E_0\oplus E_\mathrm{u}$. The distributions $ u\in \mc D'_{E_\mathrm{u}^\ast}(\M)$ satisfying $ (X_H + \lambda(H))u=0$ for all $H\in\mf a$ are the resonant states of $\lambda$ and the dimension of the space of all such distributions is called the \emph{multiplicity} $m(\lambda)$ of the resonance $\lambda$. According to Theorem~\ref{Theo2intro}, $\sigma_{\RT} \subset \mf a_\C^*$ is discrete and all resonances have finite multiplicity. Moreover, the real part of the resonances are located in the negative dual cone $-{\mf a}_+^\star\subset \mf a^*$ of the positive Weyl chamber  $\mf a_+$.

We have the following correspondence between the space of classical first band resonant states and the joint quantum eigenspace (see Subsection~\ref{subsubsec: quantum spectrum cocompact}):

\begin{thm}[\small{\cite[Thm.~1.3]{HWW21}}] 
~Let $\lambda\in \mf a_\C^\ast$  be outside the \emph{exceptional set} $\mc E\coloneqq \{\lambda\in \mf a_\C^\ast\mid \exists \alpha\in \Sigma^+:\ \frac {2\langle \lambda+\rho,\alpha\rangle}{\langle\alpha,\alpha\rangle} \in -\N_{>0}\}$. Then there is a bijection between the finite dimensional vector spaces
 \[
  \pi_\ast : \textup{Res}^0_X(\lambda) \to {}^\G E_{-\lambda-\rho},
 \]
where $\pi_\ast$ is the push-forward of distributions along the canonical projection $\pi:\G\backslash G/M\to \G\backslash G/K$.
\label{thm:qcc}
\end{thm}

Using the 1:1-correspondence in Theorem~\ref{thm:qcc} and results about the quantum spectrum, one can obtain obstructions and existence results on the Ruelle-Taylor resonances. Notably,  in \cite{HWW21}  we used results of Duistermaat-Kolk-Varadarajan \cite{DKV79} on the spectrum $\sigma_Q$ but we also use refined information on the quantum spectrum such as   $L^p$-bounds for spherical functions obtained from asymptotic expansions \cite{vdBanSchl87}  and $L^p$-bounds for matrix coefficients based on work by Cowling and Oh \cite{cowling, oh2002}.
Theorem~\ref{thm:Weyl} and Theorem~\ref{thm:gap} as stated below give only a rough version of the information on the Ruelle-Taylor resonances that we can actually obtain. As the full results require some further notation we refrain from stating them here, but rather refer to \cite[Thm.~5.1]{HWW21}.  

The first application says that for any Weyl chamber flow there exist infinitely many Ruelle-Taylor resonances by providing a Weyl lower bound on an appropriate counting function.

\begin{thm}[\small{\cite[Thm.~1.1]{HWW21}}]  
~Let $\rho$ be the half-sum of the positive restricted roots, $W$ the Weyl group, and for $t>0$ let
 \[
  N(t):=\sum_{\lambda \in \sigma_{\RT}, \Re(\lambda)=-\rho, \|\Im(\lambda)\|\leq t}m(\lambda).
 \]
 Then for $d:=\textup{dim}(G/K)$ 
\[
 N(t) \geq |W|\textup{Vol}(\G\backslash G/K) \left(2\sqrt{\pi}\right)^{-d}\frac{1}{\Gamma(d/2+1)}t^d+ \mathcal O(t^{d-1}).
\]
More generally, let $\Omega \subseteq \mf a^\ast$ be open and bounded such that $\partial \Omega$ has finite $(n-1)$-dimensional Hausdorff measure. Then 
\[\sum_{\lambda \in \sigma_{\RT}, \Re(\lambda)=-\rho, \Im(\lambda)\in t\Omega} m(\lambda) \geq  |W|\textup{Vol}(\G\backslash G/K) \left(2\pi \right)^{-d} \textup{Vol}(\Ad(K)\Omega) t^d +\mc O(t^{d-1}).
\]
\label{thm:Weyl}
\end{thm}

The second application is the existence of  a uniform spectral gap.

\begin{thm}[\small {\cite[Thm.~1.2]{HWW21}}] 
Let $G$ be a real semisimple Lie group with finite center, then for any cocompact torsion-free discrete subgroup $\G\subset G$ there is a neighborhood $\mc G\subset\mf a^*$  of $0$ such that 
\[
 \sigma_{\RT} \cap (\mc G \times i\mf a^*) = \{0\}.
\]
If $G$ furthermore has Kazhdan's property (T) (e.g. if $G$ is simple of higher rank), then  the spectral gap $\mc G$ can be taken uniformly in $\Gamma$ and only depends on the group $G$.
\label{thm:gap}
\end{thm}

This theorem has an extension to vector bundles based on the lifted version of Ruelle-Taylor resonances discussed in \cite{GGHW21}, see \cite[Thm.~5.1]{HWW21}. In the following remark we comment on the reasons why such extensions are relevant.

\begin{rem}[Homogeneous vector bundles]~
One of the key features setting hyperbolic surfaces apart from higher dimensional, let alone higher rank, locally symmetric spaces is that the sphere bundle $\mathrm{SL}_2(\mathbb R)/\{\pm\mathbf 1\}$ of its simply connected covering $\mathbb H=\mathrm{SL}_2(\mathbb R)/\mathrm{SO}_2$ is actually a group. Thus compared to the general rank one case, where the sphere bundle is $G/M$, it has extra symmetries, which simplify the analysis in many ways. In particular, the unipotent Iwasawa group $N$ acts from the right on $G/M$ giving the horocyclic flow which in general does not live on the sphere bundle of the locally symmetric space but only on $\Gamma\backslash G$.

The  tensor bundles relevant for the differential analysis on $G/M$ are $G$-homogeneous bundles associated with $M$-representations. Their sections can be described as functions on $G$ with $M$-equivariance properties. This suggests to view $G$ as a principal bundle over $G/M$ with fiber $M$ and analyze the right $G$-action in terms of its effects on the tensor algebra. This is what was done in \cite{DFG15}, where the authors determined the band structure of the Ruelle resonances for compact hyperbolic spaces and in \cite{Em14}, where Emonds worked out a dynamical interpretation of Patterson-Sullivan distributions  generalizing observations from \cite{AZ} for hyperbolic surfaces. Similarly, it was used in the study \cite{AH23} of quantum-classical correspondences for exceptional spectral parameters.  It is also implicit in the work of Bunke and Olbrich. 

For compact locally symmetric spaces of rank one Küster and Weich in \cite{KW21} undertook a systematic study of Ruelle resonant states consisting of distributions with values in $M$-representations leading up to a quantum-classical correspondence with eigenspaces of the Bochner Laplacian on $G$-homogeneous vector bundles over $K$ associated with $K$-representations intertwining nontrivially with the $M$-representation used for the Ruelle resonances. In particular, they define horocyclic operators generalizing the ones from \cite{DFG15}, thus opening a path to clarify the precise band structure of Ruelle resonances for general compact locally symmetric spaces of rank one. 

In the companion paper \cite{KW20} Küster and Weich apply special cases of their results to write the first Betti number of negatively curved compact Riemannian manifolds obtaind by deformation from hyperbolic ones in terms of Ruelle resonances.
\label{rem:vector bundles}
\end{rem}

\subsubsection{Results for noncompact locally symmetric spaces}\label{subsubsec:QCC-noncompact}

The first quantum-classical correspondences involving quantum resonant states were the ones on convex cocompact hyperbolic spaces alluded to in Subsection~\ref{subsubsec:Resolvent poles}. We state explicitly the correspondence established in \cite{GHW18} for hyperbolic surfaces and refer to  \cite{Ha17,Ha20} for extensions to higher dimension.

Let the manifold $\X=\Gamma\backslash \mathbb{H}^{2}$ be a noncompact complete smooth hyperbolic suface with infinite volume but finitely many topological ends. 
$\X$ can be compactified to the smooth manifold 
$\bbar{\X}=\Gamma\backslash (\mathbb{H}^2\cup \Omega_\Gamma)$, where
$\Omega_\Gamma\subset \mathbb{S}^1$ is the \emph{set of discontinuity} of $\Gamma$, the complement of  the limit set  $\Lambda_\Gamma\subset \mathbb{S}^1$  of $\Gamma$.
Note that $\X$ is \emph{conformally compact} in the sense of Mazzeo-Melrose \cite{MM87}: 
there is a smooth boundary defining function $r$ such that $\bar{g}:=r^2g$ 
extends as a smooth metric on $\bbar{\X}$. Here $g$ denotes the Riemannian metric on $\mc X$.

In this case the group $\Gamma$ is a subgroup of ${\rm PSL}_2(\mathbb{C})$ and acts on the Riemann sphere  $\overline{\mathbb{C}} := \mathbb{C}\cup\{\infty\}$ as conformal  transformations, it preserves the unit disk $\mathbb{H}^2$ and its complement $\overline{\mathbb{C}}\setminus \mathbb{H}^2$.  Equivalently, by conjugating by $(z-i)/(z+i)$, $\Gamma$ acts by conformal transformations on $\overline{\mathbb{C}}$ as a subgroup of ${\rm PSL}_2(\mathbb{R})\subset {\rm PSL}_2(\mathbb{C})$ and it preserves  the half-planes $\mathbb{H}^2_\pm:=\{z\in \mathbb{C}\mid  \pm{\rm Im}(z)>0\}$. The half-planes are conformally equivalent through $z\mapsto \bar{z}$ if we put opposite orientations on $\mathbb{H}^2_+$ and $\mathbb{H}^2_-$. In this model the boundary is the compactified real line $\partial \mathbb{H}_\pm= \overline{ \mathbb{R}} := \mathbb{R}\cup\{\infty\}$ and the limit set is a closed subset $\Lambda_\Gamma$ of $\overline{\mathbb{R}}$. Its complement in $\overline{\mathbb{R}}$ is still denoted by $\Omega_\Gamma$. Since $\bbar{\gamma(\bar{z})}=\gamma(z)$ for each $\gamma\in \Gamma$, the quotients $\X_\pm:=\Gamma\backslash (\mathbb{H}^2_\pm\cup \Omega_\Gamma)$ are smooth surfaces with boundaries, equipped with a natural conformal structures and $\X_+$ is conformally equivalent to $\X_-$. The surface $\Gamma\backslash (\overline{\mathbb{C}}\setminus \Lambda_\Gamma)$ is a compact surface diffeomorphic to the gluing  $\X_2:=\X_+\cup \X_-$ of $\X_+$ and $\X_-$ along their boundaries, moreover it is  equipped with a smooth conformal structure which restricts to that of $\X_\pm$. 

We denote by $\mc{I}:\X_2\to \X_2$ the involution fixing $\partial\bbar{\X}$ and derived from $z\mapsto \bar{z}$ when viewing $\Gamma$  as acting in $\overline{\mathbb{C}}\setminus \Lambda_\Gamma$.   The interiors of $\X_+$ and $\X_-$ are isometric if we put the hyperbolic metric  $|dz|^2/({\rm Im}(z))^2$ on $\mathbb{H}^2_\pm$, and they are isometric to the hyperbolic surface $\X$. The conformal class of $\X_\pm$ corresponds to the  conformal class of $\bar{g}$ on $\bbar{\X}$ as defined above. We identify $\X_+$  with $\bbar{\X}$ and define $H_{\pm n}(\X)$ as the following finite dimensional real vector spaces 
\begin{equation}\label{HnH-n}
\begin{gathered}
 H_n(\X):=\{ f|_{\X_+}\mid f\in C^\infty(\X_2; \mc{K}^n),\, \bbar{\partial}f=0, \, \mc{I}^*f=\bar{f}\}, \\
 H_{-n}(\X):=\{ f|_{\X_+}\mid f\in C^\infty(\X_2; \mc{K}^{-n}),\, \partial f=0,\, \mc{I}^*f=\bar{f}\},
\end{gathered} 
\end{equation}
where $ \mc{K}$ is the canonical line bundle, see \cite[\S~1.1]{GH78}.

With the notation established in Subsection~\ref{subsubsec:Ruelle Res for convec ccpt surfaces} we have the following  quantum-classical correspondence. 

\begin{thm}[\small cf. {\cite[Thm.~4.2]{GHW18}}] ~
Let $\X=\Gamma\backslash \mathbb{H}^{2}$ be a smooth oriented convex cocompact hyperbolic 
surface and let $\M=S\X$ be its unit tangent bundle.
\begin{enumerate}
\item[{\rm (i)}] For each $\lambda\in \mathbb C\setminus (-\demi-\demi \mathbb N_0)$  the pushforward map ${\pi}_*: \mc{D}'(S\X)\to \mc{D}'(\X)$ for the projection $\pi:\mc M\to\mc X$ restricts to a linear isomorphism of complex vector spaces for each $j\geq 1$
\begin{equation}\label{pi*iso}
{\pi_0}_* : V_0^j(\lambda) \to {\rm Res}_{\Delta_\X}^j(\lambda_0+1),
\end{equation}
where $\Delta_\X$ is the Laplacian on $\X$ acting on functions.
\item[{\rm (ii)}] For each  $\lambda=-\demi-k$ with $k\in\mathbb N$, $V_0^j(\lambda)=0$ and ${\rm Res}^j_{\Delta_\X}(\lambda+1)=0$ 
for all $j\in\N$.
\item[{\rm (iii)}] For $\lambda=-\demi$, there are no Jordan blocks, i.e. $V_0^j(-\demi)=0$ for $j>1$, 
and the map 
\begin{equation}\label{pi*iso2} 
{\pi_0}_* : V_0^1(-\demi)\to {\rm Res}_{\Delta_\X}^1(1/2)
\end{equation}
is a linear isomorphism of complex vector spaces.
\item[{\rm (iv)}]  For $\lambda=-n\in -\mathbb N$, if $\Gamma$ is non-elementary (i.e. there are no invariant sets of cardinality $1$ or $2$), there are no Jordan blocks, 
i.e. $V_0^j(-n)=0$ if $j>1$, and the following map is an isomorphism of real vector spaces
\begin{equation}\label{pi*iso3}
 i^{n+1}{\pi_n}_* : V_0^1(-n)\to H_{n}(\X)
 \end{equation}
where $H_n(\X)$ is defined by \eqref{HnH-n}. 
\end{enumerate} 
\label{classicquantic}
\end{thm}

\begin{rem}[QCC via transfer operators] 
The Lewis-Zagier correspondence from Example~\ref{ex:Lewis-Zagier} is an example of a quantum-classical correspondence established via transfer operator methods. It has been generalized in various directions. As these methods are not in the focus of this article, I do not present further details here, but refer to the introduction of \cite{Po1} for a brief discussion of further developments as well as references. 
\end{rem}

\begin{rem}[Cuspidal $\Gamma$-cohomology]~
Theorem~\ref{thm:DH05-0.2}  may be viewed as a weak quantum-classical correspondence since it proves for irreducible principal series representations $\pi$ the multiplicity space $\mathrm{Hom}_G(\pi,L^2_\mathrm{cusp}(\Gamma\bs G))$ is isomorphic to $H_\mathrm{cusp}^{r}(\Ga,\pi^{-\omega})$. However, the isomorphism is only guaranteed by the equality of dimensions and not described in any natural way.
\label{rem:cusp-cohom}
\end{rem}

\subsubsection{Exceptional spectral parameters}\label{subsubsec:exceptional spectral parameters}

In the quantum-classical correspondences described in Sections~\ref{subsubsec:HWW21} and \ref{subsubsec:QCC-noncompact} one always had to make restrictions on the spectral parameters describing classical resonances. The reason is that the method of proof in all these examples depends on (scalar) Poisson transforms which for generic parameters are bijective. For the exceptional parameters leading to Poisson transforms which are not bijective the proofs break down and one has to find a replacement for the scalar Poisson transforms. 

The strategy for an extension of the quantum-classical correspondence to exceptional spectral parameters that has been succesfully applied to compact rank one locally symmetric spaces in \cite{AH23} is as follows. As in the generic case (see \cite[\S~3.2]{GHW21}) we start by lifting the first band Ruelle resonances to $\Gamma$-invariant distributions on the global symmetric space. The lifted spaces can be interpreted in terms of spherical principal series (that part works for all spectral parameters, see \cite[Prop.~3.8]{GHW21}) and the first band resonant states $\mathrm{Res}_X^0(-\lambda-\rho)$ correspond to the space $^\Gamma H_\lambda^{-\infty}$ of $\Gamma$-invariant distribution vectors of the corresponding principal series. For an exceptional spectral parameter $\lambda$ the corresponding principal series $H_\lambda$ is no longer irreducible. But it has a managable composition series and it turns out that the $\Gamma$-invariant distribution vectors are all contained in the \emph{socle} (i.e. the sum of all irreducible subrepresentations) of the representation. In each of the rank one cases except $\mathrm{SO}_0(2,1)$ (the case of surfaces, see \cite{GHW18}) the socle turns out to be irreducible with a unique minimal $K$-type $\tau_\lambda$ %(see Thm.~\ref{thm:socle}) 
and one can show that the vector valued Poisson transform associated with this $K$-type (sum of $K$-types in the case of surfaces) is injective. %, see Prop.~\ref{prop:Poisson_sum}. 
The image consists of spaces of $\Gamma$-invariant sections of vector bundles over $\Gamma\backslash G/K$ and we have a quantum-classical correspondence as soon as we have characterized the image of this Poisson transform.

\begin{thm}[{\small\cite[Thm.~6.1]{AH23}}]~
Consider the case of a cocompact rank one locally symmetric space with $G=\mathrm{SO}_0(n,1)$, $n\ge3$. Given an exceptional parameter $\lambda=-(\rhoa+\ell\alpha),\ \ell\in\N_0,$ the socle $\on{soc}(\pst{\lambda})$ of $\pst{\lambda}$ is irreducible, unitary and its $K$-types are given by the spaces $Y_k$ of harmonic homogeneous polynomials of degree $k$ for $k\geq\ell+1$. The minimal $K$-type is $Y_{\ell+1}$ and the corresponding Poisson transform induces an isomorphism
\begin{gather*}
P_\lambda^{Y_{\ell+1}}\colon{}^\Gamma(\on{soc}(\pst{\lambda}))^{-\infty}\overset{\cong}{\longrightarrow}{}^\Gamma\{f\in C^\infty(G\times_KY_{\ell+1})\mid\mathrm{d}_-f=0,\ \mathrm{D}f=0\},
\end{gather*}
with  differential operators $\mathrm{d}_-$ and $\mathrm{D}$ which can be described explicitly.
\label{thm:spectral_corres_SOn}
\end{thm}

From the study of the zeros of dynamical zeta functions one expects the multiplicities of exceptional resonances to carry topological information. In \cite{GHW18} this expectation was confirmed for hyperbolic surfaces defined by cocompact and convex cocompact discrete subgroups of $\mathrm{PSL}_2(\mathbb R)$. 

\begin{ex}[Topological interpretation {\small\cite{GHW18}}]~Let $\Gamma\backslash G/K$ be a compact hyperbolic surface. Then
\[\dim\mathrm{Res}^0_X(-n)=\begin{cases}
|\chi(\Gamma\backslash G/K)|+2&\text{if } n=1, \\
(2n-1)|\chi(\Gamma\backslash G/K)|&\text{if } 1<n\in\mathbb N, 
\end{cases}\]
where $\chi(\Gamma\backslash G/K)$ is the Euler characteristik of $\Gamma\backslash G/K$.
\end{ex}

In higher dimension we do not have a topological interpretation of the resonance multiplicities except in very special examples boiling down to Hodge theory.

\begin{ex}[Topological interpretation {\small \cite{AH23}}]~
In the situation of Theorem~\ref{thm:spectral_corres_SOn} for the first exceptional parameter $\lambda=-\rhoa$ we get ($Y_1\cong\mf p^*$)
\begin{gather*}
P_{-\rhoa}^{Y_{1}}\colon{}^\Gamma(\on{soc}(\pst{-\rhoa}))^{-\infty}\overset{\cong}{\longrightarrow}\{f\in C^\infty(\Lambda^1(\Gamma\backslash G/K))\mid\delta f=0,\ d{ }f=0\},
\end{gather*}
where $\Lambda^1(\Gamma\backslash G/K)$ denotes the bundle of one forms and ($\delta$ resp.\@) $d$ is the (co)-differential. The dimension is given by the first Betti number $b_1(\Gamma\backslash G/K)$ in this case.
\end{ex}

The quantum-classical correspondences for exceptional spectral parameters get more complicated for the other compact rank one locally symmetric spaces. This is due to the fact that the system of differential equations  characterizing the image of the vector valued Poisson transform is more involved.  

\begin{rem}[Cocompact locally symmetric spaces of rank one]
Results analogous to Theorem~\ref{thm:spectral_corres_SOn} are available for all the remaining cocompact rank one locally symmetric spaces, i.e. for the cases
\begin{enumerate}
\item[(i)]~ $G=\mathrm{SU}(n,1),\ n\geq 2$, see \cite[Theorems~6.4 \& 6.6]{AH23},
\item[(ii)]~ $G=\mathrm{Sp}(n,1),\ n\geq 2$, see \cite[Thm.~6.8]{AH23},
\item[(iii)]~ $G=\mathrm{F}_{4(-20)}$, see \cite[Thm.~6.10]{AH23}.
\end{enumerate}
In the symplectic and the exceptional cases, however, the description of the image of the vector valued Poisson transformation is a little more complicated.
\label{rem:QCC rank 1}
\end{rem}

The representations showing up as socles of spherical principal series representation for exceptional parameters can be determined explicitly in terms of their Langlands parameters. In fact, they have a uniform description as was observed by Jan Frahm when he saw the list of Langlands parameters.

\begin{thm}[{\small\cite[Thm.~4.7]{AH23}}]~
There is a one-to-one correspondence between the representations $\on{soc}(\pst{\lambda_\ell})$ with $\lambda_\ell$ exceptional and the relative discrete series of the associated pseudo-Riemannian symmetric spaces  $G/H$. More precisely, each of the representations corresponds to a minimal closed invariant subspace of $L^2(G/H)$ with $H=\mathrm{SO}_0(n-1,1)$, $\mathrm{S}(\mathrm{U}(1)\times\mathrm{U}(n-1,1))\cong \mathrm{U}(n-1,1)$ ,  $\mathrm{Sp}(1)\times\mathrm{Sp}(n-1,1)$, or  $\mathrm{Spin}(1,8)$.
\label{thm:discrete series}
\end{thm}

\subsubsection{Applications}\label{subsubsec:qcc-applications}

Most applications of quantum-classical correspondences as described in this section that are on record so far, make use of the possibility to compute ``classical'' multiplicities from ``quantum''' multiplicites. First examples of such applications were already given in Subsection~\ref{subsubsec:HWW21} (see \cite[Theorems~1.1\ \&\ 1.2]{HWW21}). 

In \cite{DGRS20}, which gives a proof of the Fried conjecture (in dimension $3$) relating analytic torsions to twisted dynamical zeta-values, the authors need to exclude $0$ as a Ruelle resonance for the geodesic flow in compact real hyperbolic manifolds of dimension $3$. This allows them to perturb the generator of the geodesic flow and thus prove Fried's conjecture for some $3$-manifolds with variable curvature. The proof given in \cite[Prop.~7.7]{DGRS20} for the fact that $0$ is not a Ruelle resonance uses the Selberg trace formula, but the authors point out in  \cite[Rem.~10]{DGRS20} that they could also use \cite{DFG15,KW21}.

Extending ideas from \cite{KW20} on the deformation of geodesic flows (see Section~\ref{subsubsec:HWW21}) Ceki\'c et al. in  \cite{CDDP22} showed that the order of vanishing at zero of the Ruelle zeta function is not stable under generic  deformations and thus provided counterexamples to Fried's conjecture for compact hyperbolic $3$-manifolds with nonvanishing first Betti number. In this work the quantum-classical correspondence is used to calculate  the vanishing order (which can be expressed in terms of multiplicities of Ruelle resonances) for the geodesic flow which is being deformed.

Finally, we mention an application due to Schütte, Weich and Barkhofen who in \cite{SWB23} introduce weighted zeta functions for flows which are hyperbolic on their trapped set and prove that they continue meromorphically to all of $\mathbb C$. They show show that poles of the zeta functions are Ruelle resonances of the flows and describe the residues of the zeta functions in terms of the residue operators of the Ruelle resonance (cf. Example~\ref{ex:DFG}) associated with the resolvent of the generating vector fields. Combining these residue formulas with \cite[Cor.~6.1]{GHW21} which expresses the \emph{flat trace} of residue operators in terms of Patterson-Sullivan distributions, they are able to extend the scope of a result of Anantharaman-Zelditch \cite{AZ} from compact hyperbolic surfaces  to all compact Riemannian locally symmetric spaces of rank one (see \cite[Thm.~4.1]{SWB23}). Note that \cite[Cor.~6.1]{GHW21} is a corollary of the pairing formula \cite[Thm.~6.1]{GHW21} relating the product of a resonant and a coresonant state of a Ruelle resonance to the pairing of their push-forwards along the projection to the base space (cf. Theorem~\ref{thm:qcc}). 

Actually, versions of the aforementioned results of from \cite{AZ} and \cite{SWB23} were formulated by physicists long ago (see \cite{EFMW92}). Arguing with semiclassical trace formulae they predicted residues of weighted zeta functions to be given by quantum phase space distributions in a semiclassical limit. In the companion paper \cite{BSW22} to \cite{SWB23} the same authors explain how their  \cite[Thm.~4.1]{SWB23} is a mathematically rigorous version of this prediction.

\section{Quotients of trees and affine buildings}

Trees and affine buildings come up in the context of $p$-adic Lie groups where they play a role analogous to symmetric spaces for real Lie groups (see e.g. \cite{BT72,Ca73,FTN91,Se80,Pa06a}). Homogeneous trees correspond to rank one situations, the affine buildings may be viewed as suitable higher rank versions of homogeneous trees. The quotients of trees and buildings then correspond to the locally symmetric spaces which are in the focus of this article. 

It should be noted that quotients of trees amount to connected graphs which are used very often as approximations or toy models for manifolds. This is in particular the case in dynamical systems (see e.g. \cite{LP16,An17}) but also in microlocal analysis (see e.g. \cite{ALM15,LM14}). It is thus reasonable to expect that meaningful results about buildings hold, which are analogous to the quantum-classical correspondences we have discussed in Section~\ref{subsec:QCC}. At the same time one can expect that the necessary analysis is less involved compared to the real case. One can hope for clues how to proceed in the real case in places where progress has been impeded by technical difficulties. Moreover, recent results (see \cite{BHW22,BHW23,AFH23a}) support the expectation that one can weaken the symmetry conditions and still obtain nontrivial correspondences.

\subsection{Graphs of bounded degree}

A natural starting point for research in this direction are homogeneous trees. For these, key tools like the Poisson transform have been available for a long time, see e.g. \cite{FTN91}. It turns out that in the case of trees one can weaken the symmetry assumptions quite a bit and still get Poisson transforms as well as quantum-classical correspondences.

We consider graphs $\mf G=(\mf X,\mf E)$ consisting of a set $\mf X$ of vertices and a set $\mf E\subseteq \mf X^2$ of directed edges.  We assume that the graph is  \emph{symmetric}, i.e.\@ the set $\mf E\subseteq \mf X^2$ is invariant under the switch of coordinates. Further we assume that the graph contains \emph{no loops}, i.e.\@ $(x,x)\not\in\mf E$ for all $x\in\mf X$, and has \emph{no dead ends}, which means that each vertex has at least two neighbors.

For each directed edge $\vec{e}=(a,b)$ we call $a=\iota(\vec{e})$ the initial and $b=\tau(\vec{e})$ the terminal point of $\vec{e}$. A \emph{path} in $\mf G$ is a (finite or infinite) sequence $(\vec{e}_j)$ of edges such that $\tau(\vec{e}_j)=\iota(\vec{e}_{j+1})$, but $\iota(\vec{e}_j)\not=\tau(\vec{e}_{j+1})$ (no backtracking!). We assume that any two vertices in $\mf X$ can be connected by a finite path of edges. 

As an analog of geodesic rays in Riemannian manifolds we introduce the space 
\[\mf P_+:= \{(\vec{e}_j)_{j\in \mathbb N}\}\]
of infinite paths starting at some vertex. We refer to such paths as \emph{geodesic rays} in $\mf G$. The assumption that $\mf G$ has no dead ends together with the symmetry implies that any finite path can be extended to a geodesic ray. This gives a kind of geodesic completeness of the graph. 

The shift  $(\vec e_1,\vec e_2,\ldots) \mapsto (\vec e_2,\vec e_3,\ldots)$ of geodesic rays is one possible analog of the geodesic flow in the graph. To obtain spectral invariants for this geodesic flow one can consider the transfer operator 
$\mathcal L: \mathbb C^{\mathfrak{P}_+}\to  \mathbb C^{\mathfrak{P}_+}$ defined by 
\[\mathcal L f(\vec e_1,\vec e_2,\ldots):=\sum_{\tau(\vec e_0)=\iota(\vec e_1), \iota(\vec e_0)\not=\tau(\vec e_1)}f(\vec e_0,\vec e_1,\vec e_2,\ldots), 
\] 
where $\mathbb C^{\mathfrak{P}_+}$ is the space of all complex valued functions on $\mf P_+$. The space $\mf P_+$ carries a natural metric topology (see e.g. \cite{BHW23}), from which one can derive a number of natural function spaces. One of them is the space $C^\mathrm{lc}(\mathfrak P_+)$ of locally constant functions, which can be viewed as an analog of the space of smooth functions on a smooth manifold.

The analog of the Laplace-Beltrami operator for the graph $\mf G$ is the Laplacian $\Delta_\mathfrak{X}:\mathbb C^\mathfrak X\to \mathbb C^\mathfrak X$ defined by 
\[\Delta f(x):=\frac{1}{1+q_x}\sum_{d(x,y)=1}f(y),
\]
where $q_x$ for a vertex $x\in \mf X$ is the number of neighboring vertices in the graph and $d(x,y)$ is the minimal numbers of edges needed to connect the vertices $x$ and $y$ by a path. 

In order to describe a graph version of a quantum-classical correspondence, for any spectral parameter $0\not=z\in \mathbb C$  we introduce the \emph{multiplier functions}
$\delta_z:\mathbb C^\mathfrak X\to \mathbb C^\mathfrak X$ given by  
\[\delta_z f(x):=\frac{z+z^{-1}q_x}{1+q_x}f(x).\]
Then we have the following theorem.

\begin{thm}[{{\small\cite[Thm.~11.5, Cor.~9.4]{BHW23}, \cite[Thm.~4.4]{AFH23a}}}] Let $\mf G$ be a finite graph satisfying the hypotheses explained above. Then the restriction $\mc L_\mathrm{lc}$ of $\mc L$ to $C^\mathrm{lc}(\mathfrak P_+)$ satisfies
\[\forall  z\not\in \{0,\pm1\}:\quad \{\mathcal L_\mathrm{lc}=z\}\cong 
\{\Delta_\mathfrak{X} =\delta_z\},
\]
where $\cong$ abbreviates an explicitly known isomorphism of vector spaces. 
\label{thm:graph-qcc_generic}
\end{thm}

Note that while the algebraic equalizer $\{\mathcal L_\mathrm{lc}=z\}$ of $\mathcal L_\mathrm{lc}$ and $z\cdot \mathrm{id}$ on $C^\mathrm{lc}(\mathfrak P_+)$ is simply the eigenspace of $\mathcal L_\mathrm{lc}$ for the eigenvalue $z$, the algebraic equalizer $\{\Delta =\delta_z\}$ of $\Delta$ and the multiplication operator $f\mapsto \delta_z f$ is an eigenspace only if $\delta_z$ is a constant function which is the case if the graph is \emph{homogeneous}, i.e. if the function $\mf X\to\mathbb N,\ x\mapsto q_x$ is constant.

While the isomorphism in Theorem~\ref{thm:graph-qcc_generic} can be obtained directly by combinatorial considerations (see \cite{BHW23}), it can also be obtained in a fashion parallel to the case of compact locally symmetric spaces (see \cite{AFH23a}). The strategy is then to consider the tree 
$\widetilde{\mathfrak G}=(\widetilde{\mathfrak X},\widetilde{\mathfrak E})$
which is obtained as the ``universal covering'' of  $\mathfrak G=(\mathfrak X,\mathfrak E)$ and the (scalar) Poisson transforms defined for that tree.

To make this precise, note first that symmetric graphs of the type we consider here are in one-to-one correspondence with undirected graphs having corresponding properties (connected, no loops, no dead ends). We simply replace pairs of directed edges $(x,y),(y,x)\in \mf X^2$ by the corresponding two-element subset $\{x,y\}\subseteq \mf X$ defining an undirected edge. In this way we can use concepts defined for undirected graphs also in our setting. For example, a connected undirected graph is called a \emph{tree} if it does not contain any cycles. Given an undirected graph we can construct a true universal covering which is a tree. Then we replace the undirected edges by pairs of opposite directed edges und thus arrive at $\widetilde{\mathfrak G}=(\widetilde{\mathfrak X},\widetilde{\mathfrak E})$. 

Recall that one can define the \emph{boundary} of the tree as the set of equivalence classes of geodesic rays, where we call two geodesic rays equivalent if they have infinitely many edges in common. We call the boundary point defined from a geodesic ray its \emph{limit} and denote the boundary of $\widetilde{\mf G}$ by $\widetilde\Omega$. Then there exists a \emph{horocycle bracket} 
\[\langle\cdot ,\cdot \rangle: \widetilde{\mathfrak X}\times \widetilde{\Omega}\to \mathbb Z, \ (x,\omega)\mapsto d(o,y)-d(x,y)\]
where $o,x,y\in \widetilde{\mf X}$ are such that the uniquely determined geodesic rays $[o,\omega[, [x,\omega[, [y,\omega[$ representing $\omega\in\widetilde{\Omega}$ and starting respectively at $o$, $x$, and $y$ satisfy $[o,\omega[\cap [x,\omega[{} = [y,\omega[$.
Then for any $z\in \mathbb C\setminus \{0\}$ we obtain a \emph{Poisson kernel}
\[\tilde{p}_z: \widetilde{\mathfrak X}\times \widetilde{\Omega}\to \mathbb C,\quad (x,\omega)\mapsto z^{\langle x,\omega\rangle}.\]
The boundary $\widetilde \Omega$ inherits a natural topology from the space of geodesic rays in $\widetilde{\mf G}$, so it makes sense to talk about the space $C^\mathrm{lc}(\widetilde \Omega)$ of locally constant functions on $\widetilde\Omega$. Following traditions in harmonic analysis on $p$-adic groups we denote the \emph{algebraic} dual space of $C^\mathrm{lc}(\widetilde \Omega)$ by $\mc D'(\widetilde\Omega)$. By \cite[Prop.~3.9]{BHW22} the space $\mc D'(\widetilde\Omega)$ is naturally isomorphic to the space of finitely additive measures on $\widetilde\Omega$. Locally constant functions can be integrated against such measures, so that we obtain the \emph{scalar Poisson transforms}
\[\widetilde{\mathcal P}_z:\mathcal D'(\widetilde\Omega)\to \mathbb C^{\widetilde{\mathfrak X}},\quad \mu\mapsto \int_{\widetilde\Omega} \tilde p_z(\cdot,\omega)\,\mathrm d\mu(\omega),\]
which take their values in  the respective algebraic equalizer $\{\Delta_{\widetilde{\mf X}} =\tilde\delta_z\}$.
Here $\tilde\delta_z:\widetilde{\mf X}\to\mathbb C$ is the analog of $\delta_z$ for the tree $\widetilde{\mf G}$. In analogy to the symmetric space case the Poisson transforms have natural factorizations of the form $\widetilde{\mathcal P}_z= \tilde\pi_*\circ \tilde p_\mu\circ \widetilde B^*$, where $\tilde\pi_*: \mathcal D'(\widetilde{\mathfrak P}_+)\to \mathrm{Maps}(\widetilde{\mathfrak X})'$ is the pushforward of the canonical projection 
\[\tilde\pi: \widetilde{\mathfrak P}_+\to \widetilde{\mathfrak X},\quad (\vec{e}_j)_{j\in \mathbb N}\mapsto \iota(\vec{e}_1),\] 
$\tilde p_\mu$ is viewed as a multiplier, and $\widetilde B^*$ is the pullback by the \emph{endpoint map} $\widetilde B:\widetilde{\mf P}_+\to \widetilde \Omega$ mapping a geodesic ray to its limit in the boundary. Then the strategy to prove the quantum-classical correspondence is as follows.
\begin{itemize}
\item Work with $\Gamma$-invariant lifts, where $\Gamma$ is the group of deck transformations for the covering $\widetilde{\mf G}\to \mf G$.
\item Apply the factorization for bijective $\mathcal P_z$, i.e. for $z\not\in \{0,\pm 1\}$.
\item  Show that $\tilde \pi_*$, when restricted to $\Gamma$-invariant objects, defines the linear isomorphisms of the quantum-classical correspondence.
\end{itemize}
Apart from being parallel to the case of compact locally symmetric spaces this strategy of proving a quantum-classical correspondence for finite graphs has the advantage of generalizing also to \emph{exceptional spectral parameters}, i.e. to $z=\pm1$, for which the scalar Poisson transform fails to be bijective. In fact, for $\vec e\in \widetilde \Omega$ we consider the set
\[\partial_+\vec e :=\{\omega\in \widetilde \Omega\mid \omega = \widetilde B(\vec e,\vec e_2,\ldots)\}\]
of all boundary points which are limits of geodesic rays starting with $\vec e$. Then for $z\in \mathbb C\setminus \{0\}$ we define the \emph{edge Poisson kernel}
\[\tilde p^\mathrm{e}_z(\vec e,\omega):=z^{\langle \iota(\vec e),\omega\rangle}I_{\partial_+\vec e}(\omega),\]
where $I_{\partial_+\vec e}$ is the indicator function of the set $\partial_+\vec e$. Integration against the edge Poisson kernel gives the \emph{edge Poisson transform}
\[\widetilde{\mathcal P}^\mathrm{e}_z:\mathcal D'(\widetilde\Omega)\to \mathbb C^{\widetilde{\mathfrak E}},\quad \mu\mapsto \int_{\widetilde{\Omega}} \tilde p^\mathrm{e}_z(\cdot,\omega)\,\mathrm d\mu(\omega).\]
The images of the edge Poisson transforms are eigenspaces of the \emph{edge Laplacian}  $\Delta_{\widetilde{\mathfrak E}}:\mathbb C^{\widetilde{\mathfrak E}}\to \mathbb C^{\widetilde{\mathfrak E}}$ defined by 
\[
\Delta_{\widetilde{\mathfrak E}}f(\vec e):=\sum_{\iota(\vec e')=\tau(\vec e), \tau(\vec e')\not=\iota(\vec e)} f(\vec e').
\]

\begin{lem}[{\small\cite[Prop.~2.11\ \& \ 2.18]{AFH23a}}]
For $0\not=z\in \mathbb C$ we have 
\begin{description}
\item[{\rm (i)}] $\widetilde{\mathcal P}^\mathrm{e}_z=\tilde\pi^{\widetilde{\mathfrak E}}_* \circ \tilde p_z\circ \widetilde B^*$ for $\tilde\pi^{\widetilde{\mathfrak E}}: \widetilde{\mathfrak P}_+ \to \widetilde{\mathfrak E},\ (\vec{e}_j)_{j\in \mathbb N}\mapsto \vec{e}_1$.\\
\item[{\rm (ii)}]  $\widetilde{\mathcal P}^\mathrm{e}_z:\mathcal D'(\widetilde\Omega)\to \{\Delta_{\widetilde{\mathfrak E}}=z\}$  is a linear isomorphism. 
\end{description} 
\end{lem}

Replacing $\mathcal P_z$ for $z\not\in \{0,\pm 1\}$ in the above strategy by $\widetilde{\mathcal P}^\mathrm{e}_z$ for $z\not= 0$ then yields a quantum-classical correspondence for finite graphs which is valid also for exceptional parameters.

\begin{thm}[{\small\cite[Thm.~4.5]{AFH23a}}]~Let $\mf G$ be a finite graph satisfying the hypotheses explained above. Then 
\[\forall  z\not=0:\quad \{\mathcal L_\mathrm{lc}=z\}\cong 
\{\Delta_{\mathfrak{E}} =z\},
\]
where $\Delta_\mathfrak{E}$ is the edge Laplacian for $\mf G$ and $\cong$ abbreviates an explicitly known isomorphism of vector spaces. 
\label{thm:graph-qcc_exceptional} 
\end{thm}
 
In the case of finite graphs we can substantiate the expectation that eigenspaces for exceptional spectral values carry topological information as suggested by examples in the case of locally symmetric spaces (see Section~\ref{subsubsec:exceptional spectral parameters}).  

Let $c(\mathfrak G)$ the \emph{cyclomatic number} of the finite graph $\mf G$, i.e. the minimal number of edges which need to be removed from the graph in order to break all cycles.

\begin{thm}[{\small\cite[Theorems~C\ \&\ D]{AFH23a}}]~Let $\mf G$ be a finite graph satisfying the hypotheses explained above.
\begin{description}
\item[{\rm (i)}]  $\dim\{\Delta_{{\mathfrak E}}=1\}=\begin{cases}
c(\mathfrak G)&\text{if } c(\mathfrak G)\not=1, \\
c(\mathfrak G)+1&\text{if } c(\mathfrak G)=1. 
\end{cases}$
\item[{\rm (ii)}] If $c(\mathfrak G)\not=1$ or $\mathfrak G$ is not bipartite, then
\[\dim\{\Delta_{{\mathfrak E}}=-1\}=\begin{cases}
c(\mathfrak G)&\text{if } \mathfrak G \text{ is bipartite,}\\
c(\mathfrak G)-1&\text{if } \mathfrak G \text{ is not bipartite.}
\end{cases}\]
\item[{\rm(iii)}] If $c(\mathfrak G)=1$ and $\mathfrak G$ is bipartite, then
$\dim\{\Delta_{{\mathfrak E}}=-1\}=2$.

\end{description}
\end{thm}

\subsection{Affine buildings}\label{subsec: buildings}

Affine buildings of higher rank occur for instance as Bruhat-Tits buildings of $p$-adic algebraic groups beyond $\mathrm{SL}_2(\mathbb K)$. In contrast to the case of trees there are few affine buildings of higher rank which are \emph{not} Bruhat-Tits buildings. Thus if makes sense to make use of group theory and harmonic analysis to study quantum and classical analogs of geodesic and Weyl chamber flows for higher rank affine buildings. 

The quantum side is well established with Hecke algebras taking the place of invariant differential operators, see e.g. \cite{Ma71}. Geometric versions trying to avoid group theory appeared later, see e.g. \cite{Pa06a}. They add new points of view but do not extend the scope of the theory a lot.   

While there is some literature on higher rank Cartan actions (see e.g. \cite{Mo95}), analogs of Weyl chamber flows suitable to build quantum-classical correspondences in the spirit of \cite{HWW21} do not seem to be available as of now. What is avaivalable are analogs of the representation theoretic tools used in the context of locally symmetric spaces. For instance S. Kato's papers \cite{Ka81,Ka82} contain scalar Poisson transforms from unramified principal series representations to Hecke eigenspaces and a characterization of the spectral parameters for which these Poisson transforms are bijective. Geometric versions of Poisson transforms were studied by Mantero and Zappa (see \cite{MZ17} and the references given there), mostly for rank two buildings.

Thus, in order to establish quantum-classical correspondences in the context of higher rank affine buildings, one first has to do some more foundational work on the underlying dynamical systems. 

It should be noted, however, that there is related material that may become relevant in our context at some point, such as for instance the multivariate geometric zeta functions for higher rank $p$-adic groups Deitmar and Kang introduced in \cite{DK14}.

\section{Open problems}

The existing results on quantum-classical correspondences for locally symmetric spaces were obtained with a variety of different methods and have very different degrees of generality. The most complete set of results is available for hyperbolic surfaces. They range from trivial and thin to cocompact fundamental groups. The methods employed range from symbolic dynamics and classical Fourier analysis to noncommutative harmonic analysis and microlocal analysis. The known extensions to higher dimension and higher rank work for different settings which depend on the methods employed. Trying to further extend the scope one meets a number of difficulties. Some are more of a technical nature, such as explicitly determining composition series of principal series representations, some are truly conceptual such as finding a suitable replacement for resolvents when going from a single operator to a set of commuting operators. 

In this section we formulate a number of concrete challenges that could be seen as possible next steps in realizing the program outlined in Section~\ref{subsubsec: program}. The focus will be on problems in establishing quantum-classical correspondences via microlocal methods. Given the results described in the previous sections there will be a number of problems asking for various types of generalizations such as  

\begin{itemize}
\item globally  vs. locally symmetric spaces
\item higher rank vs. rank one
%\item general parabolic vs. minimal parabolic
\item higher dimension vs. surfaces
\item $\Gamma$ general vs. $\Gamma$ thin or $\Gamma$ cocompact
\item Orbits or parameters general vs. orbits or parameters generic (regular) 
\end{itemize}

There are also some problems asking for clarification of the relation between the methods and results laid out here with the ones contained in previous work.

\subsection{Quantum resonances in higher rank}

The Laplacian of a noncompact Riemannian symmetric space $G/K$ can be diagonalized by the non-Euclidean Fourier transform. Thus its resolvent can be written as an explicit integral operator with the integrand depending on a complex parameter. For spaces of low rank (actually $1$ and $2$) this representation can be used to construct meromophic continuations via suitable deformations of contours, \cite{HP09,HPP16,HPP17a,HPP17b}. Already for rank $2$ spaces the construction of suitably deformed contours turns out to be quite tricky with no obvious pattern of generalization. 

\begin{OP}[Resolvent poles for Laplacians]~ 
Extend the method used in \cite{HP09,HPP16,HPP17a,HPP17b} to determine quantum resonances for the Laplacian of symmetric spaces of arbitrary rank (cf. Subsection~\ref{subsubsec:Resolvent poles}).
\begin{enumerate}
\item[(i)]~Find a systematic way to do the contour deformations working in arbitrary rank\footnote{Note that Mazzeo and Vasy in \cite{MV05} studied analytic continuations of the resolvent of the Laplacian in arbitrary rank. They do, however, not solve the problem formulated here, since they do not go around possible singularities but introduce cuts and branched coverings wherever a pole might occur.}.
\item[(ii)]~Describe the poles of the resulting meromorphic continuation, i.e. the quantum resonances of the Laplacian.
\item[(iii)]~Determine the $G$-representations obtained from the residue operators at the quantum resonances.
\end{enumerate}
\end{OP}

In view of the Plancherel decomposition of the regular representaion of $G$ on $L^2(G/K)$ it seems plausible that one should aim for quantum resonances which are spectral invariants not just for the Laplacian, but rather for the entire commutative algebra $\mathbb D(G/K)$ of invariant differential operators. In rank one this requires no essential change since then $\mathbb D(G/K)$ is generated by the Laplacian. If the rank is $\ell$, the algebra $\mathbb D(G/K)$ is isomorphic a polynomial algebra in $\ell$ variables due to results of Harish-Chandra and Chevalley, see Remark~\ref{Chevalley}. In fact, one can still use the non-Euclidean Fourier transform to diagonalize $\mathbb D(G/K)$ simultaneously. The resulting eigenvalues are provided by the Harish-Chandra isomorphism, i.e. the spectral parameters are elements of $\mf a^*_\C$. 

\begin{OP}[Resolvent cohomology classes for invariant differential operators]~ 
Modify the definition of the Taylor spectrum so that it applies to the commuting family $\mathbb D(G/K)$ of differential operators. A good starting point might be $G=\mathrm{SL}_3(\R)$, since in that case $\mathbb D(G/K)$  can be described very explicitly, see \cite{BCH21}.
\begin{enumerate}
\item[(i)] Use the Harish-Chandra isomorphism to express the Taylor spectrum of $\mathbb D(G/K)$ in a way independent of a choice of generators.
\item[(ii)] Follow \cite[\S\,2]{Sa03} to construct a resolvent cohomology class for $\mathbb D(G/K)$ depending on a parameter in $\mf a^*_\C$.
\item[(iii)] Work out the $G$-invariance and equivariance properties of the spaces involved in the construction of the resolvent class. 
\item[(iv)] Build a framework in which one can talk about meromorphic continuation of the resolvent class and show that the resolvent class has meromorphic continuation to all of $\mf a^*_\C$.
\item[(v)] Determine the divisor of the meromorphic continuation of the resolvent class.
\item[(vi)] Follow the methods of \cite[Chap.~5]{GH78} to calculate suitable residues of the meromorphically continued resolvent classes.
\item[(vii)] Give the residues from (vi) a representation theoretic interpretation.
\end{enumerate}
\label{OP:higher rank resolvents1}
\end{OP}

There are other conceivable ways to describe resolvents of (commuting) families of operators and their meromorphic continuations.

\begin{OP}[Multivariate resolvents via Laplace transform]~
Take the description of the resolvent of a single operator $A$ via the Laplace transform of the one-parameter semigroup $S(t)$ generated by $A$ (cf. Example~\ref{ex:DFG} and \cite[\S~34.1]{La02}) as a model and construct higher rank resolvents as multivariate Laplace transforms. Then study meromorphic continuations and residues. 
\label{OP:higher rank resolvents2}
\end{OP}

The view that the ``correct'' approach to describe quantum resonances for symmetric spaces is to  study the spectral theory not only of the Laplacian but the simultaneous spectral theory of $\mathbb D(G/K)$ is supported by the relation between resonances and scattering poles predicted by physics. In fact, following \cite{STS76} one may view the \emph{standard intertwining operators} for spherical principal series representations as scattering matrices for the multitemporal wave equation 
\[\forall D\in \mathbb D(G/K):\quad D_xu(tx) = \Gamma(D)_{\partial t} u(t,x)\]
on the symmetric space coming from the Harish-Chandra isomorphism $HC$ (see Section~\ref{sec:harishchandra}) by letting $\Gamma(D)_{\partial t}$ be the constant coefficient differential operator with symbol $HC(D)$. The standard  intertwining operators depend on parameters in $\mf a^*_\C$ and it is well-known (see e.g. \cite[Chap.~10]{Wa92}) that they have meromorphic continuation to all of $\mf a^*_\C$.

\begin{OP}[Resolvent poles vs. scattering poles for higher rank symmetric spaces]~
 Extend the results of \cite{HHP19} (see Theorem~\ref{thm-scatt-poles-are-resonances}) to higher rank using standard intertwining operators as scattering matrices and the resolvent poles of Problems~\ref{OP:higher rank resolvents1} or \ref{OP:higher rank resolvents2} as quantum resonances.
\end{OP}

In \cite[\S\,6]{BO12} Bunke and Olbrich study meromorphic continuations of resolvent kernels for the Laplacian on convex cocompact locally symmetric spaces. In the process they compare the resolvent kernel of the locally symmetric space with the resolvent kernel of the corresponding symmetric space. Moreover, they make use of the scattering matrix for convex cocompact locally symmetric spaces they had defined in \cite{BO99,BO00}.

\begin{OP}[Resolvent poles vs. scattering poles for convex cocompact locally symmetric spaces] Use the techniques of \cite{BO99,BO00,BO12} to extend the results of \cite{HHP19} to
convex cocompact locally symmetric spaces of rank $1$.
\end{OP}

Convex cocompact locally symmetric spaces are rare in higher rank (basically there are only products of rank $1$ spaces). As a replacement various authors (see e.g. \cite{GW12,KLP18}) have studied \emph{Anosov representations} of discrete groups.

\begin{OP}[Bunke-Olbrich for higher rank]~ Generalize the techniques of \cite{BO99,BO00,BO12} to locally symmetric spaces whose fundamental group is the image of an Anosov representation.
\end{OP}

A detailed scattering theory for locally symmetric spaces generalizing \cite{STS76} would be of interest in the study of automorphic forms, see e.g. \cite{LP76} for the case of hyperbolic surfaces. In \cite{JZ01} Ji and Zworski take up the subject for locally symmetric spaces of $\mathbb Q$-rank $1$.

\begin{OP}[Automorphic scattering in higher rank]~
Extend the scattering theory of Semenov-Tian-Shansky laid out in \cite{STS76} for symmetric spaces to locally symmetric spaces (see \cite{LP76} for the case of hyperbolic surfaces).
\end{OP}

\subsection{Ruelle-Taylor resonances}

The obvious problems in this context are
\begin{itemize}
\item Construct Ruelle-Taylor resonances for the $A$-action from Remark~\ref{rem:A-action on G/M} for any locally symmetric space $\mc X=\Gamma\bs G/K$.
\item Show that these Ruelle-Taylor resonances have a natural band structure.
\item Establish a quantum-classical correspondence for the corresponding resonant states.
\end{itemize}

Attacking these problems head-on is way too ambitious for the moment. We formulate some more modest problems.

\begin{OP}[Band structure of Ruelle-Taylor resonances]~ 
Give a precise description of the band structure of Ruelle-Taylor resonances in the spirit of \cite{DFG15} (cf. the comments in Remark~\ref{rem:vector bundles}).  
\end{OP}

\begin{OP}[Quantum-classical correspondences for real hyperbolic manifolds]~
As mentioned in Remarks~\ref{rem:Ruelle resonance - Hadfield} and \ref{rem:Ruelle resonance - Bonthonneau-Weich} Bonthonneau and Weich constructed Ruelle resonances for Riemannian manifolds with finite volume, negative curvature and the unbounded part consisting of finitely many real hyperbolic cusps, whereas Hadfield constructed Ruelle resonances for convex cocompact  real hyperbolic manifolds, see \cite{BW22,Ha17}. Moreover, in \cite{Ha20} Hadfield established quantum-classical correspondences for the latter manifolds.
\begin{enumerate}
\item[(i)] Construct Ruelle resonances for general real hyperbolic manifolds.
\item[(ii)] Extend Hadfield's quantum-classical correspondences to general real hyperbolic manifolds.
\end{enumerate}
\end{OP}

\begin{OP}[Ruelle resonances for convex cocompact rank one locally symmetric spaces]~ Construct Ruelle resonances for general convex cocompact rank one locally symmetric spaces.
\end{OP}

\begin{OP}[Ruelle-Taylor resonances for thin $\Gamma$ in higher rank]~ Construct Ruelle-Taylor resonances for  locally symmetric spaces with fundamental group in the image of an Anosov representation.
\end{OP}

\subsection{Dynamical zeta functions, $\Gamma$-cohomology and transfer operators}
\label{subsubsec:Juhl-Bunke-Olbrich}

Following pioneering work of Patterson \cite{Pa89} and Juhl \cite{Ju93,Ju01}, in the late 1990s Bunke and Olbrich as well as Patterson and Perry established links between the dynamical zeta functions and $\Gamma$-cohomology with coefficients in $G$-representations (see Subsections~\ref{subsubsec: Gamma cohomology} \& \ref{subsec:zeta} and \cite{BO99,PP99,BO00,BO12}).

%\paragraph{Precise relation with the work of Juhl and Bunke-Olbrich}\label{subsubsec:Juhl-Bunke-Olbrich}

The focus of  Juhl's books \cite{Ju93,Ju01} is on dynamical and Selberg zeta functions for the geodesic flow on compact hyperbolic manifolds. In \cite[Chap.~3]{Ju01} Juhl establishes a Lefschetz-type formula, where the geometric side is expressed in terms of the contractive part of the Poincar\'e map of closed geodesics and the spectral side is expressed in terms of the $\mf n$-cohomologies of the irreducible components of $L^2(\Gamma\backslash G)$. As an application he gets formulas for the divisors of (twisted) Selberg zeta functions in terms of $\mf n$-cohomologies, which he then rephrases in terms of $\Gamma$-cohomologies with coefficients in principal series representations. To that end he uses de Rham type complexes to compute such $\Gamma$-cohomologies. The resulting characterization of the zeta divisors in terms of the spectral decomposition of $L^2(\Gamma\backslash G)$ obscures the dynamical meaning, so Juhl sets out in  \cite[Chap.~4]{Ju01} to develop a kind of Hodge theory for the sphere bundle $S(\Gamma\backslash G)=\Gamma\backslash G/M=\M$ using the Anosov property of the geodesic flow. It is in the detailed description of the complexes and operators on $\M$ given in \cite[Chap.~5]{Ju01} that Poisson transforms make their appearance.

The published work of Bunke and Olbrich exclusively deals with locally symmetric spaces of rank one.  It in particular describes zeros of dynamical zeta functions in terms of $\Gamma$-cohomology with coefficients in spaces of distributions on the boundary (see e.g. \cite[Thm.~1.3]{BO99}), which again might be viewed as a result purely on the classical side. But even more than Juhl, Bunke and Olbrich  make extensive use not only of the Poisson transform, but also of the scattering matrix in establishing such results. 

So far the wealth of information Juhl, Bunke and Olbrich have generated in this context has not found its way into the efforts of the more PDE and dynamical systems oriented communities in spectral geometry\footnote{Neither has it been taken up a lot in the literature on noncommutative harmonic analysis and representation theory.}. One reason for this may be that they use heavy group and representation theoretic machinery as well as results originating in the algebraic analysis in the  spirit of \cite{SKK73} which to this day remains to be presented in a form accessible to a wider mathematical audience.

\begin{OP}[Relation with the works of Juhl and Bunke-Olbrich]
~Work out the precise relation between the cohomological methods of Juhl and Bunke-Olbrich on the one side and the dynamical respectively microlocal methods leading to resonances on the other side. A starting point could be to compare the descriptions of the zeta divisors given in \cite[\S~5]{BO99} and \cite[\S~5]{GHW18} (the latter is only for surfaces).
\end{OP}

%\paragraph{Multivariate zeta functions in higher rank}

Considering the Weyl chamber action as the suitable analog of the geodesic flow, the dynamical zeta functions should count compact $A$-orbits (see \cite{Mo1,Mo2} for results on the frequency of such orbits). One may follow various approaches:

\begin{itemize}
\item[1.] Count according to volume. 

This leads to one-variable zeta functions  (see \cite{De,DKV79} for background).

\item[2.] Count according to the geometry. 

A compact $A$-orbit $c=\xi\cdot A$ is diffeomorphic to a torus and the stabilizer of a point is a lattice in $A$.
Let $\ell(c)=\big(\ell_1(c),\ldots, \ell_r(c)\big)$ be successive length minima of the lattice (use the Killing form to fix a metric on $A$). Then one might define
$$Z_{\mathrm R}(\lambda)= \prod_c \big(1-e^{-\lambda\ell(c)}\big)^{\pm 1}$$
for $\lambda\in \mathbb C^r$. In order to generalize the other zeta functions from Subsection~\ref{subsec:zeta} one needs to have a reasonable Poincar\'e map. Given a lattice $L\subseteq A$ one might try to take the following steps:
\begin{itemize}
\item[(a)]~ Lift the action to $T(\Gamma\backslash G/M)$.
\item[(b)]~ Define a Poincar\'e map for each element of $L$.
\item[(c)]~ For $a\in L$ let $\xi\cdot a=\xi\in\Gamma gM$. Then there exist $\gamma\in  \Gamma$ and $m\in M$ such that $gam=\gamma g$ and right multiplication $\rho_{am}$ by $am$ yields
the identification of $\rho_{am}$ with $\Ad(am)$. 
%$$\xymatrix{ T_g(G)\ar[rr]^{\rho_{am}}\ar[d]_{\cong}&& T_{\gamma g}(G)\ar[d]^{\cong}\\ \mathfrak g\ar[rr]_{\mathrm{Ad}(am)}&&\mathfrak g }$$
Thus the monodromy on $T(\Gamma\backslash G/M)= \Gamma\backslash G\times_M (\mathfrak n+\mathfrak a+\overline{\mathfrak n})$ is given by $\mathrm{Ad}(am)$. Then follow  \cite{BO95,De}.
\end{itemize}

\item[3.] Write the zeta function as a Fredholm determinant of a suitable (transfer) operator. 

This works fine for subshift dynamics and requires symbolic dynamics in the case of geodesic flows (known for specific hyperbolic surfaces, see \cite{Po1,Po2}). 

\end{itemize}

Note that the possibility to write zeta functions as (Fredholm) determinants is the reason why symbolic dynamics is used in the effort to construct meromorphic extensions of dynamical zeta functions.

The following three problems take up aspects of the approach to dynamical zeta functions for Weyl chamber flows just outlined.

\begin{OP}[Geometric invariants for $A$-orbits]
~Find other geometric invariants of compact $A$-orbits according to which one can count. Is it possible to use the restriction of the Sasaki metric on $\Gamma\backslash G/M\subseteq T(\Gamma\backslash G/K)$ to the orbit to construct suitable invariants?  
\end{OP}

\begin{OP}[Multivariate dynamical zeta functions]
~Work out a basic theory of multivariate dynamical zeta functions
\begin{itemize}
\item[(a)]~ for flows and torus orbits
\item[(b)]~ for higher rank subshifts of finite type
\end{itemize}
\end{OP}

\begin{OP}[Symbolic dynamics for Weyl chamber flows]
~Construct symbolic dynamics for Weyl chamber flows in rank $2$ examples. The product of Schottky surfaces is dealt with in \cite{Po20}. One could start with $\mathrm{SL}_3(\mathbb Z)\backslash\mathrm{SL}_3(\R)/\mathrm{SO}_3$.
\end{OP}

The treatment of the transfer operator for the modular surface by Chang and Mayer in \cite{CM99} shows a clear connection with the hyperbolic nature of the geodesic flow. This is a lot less transparent for the transfer operators derived from the symbolic dynamics other discrete subgroups of $\mathrm{PSL}_2(\R)$.

\begin{OP}[Hyperbolic nature of Pohl's symbolic dynamics] 
~Make the connection between Pohl's reduced symbolic dynamics for the geodesic flow and the hyperbolic nature of the geodesic flow explicit (cf. the discussion in Subsection~\ref{subsuubsec: transfer operators}).
\end{OP}

The next problem is suggested by the similarities between Maass and modular forms.

\begin{OP}[Quantum interpretation of modular cusp form]
~Give a reasonable quantum system having modular cusp forms in the sense of Example~\ref{ex:period-polynomials} as states and such that it quantizes a classical system to which one can naturally associate period polynomials. 
\end{OP}

The last problem in this section describes a strategy to extend the Lewis-Zagier correspondence to higher rank with Remark~\ref{rem:cusp-cohom} as point of departure.

\begin{OP}[Cuspidal cohomology]
\begin{enumerate}
\item[(i)] Use the multiplicity space from Remark~\ref{rem:cusp-cohom} and $K$-invariant vectors of spherical principal series representations to relate the cuspidal $\Gamma$-cohomology spaces for spherical principal series representations to Maass cusp forms. 
\item[(ii)] Follow the isomorphisms proving Theorem~\ref{thm:DH05-0.2} through to construct a natural isomorphism of vector spaces between Maass cusp forms and cuspidal $\Gamma$-cohomology spaces.
\item[(iii)] Clarify the relation between the various parabolic $\Gamma$-cohomologies from \cite{BLZ15} and cuspidal $\Gamma$-cohomology in the case of $G=\mathrm{PSL}_2(\R)$.
\item[(iv)] Give conditions under which cuspidal $\Gamma$-cohomology  and parabolic $\Gamma$-cohomology agree.
\end{enumerate}
\end{OP}

\subsection{Horocycle flow}

Let $\X=\Gamma\backslash G/K$ be a locally symmetric space of finite volume. The \emph{horocycle flow} of $N$ on $\Gamma\backslash G$ is given by right multiplication. In contrast to the Weyl chamber flow it does not descend to $\Gamma\backslash G/M$. Following \cite{FF03} (see Example~\ref{ex:Flaminio-Forni})  one can define the space $\mc J(\Gamma\backslash G)$ of compactly supported distributions on $\Gamma\backslash G$ which are invariant under the horocycle flow. In addition we consider the subspace $\mc J(\Gamma\backslash G)^M$ of  $M$-invariant distributions in $\mc J(\Gamma\backslash G)$. This amounts to taking $MN$-invariant distributions on $\Gamma\backslash G$.

\begin{OP} Check the following claims:
\begin{enumerate}
\item[(i)] $\mc J(\Gamma\backslash G)$ is invariant under the natural $MA$-action on the distributions on $\Gamma\backslash G$ coming from the right $MA$-action on $\Gamma\backslash G$.
\item[(ii)]  $\mc J(\Gamma\backslash G)$ and $\mc J(\Gamma\backslash G)^M$ have compatible direct integral decompositions coming from the right regular representation of $G$ on $C_c^\infty(\Gamma\backslash G)$.
\item[(iii)] The summands in the decompositions from (ii) are invariant under the $MA$-action.
\end{enumerate}
\label{OP:Flaminio-Forni1}
\end{OP}

\begin{OP} 
\begin{enumerate}
\item[(i)] Determine the decompositions from Problem~\ref{OP:Flaminio-Forni1}(ii) in terms of conical distributions, \cite[\S~II.4]{He94}, and describe the $MA$-actions on summands.
\item[(ii)] Derive a quantum-classical correspondence for $\X$ from (i) and compare the result in the cocompact case with \cite[Thm.~5.1]{HWW21}.
\end{enumerate}
\label{OP:Flaminio-Forni2}
\end{OP}

It may very well turn out that instead of just considering conical distributions in Problem~\ref{OP:Flaminio-Forni2}(ii) one will have to consider $M$-finite vectors in the space of $N$-invariant distributions as in \cite{Em14}. Also one should be looking for a connection with the $\mf n$-cohomology techniques applied by Juhl and Bunke-Olbrich (cf. the comments in Subsection~\ref{subsubsec:Juhl-Bunke-Olbrich}).

\subsection{Patterson-Sullivan distributions}
\label{sec:PS-distributions}

In \cite[Thm.~5.2]{GHW21} it is shown that the Patterson-Sullivan distributions on compact rank  one locally symmetric spaces can be easily obtained from the intertwiners between Laplace eigenfunctions and the generalized Ruelle resonant states from the first band. Of course one would like to have a higher rank generalization.

\begin{OP}
Use the first band Ruelle resonant states in compact locally symmetric spaces to generalize \cite[Thm.~5.2]{GHW21} to higher rank.
\end{OP}

To generalize the realizations of Patterson-Sullivan distributions as residues of their zeta functions from hyperbolic surfaces to rank one locally symmetric spaces turned out to be a nontrivial task as  neither the symbolic dynamics nor the weighted trace formula used in \cite{AZ} have immediate generalizations in higher dimension. For compact hyperbolic manifolds steps in this direction have been taken in \cite{Em14}. As was mentioned in Section~\ref{subsubsec:qcc-applications} the approach taken in \cite{SWB23} using \cite{GHW21} and \cite{DyZw} was actually more successful.

\begin{OP} Extend the results of \cite{SWB23} to general Riemannian locally symmetric spaces.
\end{OP}

Even though the microlocal approach turned out to be more efficient in the problem at hand, there remain open questions about the harmonic analysis approach taken in   \cite{AZ} and \cite{Em14}.

\cite{AZ} makes use of a  pseudo-differential calculus due to Zelditch (see \cite{Ze84}) on the basis of the non-Euclidean Fourier-Helgason transform.  Although lifted quantum limits do not depend on the specific \psdiff\ calculus chosen for their definition, it is useful, for establishing invariance properties, to have an equivariant calculus. For hyperbolic surfaces, based on the non-Euclidean Fourier analysis  and closely following the Euclidean model, such a calculus was provided by  Zelditch \cite{Ze84}. In \cite{Sch10} this calculus was extended to rank one symmetric spaces. Using this calculus the construction of the Patterson--Sullivan distributions and the proof of the asymptotic equivalence from \cite{AZ} can be generalized. However, due to singularities arising from Weyl group invariance, it is difficult to construct an equivariant non-Euclidean \psdiff\ calculus in higher rank. Silberman and Venkatesh \cite{SV07}, generalizing work of Zelditch and Wolpert for surfaces to compact locally symmetric spaces, introduced a representation theoretic lift as a replacement for a microlocal lift. They sketch, in \cite[Remark 1.7(4) and \S 5.4]{SV07}, a proof that the representation theoretic lift asymptotically gives the same result as a microlocal lift using \psdiff\ operators.

If one had a formula intertwining Patterson-Sullivan distributions $PS^\Gamma_h$
into lifted Wigner distributions $W_h$, one might be able to deduce \eqref{W-PS-diagonal} as a corollary. Presumably, an intertwining formula holds only for special \psdiff\ calculi.

\begin{OP}
Establish a non-Euclidean pseudo-differential calculus for general Riemannian symmetric spaces of noncompact type. First steps were done in \cite{Sch10}, but the problems encountered there lead to the use of a different geometric calculus in \cite{HHS}.
\end{OP}

\begin{OP}
Silberman and Venkatesh introduced a \emph{representation theoretic lift} as an alternative to the microlocal lift in proving invariance properties of limit measures, see \cite{SV07,Si15} and \cite{BO06}. As pointed out in \cite[\S~5]{Si15} the representation theoretic lift works for finite volume and does not need compactness of the locally symmetric space.
\begin{enumerate}
\item[(i)] Work out the argument for the equivalence of microlocal and representation theoretic lifts sketched in \cite[\S~5.2]{Si15} for compact locally symmetric spaces.
\item[(ii)] Extend the theory of Patterson-Sullivan distributions to locally symmetric spaces of finite volume.
\end{enumerate}
\end{OP}

 \subsection{Singular situations}

We have discussed several singular situations in this article. The first was singular $G$-orbits $G\cdot[\mathbf 1,\xi]\cong G/K_\xi$ in $T^*(G/K)$ which results in the right $A_\xi$-action on $G/K_\xi$, see Remark~\ref{rem:singular flows} and \cite{Si15}. This action descends to $\M_\xi:=\Gamma\backslash G/K_\xi$. Theorem~\ref{thm:Bialo8.2} provides this action with a symplectic, i.e. classical, interpretation.

\begin{OP}[Singular generalized energy shells]~ Extend the various results presented in this article for the regular generalized energy shells $\M$ to singular generalized energy shells $\M_\xi$.
\label{OP:singular flow}
\end{OP}

The second kind of singular situation we have encountered are exceptional spectral parameters, for which Poisson transforms fail to be bijective. These exceptional parameters required separate treatment in the search for quantum-classical correspondences, see Section~\ref{subsubsec:exceptional spectral parameters}. 
Some examples considered there suggest that the multiplicities of exceptional Ruelle resonances have topological interpretations. 

\begin{OP}[Topological interpretation of exceptional multiplicities]~~
Give a topological interpretation of  $\dim({}^\Gamma(\on{soc}(\pst{\mu}))^{-\infty})$ for the exceptional spectral parameters $\mu$ (see Theorem~\ref{thm:spectral_corres_SOn} and Remark~\ref{rem:QCC rank 1}).

In view of the cohomological constructions done in \cite{TW89}, the observation from Theorem~\ref{thm:discrete series} that the representations associated with the exceptional spectral parameters are discrete series representation of associated non-Riemannian symmetric spaces, may be helpful in attacking this problem.
\end{OP}

All the results described in Section~\ref{subsubsec:exceptional spectral parameters} are restricted to rank one spaces. One reason for this was that in \cite{AH23} we needed explicit information on the spherical principal series representations involved. So the first problem to solve when trying to extend the results from \cite{AH23} to higher rank are purely representation theoretic.

\begin{OP}[Spherical principal series for exceptional parameters]~
Let $\pi_\lambda$ be a spherical principal series representation of $G$ with exceptional spectral parameter $\lambda$.
\begin{enumerate}
\item[(i)] It is known that $\pi_\lambda$ is not irreducible. Determine the composition series and the socle of $\pi_\lambda$.
\item[(ii)] Determine the kernel of the scalar Poisson transform $\mc P_\lambda$.
\item[(iii)] Determine the minimal $K$-types of $\pi_\lambda$.
\item[(iv)] Show that the vector valued Poisson transform $\mc P^\text{v}_\lambda$ associated with $\lambda$ and the sum of all minimal $K$-types is injective on the kernel of $\mc P_\lambda$. 
\item[(v)] Show that $\mathcal{P}^\text{v}_\lambda$  is injective on the socle of $\pi_\lambda$.
\item[(vi)] Determine the image of the socle of $\pi_\lambda$ under $\mc P^\text{v}_\lambda$.
\end{enumerate}
\label{OP:exceptional1}
\end{OP}

In Problem~\ref{OP:exceptional1} the discrete subgroup $\Gamma$ of $G$ does not appear. In the next problem it does and we do not specify which further assumptions one should make. Probably the simplest case is to choose $\Gamma$ cocompact.

\begin{OP}[$\Gamma$-invariant distribution vectors in spherical principal series representations]~ Let $\pi_\lambda$ be a spherical principal series representation of $G$ with exceptional spectral parameter $\lambda$.
\begin{enumerate}
\item[(i)] Show that all $\Gamma$-invariant distribution vectors for $\pi_\lambda$ are in the socle of $\pi_\lambda$.
\item[(ii)]  Show that the vector valued Poisson transform $\mc P^\text{v}_\lambda$ associated with $\lambda$ and the sum of all minimal $K$-types is injective on the socle of $\pi_\lambda$. 
\item[(iii)] Determine the image of the space of $\Gamma$-invariant distribution vectors for $\pi_\lambda$ under $\mc P^\text{v}_\lambda$.
\end{enumerate}
\label{OP:exceptional2}
\end{OP}

When we discussed Patterson-Sullivan distributions in Subsection~\ref{Lifts and PS-distributions} we saw that the asymptotic equivalence of lifted quantum limits and Patterson-Sullivan distributions depended on the lifted quantum limit being regular.

\begin{OP}[Singular  lifted quantum limits]
Extend \cite[Thm.~7.4]{HHS} to singular lifted quantum limits. 
\end{OP}

The final problem we want to state here showed up implicitly in many places in this article. It is related to the boundary values of functions on a locally symmetric space which occur for instance in the inverse of regular Poisson transforms. The most general construction of these boundary value maps has been given in \cite{KKMOOT}, which depends massively on \cite{SKK73}. While the main result of \cite{KKMOOT}, the solution of the so-called \emph{Helgason conjecture} characterizing the bijectivity of the Poisson transform $\mc P_\lambda$ in terms of $\lambda$, has been reproved several times in ways independent of \cite{SKK73}, except for the rank one case (done already by Helgason, see also \cite{HHP19}) there is still no elementary construction of the boundary value map (see \cite{HHP18} for a discussion and a PDE-construction of boundary value maps for generic parameters).

\begin{OP}[Boundary value maps] Give a PDE-construction of the boundary value maps from \cite{KKMOOT}. 
\end{OP}

\subsection{Graphs and Buildings}

So far the quantum-classical correspondences established for graphs (see \cite{BHW23,AFH23a}) deal only with finite graphs. One reason is that resonances for the graph Laplacians have not been studied so far.

\begin{OP}[Quantum resonances for infinite graphs] Establish meromorphic continuation of the resolvent of graph Laplacians for infinite graphs and determine its poles (resonances) and residue operators.
\end{OP} 

The geodesic flow acting on geodesic rays is only one possible analog of the geodesic flow on Riemannian manifolds. Another possibility would be to consider paths $(\vec{e}_j)_{j\in \mathbb Z}$ as \emph{geodesic lines}. In that case the geodesic flow is generated by the twosided shift. The analogy to the classical side of the quantum-classical correpondence for locally smmetric spaces suggests that one has to find a way to formulate hyperbolic properties of the shift dynamics on the space of geodesic lines. Note here that the use of geodesic rays and the transfer operator rather mimick the procedure given in \cite{CM99} for the modular surface, which throws out the stable part of the dynamics altogether and only keeps the unstable part. 

\begin{OP}[Hyperbolicity of geodesic flows on graphs] Formulate hyperbolicity properties of the twosided shift on the space of geodesic lines in a graph. 

A possible approach to this question is to build stable and unstable foliations using horocycles of the covering tree. Using these one might try to construct anisotropic Sobolev spaces analogous to the ones used for locally symmetric spaces (see Subsection~\ref{subsubsec: RT-resonances}).

Further goals could be to introduce resolvents of the shift operator and prove meromorphic continuation. The residues could then be used to introduce graph analogs of \emph{invariant Ruelle densities}, which in the case of compact rank one locally symmetric spaces occur as flat traces of the residue operators of Ruelle resonances (cf. \cite[\S~2]{GHW21}). 
\label{problem:hyperbolicity of geoflow - graph}
\end{OP} 

The space of geodesic lines just considered came up already in \cite{Ah88}, where one also finds a tree version of the Radon transform. Recall from  Section~\ref{sec:PS-distributions} that for compact rank one locally symmetric spaces Ruelle resonant states could be used to describe Patterson-Sullivan distributions. In that case they appear as products of distributions on the space of geodesic lines. All this suggests that it is possible to construct Patterson-Sullivan distributions for finite graphs.

\begin{OP}[Patterson-Sullivan distributions for finite graphs]~Build a theory of  Patterson-Sullivan distributions for finite graphs.

In particular, it should be clarified how the Patterson-Sullivan distributions are related to eigenfunctions of the transfer operator and invariant Ruelle densities (cf. Problem~\ref{problem:hyperbolicity of geoflow - graph}). Moreover, one would want an analog of the residue formula \cite[Thm.~4.1]{SWB23} connecting Patterson-Sullivan distributions with dynamical zeta functions. Note here that a graph analog of the pairing formula  \cite[Thm.~6.1]{GHW21} mentioned in Section~\ref{subsubsec:qcc-applications}  has already been established in \cite{AFH23b} for homogeneous graphs.
\end{OP} 

As was mentioned already in Section~\ref{subsec: buildings}, we are lacking an analog of Weyl chamber flows for affine buildings.

\begin{OP}[Weyl chamber flows for affine buildings] Establish a theory of Weyl chamber flows for affine buildings.
\begin{enumerate}
\item[(i)] Use these Weyl chamber flows to prove quantum-classical correspondences for compact quotients of affine buildings in the spirit of \cite{HWW21}.
\item[(ii)] Define multivariate dynamical zeta functions for these Weyl chamber flows and, for Bruhat-Tits buildings of $p$-adic groups of $\mathbb Q$-rank $1$, relate them to the geometric multi-variate zeta functions of Deitmar and Kang, \cite{DK14}.
\end{enumerate} 
\end{OP}   

Of course, at this point one may repeat all the questions and problems listed for the real case also for the $p$-adic situation. In fact, one could even reformulate many of them in an adelic context. We refrain from doing this but emphasize the point that the reduced analytic complexity of the non-Archimedean context may allow us to find solutions for the conceptual problems around higher rank resolvents that can then be transferred  to the Archimedean context.

%%%%%%%%%%%%%%%%%%%%%%%% referenc.tex %%%%%%%%%%%%%%%%%%%%%%%%%%%%%%
% sample references
% %
% Use this file as a template for your own input.
%
%%%%%%%%%%%%%%%%%%%%%%%% Springer-Verlag %%%%%%%%%%%%%%%%%%%%%%%%%%
%
% BibTeX users please use
% \bibliographystyle{}
% \bibliography{}
%

\end{document}